\documentclass{amsart}
\usepackage{amssymb}
\usepackage{amscd}
\usepackage{verbatim}
\usepackage{epsfig}
\begin{document}
\newcommand\rbm[1]{\mar{RBM:#1}}
\newcommand{\mar}[1]{{\marginpar{\sffamily{\scriptsize #1}}}}
\setlength{\marginparwidth}{1.12in}
\newcommand\av[1]{\mar{AV:#1}}

\newcommand\Mand{\ \text{and}\ }
\newcommand\Mwith{\ \text{with}\ }
\newcommand\Mfor{\ \text{for}\ }
\newcommand\Mst{\ \text{such that}\ }
\newcommand\Mor{\ \text{or}\ }
\newcommand\Mon{\ \text{on}\ }
\newcommand\Min{\ \text{in}\ }
\newcommand\Mif{\ \text{if}\ }
\newcommand\Miff{\ \text{iff}\ }
\newcommand\Mthen{\ \text{then}\ }
\newcommand\nin{\notin}
\newcommand\identity{\operatorname{id}}
\newcommand\Id{\operatorname{Id}}
\newcommand\Real{\mathbb{R}}
\newcommand\RR{\mathbb{R}}
\newcommand\pos{\Real^+}
\newcommand\Rnp{\Real\setminus\{0\}}
\newcommand\nzero{\setminus\{0\}}
\newcommand\Cx{\mathbb{C}}
\newcommand\Cxp{\Cx^+}
\newcommand\Cxm{\Cx^-}
\newcommand\Nat{\mathbb{N}}
\newcommand\halfNat{{\frac{1}{2}}\mathbb{N}}
\newcommand\intgr{\mathbb{Z}}
\newcommand\HH{\mathbb{H}}
\newcommand\im{\operatorname{Im}}
\newcommand\re{\operatorname{Re}}
\newcommand\sign{\operatorname{sign}}
\newcommand\codim{\operatorname{codim}}
\newcommand\End{\operatorname{End}}
\newcommand\Ker{\operatorname{Ker}}
\newcommand\Hom{\operatorname{Hom}}
\newcommand\tr{\operatorname{tr}}
\newcommand\Tr{\operatorname{Tr}}
\newcommand\ideal{{\mathcal I}}
\newcommand\Span{\operatorname{span}}
\newcommand\image{\operatorname{image}}
\newcommand\Range{\operatorname{Ran}}
\newcommand\Graph{\operatorname{graph}}
\newcommand\Card{\operatorname{Card}}
\newcommand\Hess{\operatorname{Hess}}
\newcommand\slim{\operatornamewithlimits{s-lim}}
\newcommand\spp{\operatorname{sp}}
\newcommand\sll{\operatorname{sl}}
\newcommand\sol{\operatorname{so}}
\newcommand\SL{\operatorname{SL}}
\newcommand\SO{\operatorname{SO}}
\newcommand\On{\operatorname{O}}
\newcommand\pa{\partial}
\newcommand\eff{\mathrm{eff}}
\newcommand\Rn{\Real^n}
\newcommand\Rm{\Real^m}
\newcommand\RN{\Real^N}
\newcommand\RtN{\Real^{2N}}
\newcommand\RM{\Real^M}
\newcommand\sphere{\mathbb{S}}
\newcommand\Sn{\sphere^{n-1}}
\newcommand\Sm{\sphere^{m-1}}
\newcommand\Snp{\sphere^n_+}
\newcommand\Smp{\sphere^m_+}
\newcommand\SN{\sphere^{N-1}}
\newcommand\SNp{\sphere^N_+}
\newcommand\circlep{\sphere^1_+}
\newcommand\Phom{P_{h}}
\newcommand\Shom{S_{h}}
\newcommand\distance{\operatorname{dist}}
\newcommand\cl{\operatorname{cl}}
\newcommand\interior{\operatorname{int}}
\newcommand\Fa{\operatorname{Fa}}
\newcommand\ff{\operatorname{ff}}
\newcommand\mf{\operatorname{mf}}
\newcommand\cf{\operatorname{cf}}
\newcommand\scf{\operatorname{sf}}
\newcommand\lf{\operatorname{lf}}
\newcommand\rf{\operatorname{rf}}
\newcommand\indfam{{\mathcal K}}
\newcommand\fraka{{\mathfrak a}}
\newcommand\cB{{\mathcal B}}
\newcommand\calA{{\mathcal A}}
\newcommand\calB{{\mathcal B}}
\newcommand\calR{{\mathcal R}}
\newcommand\calO{{\mathcal O}}
\newcommand\calJ{{\mathcal J}}
\newcommand\calL{{\mathcal L}}
\newcommand\calM{{\mathcal M}}
\newcommand\calN{{\mathcal N}}
\newcommand\calX{{\mathcal X}}
\newcommand\calU{{\mathcal U}}
\newcommand\calV{{\mathcal V}}
\newcommand\calF{{\mathcal F}}
\newcommand\calG{{\mathcal G}}
\newcommand\calT{{\mathcal T}}
\newcommand\calC{{\mathcal C}}
\newcommand\calCt{{\tilde {\mathcal C}}}
\newcommand\calCL{{\mathcal C}_{\text L}}
\newcommand\calCR{{\mathcal C}_{\text R}}
\newcommand\cF{{\mathcal F}}
\newcommand\cE{{\mathcal E}}
\newcommand\cH{{\mathcal H}}
\newcommand\cG{{\mathcal G}}
\newcommand\cU{{\mathcal U}}
\newcommand\Cinf{{\mathcal C}^{\infty}}
\newcommand\dist{{\mathcal C}^{-\infty}}
\newcommand\dCinf{{\dot{\mathcal C}}^{\infty}}
\newcommand\ddist{\dot\dist}
\newcommand\Cj{{\mathcal C}^j}
\newcommand\Linf{L^{\infty}}
\newcommand\phg{{\text{phg}}}
\newcommand\comp{{\text{comp}}}
\newcommand\loc{{\text{loc}}}
\newcommand\bcon{{\mathcal A}}
\newcommand\bconc{{\mathcal A}_{\text{phg}}}
\newcommand\Sch{{\mathcal S}}
\newcommand\temp{\Sch^{\prime}}
\newcommand\Diff{\operatorname{Diff}}
\newcommand\Diffb{\operatorname{Diff}_{\text{b}}}
\newcommand\Diffc{\operatorname{Diff}_{\text{c}}}
\newcommand\Diffsc{\operatorname{Diff}_{\text{sc}}}
\newcommand\DiffI{\operatorname{Diff}_{\text{I}}}
\newcommand\DiffIq{\operatorname{Diff}_{\text{I},q}}
\newcommand\sing{\text{sing}}
\newcommand\reg{\text{reg}}
\newcommand\supp{\operatorname{supp}}
\newcommand\ssupp{\operatorname{sing\ supp}}
\newcommand\csupp{\operatorname{cone\ supp}}
\newcommand\esupp{\operatorname{ess\ supp}}
\newcommand\Fr{{\mathcal F}}
\newcommand\Frinv{\Fr^{-1}}
\newcommand\bop{{\mathcal B}}
\newcommand\spec{\operatorname{spec}}
\newcommand\pspec{\spec_{pp}}
\newcommand\cspec{\spec_{c}}
\newcommand\FIO{{\mathcal I}}
\newcommand\SP{\operatorname{RC}}
\newcommand\RC{\operatorname{RC}}
\newcommand\Symc{S_c}
\newcommand\Symca{S_c^{\alpha}}
\newcommand\Symczero{S_c^{0,...,0}}
\newcommand\sci{{}^{\text{sc}}}
\newcommand\sct{\sci T^*}
\newcommand\scT{\sci T}
\newcommand\scdt{\sci \dot T^*}
\newcommand\dS{\dot S^*}
\newcommand\dT{\dot T^*}
\newcommand\dSreg{\dot\Sigma_{\text reg}}
\newcommand\scct{\sci\bar{T}^*}
\newcommand\Csc{C_{\text{sc}}}
\newcommand\SNpscd{(\SNp)^2_{\text{sc}}}
\newcommand\scdiag{\Delta_{\text{sc}}}
\newcommand\projscl{\pi^L_{\text{sc}}}
\newcommand\projscr{\pi^R_{\text{sc}}}
\newcommand\scHL{\sci H^{2,0}_{|\zeta|^2-\lambda^2}}
\newcommand\scHrg{\sci H^{2,0}_{\sqrt{g}}}
\newcommand\Hsc{H_{\text{sc}}}
\newcommand\Char{\operatorname{Char}}
\newcommand\dChar{\operatorname{\dot Char}}
\newcommand\WF{\operatorname{WF}}
\newcommand\WFb{\operatorname{WF}_{\bl}}
\newcommand\WFbz{\operatorname{WF}_{\bl}}
\newcommand\WFbd{\operatorname{WF}_{\bl}}
\newcommand\WFp{\operatorname{WF^{\prime}}}
\newcommand\WFsc{\operatorname{WF}_{\text{sc}}}
\newcommand\WFscp{\operatorname{WF_{sc}^{\prime}}}
\newcommand\WFC{\operatorname{WF}_C}
\newcommand\WFCi{\operatorname{WF}_{C_i}}
\newcommand\elliptic{\operatorname{ell}}
\newcommand\Psop{\operatorname{\Psi}}
\newcommand\Psib{\operatorname{\Psi}_{\bl}}
\newcommand\Psibc{\operatorname{\Psi}_{\text{bc}}}
\newcommand\Psiscrs{\operatorname{\Psi_{sc}^{-2,\infty}}}
\newcommand\Psiscr{\operatorname{\Psi_{sc}^{-2,0}}}
\newcommand\Psiscrm{\operatorname{\Psi_{sc}^{0,2}}}
\newcommand\PsiscHam{\operatorname{\Psi_{sc}^{2,0}}}
\newcommand\Psisci{\operatorname{\Psi_{sc}^{*,*}}}
\newcommand\Psiscid{\operatorname{\Psi_{sc}^{0,0}}}
\newcommand\Psiscis{\operatorname{\Psi_{sc}^{0,\infty}}}
\newcommand\Psiscsi{\operatorname{\Psi_{sc}^{-\infty,0}}}
\newcommand\Psiscs{\operatorname{\Psi_{sc}^{-\infty,\infty}}}
\newcommand\Psiscalg{\operatorname{\Psi_{sc}^{\infty,-\infty}}}
\newcommand\nullHam{{\mathcal N}}
\newcommand\charD{\Sigma_{\Delta-\lambda^2}}
\newcommand\charLap{\Sigma_{\Delta-\lambda}}
\newcommand\Snl{\Sn_{\lambda}}
\newcommand\SNl{\SN_{\lambda}}
\newcommand\gammat{\tilde\gamma}
\newcommand\gammasc{\gamma}
\newcommand\Tau{\mathcal{T}}
\newcommand\taut{\tilde\tau}
\newcommand\taub{\bar\tau}
\newcommand\Nout{N^+_{\lambda}}
\newcommand\Nin{N^-_{\lambda}}
\newcommand\Nio{N^{\pm}_{\lambda}}
\newcommand\El{E_{\lambda}}
\newcommand\Elt{\tilde E_{\lambda}}
\newcommand\Eil{E^i_{\lambda}}
\newcommand\Ejl{E^j_{\lambda}}
\newcommand\Eajl{E^{\alpha_j}_{\lambda}}
\newcommand\Eilt{\tilde E^i_{\lambda}}
\newcommand\Np{N^+}
\newcommand\Nm{N^-}
\newcommand\Npm{N^{\pm}}
\newcommand\Fin{F^-(\lambda)}
\newcommand\Fini{F^-_i(\lambda)}
\newcommand\Fout{F^+(\lambda)}
\newcommand\Fouti{F^+_i(\lambda)}
\newcommand\Foutj{F^+_j(\lambda)}
\newcommand\Rout{R^+_{\lambda}}
\newcommand\Routl{R^+_{\lambda^2}}
\newcommand\Routsgnl{R^{\sign\lambda}_{\lambda^2}}
\newcommand\Rin{R^-_{\lambda}}
\newcommand\Rinl{R^-_{\lambda^2}}
\newcommand\Rinsgnl{R^{-\sign\lambda}_{\lambda^2}}
\newcommand\Rio{R^{\pm}_{\lambda}}
\newcommand\Riol{R^{\pm}_{\lambda^2}}
\newcommand\Roi{R^{\mp}_{\lambda}}
\newcommand\Roil{R^{\mp}_{\lambda^2}}
\newcommand\Riob{R^{\pm}}
\newcommand\Roib{R^{\mp}}
\newcommand\Tio{T^{\pm}}
\newcommand\Tiob{T^{\pm}_{\ff}}
\newcommand\Toi{T^{\mp}}
\newcommand\Toib{T^{\mp}_{\ff}}
\newcommand\TIiob{T_I^{\pm}}
\newcommand\Rinb{R^-}
\newcommand\Rinbsgnl{R^{-\sign\lambda}}
\newcommand\Tin{T^-}
\newcommand\Tinb{T^-_{\ff}}
\newcommand\TIinb{T^-_I}
\newcommand\Routb{R^+}
\newcommand\Routbsgnl{R^{\sign\lambda}}
\newcommand\Tout{T^+}
\newcommand\Toutb{T^+_{\ff}}
\newcommand\TIoutb{T^+_I}
\newcommand\Rlkf{(|\xib|^2-(\lambda-i0)^2)^{-1}}
\newcommand\Rlk{\rho_0(\lambda)}
\newcommand\Rmlk{\rho_0(-\lambda)}
\newcommand\Rpmlk{\rho_0(\pm\lambda)}
\newcommand\Rlka{\rho_1(\lambda)}
\newcommand\Rlkb{\rho_2(\lambda)}
\newcommand\Rilk{\rho_i(\lambda)}
\newcommand\reduced{\natural}
\newcommand\Rlf{R_0(\lambda)}
\newcommand\Rla{R_1(\lambda)}
\newcommand\Rlb{R_2(\lambda)}
\newcommand\Ril{R_i(\lambda)}
\newcommand\Rlj{R_j(\lambda)}
\newcommand\Rlft{R_0(\lambda)}
\newcommand\Rflambda{R_0^{\reduced}(\sigma)}
\newcommand\RV{R^{\reduced}_V}
\newcommand\Rfsigma{R_0^{\reduced}(\sigma)}
\newcommand\Rfsigmah{R_0^{\reduced}(\sigma^{1/2})}
\newcommand\Rfzero{R_0^{\reduced}(0)}
\newcommand\RlV{R^{\reduced}_V(\sigma)}
\newcommand\RlVi{R^{\reduced}_{V_i}(\sigma)}
\newcommand\RlVt{R_V(\lambda)}
\newcommand\RlVtL{{R}_V^L(\lambda)}
\newcommand\RlVtR{{R}_V^R(\lambda)}
\newcommand\RlVit{{R}_{V_i}(\lambda)}
\newcommand\RlVta{{R}_V^{(1)}(\lambda)}
\newcommand\RlVtk{{R}_V^{(k)}(\lambda)}
\newcommand\RlVatV{{R}_{V_{\alpha}}(\lambda)V_{\alpha}}
\newcommand\RlVatVa{{R}_{V_{\alpha_1}}(\lambda)V_{\alpha_1}}
\newcommand\RlVatVb{{R}_{V_{\alpha_2}}(\lambda)V_{\alpha_2}}
\newcommand\RlVatVk{{R}_{V_{\alpha_k}}(\lambda)V_{\alpha_k}}
\newcommand\RlVatVkk{{R}_{V_{\alpha_{k+1}}}(\lambda)V_{\alpha_{k+1}}}
\newcommand\RlVaptV{{R}_{V_{\alpha'}}(\lambda)V_{\alpha'}}
\newcommand\RlVapptV{{R}_{V_{\alpha''}}(\lambda)V_{\alpha''}}
\newcommand\RlVajtV{{R}_{V_{\alpha_j}}(\lambda)V_{\alpha_j}}
\newcommand\RlVaktV{{R}_{V_{\alpha_k}}(\lambda)V_{\alpha_k}}
\newcommand\RlVakktV{{R}_{V_{\alpha_{k+1}}}(\lambda)V_{\alpha_{k+1}}}
\newcommand\Tl{T(\lambda)}
\newcommand\Tlt{\tilde\Tl}
\newcommand\Tltp{\tilde T'(\lambda)}
\newcommand\Tltpp{\tilde T''(\lambda)}
\newcommand\Tli{T_i(\lambda)}
\newcommand\Tlit{\tilde\Tli}
\newcommand\Tlip{T_i'(\lambda)}
\newcommand\Tlipp{T_i''(\lambda)}
\newcommand\Tlj{T_j(\lambda)}
\newcommand\Tla{T_{\alpha}(\lambda)}
\newcommand\Tlaa{T_{\alpha_1}(\lambda)}
\newcommand\Tlab{T_{\alpha_2}(\lambda)}
\newcommand\Tlak{T_{\alpha_k}(\lambda)}
\newcommand\Tlakt{\tilde\Tlak}
\newcommand\Tlaj{T_{\alpha_j}(\lambda)}
\newcommand\Tlajj{T_{\alpha_{j+1}}(\lambda)}
\newcommand\Tlajp{T_{\alpha_j}'(\lambda)}
\newcommand\Tlajpt{\tilde\Tlajp}
\newcommand\Tlajt{\tilde\Tlaj}
\newcommand\Tlakk{T_{\alpha_{k+1}}(\lambda)}
\newcommand\Tlakkp{T_{\alpha_{k+1}}'(\lambda)}
\newcommand\Tlap{T_{\alpha'}(\lambda)}
\newcommand\Tlapt{\tilde\Tlap}
\newcommand\Tlapp{T_{\alpha''}(\lambda)}
\newcommand\Tkl{T^{(k)}(\lambda)}
\newcommand\Tcl{T^{\flat}(\lambda)}
\newcommand\Fl{F(\lambda)}
\newcommand\BlVt{\tilde B_V(\lambda)}
\newcommand\KBlVt{K_{\BlVt}}
\newcommand\BlVaat{B_{V_{\alpha_1}}(\lambda)}
\newcommand\BV{B_V}
\newcommand\Bone{B_1}
\newcommand\Btwo{B_2}
\newcommand\Bthree{B_3}
\newcommand\Banyj{B_j}
\newcommand\PlV{P_V(\lambda)}
\newcommand\PlVc{P_V^{\flat}(\lambda)}
\newcommand\Pl{P_0(\lambda)}
\newcommand\SVl{S_V(\lambda)}
\newcommand\Sjr{S_j^{\reduced}}
\newcommand\Rkp{{\mathcal R}^k_+}
\newcommand\Rkm{{\mathcal R}^k_-}
\newcommand\Rkpm{{\mathcal R}^k_{\pm}}
\newcommand\Phys{{\mathcal P}}
\newcommand\Pc{\overline{\mathcal P}}
\newcommand\pip{\pi^{\perp}}
\newcommand\pipa{\pi_\partial}
\newcommand\gammapa{\gamma_\partial}
\newcommand\pipah{\hat\pi_\partial}
\newcommand\pit{\tilde\pi}
\newcommand\xit{\tilde\xi}
\newcommand\zetat{\tilde\zeta}
\newcommand\etat{\tilde\eta}
\newcommand\sigmat{\tilde\sigma}
\newcommand\sigmahat{\hat\sigma}
\newcommand\thetat{\tilde\theta}
\newcommand\psit{\tilde\psi}
\newcommand\phit{\tilde\phi}
\newcommand\chit{\tilde\chi}
\newcommand\rhot{\tilde\rho}
\newcommand\xib{\bar\xi}
\newcommand\zetab{\bar\zeta}
\newcommand\thetab{\bar\theta}
\newcommand\etab{\bar\eta}
\newcommand\iotal{\iota_{\lambda}}
\newcommand\rhoat{\rhot_{\alpha_1}}
\newcommand\Lambdat{\tilde\Lambda}
\newcommand\Lambdati{\tilde\Lambda^{\text{in}}}
\newcommand\Lambdato{\tilde\Lambda^{\text{out}}}
\newcommand\Lambdatp{\tilde\Lambda^{\text{prop}}}
\newcommand\Lambdai{\Lambda^{\text{in}}}
\newcommand\Lambdao{\Lambda^{\text{out}}}
\newcommand\poles{\Lambda'}
\newcommand\rpoles{\Lambda_p}
\newcommand\thresholds{\Lambda}
\newcommand\Vt{\tilde V}
\newcommand\It{\tilde I}
\newcommand\half{{\frac{1}{2}}}
\newcommand\sigmah{\sigma^{1/2}}
\newcommand\bX{\partial X}
\newcommand\bXb{\partial \Xb}
\newcommand\Deltabt{\tilde\Delta_0}
\newcommand\strip{\Omega_T}
\newcommand\Kf{K^{\flat}}
\newcommand\Gs{G^{\sharp}}
\newcommand\Gt{\tilde G}
\newcommand\Osc{\sci\Omega}
\newcommand\OSc{{}^\Scl\Omega}
\newcommand\Osch{\sci\Omega^{\half}}
\newcommand\Oscmh{\sci\Omega^{-\half}}
\newcommand\Isc{I_{sc}}
\newcommand\os{{\text{os}}}
\newcommand\Qzl{Q^0_{-\lambda}}
\newcommand\Lie{{\mathcal L}}
\newcommand\bl{{\text b}}
\newcommand\scl{{\text{sc}}}
\newcommand\sccl{{\text{scc}}}
\newcommand\Scl{{\text{sc}}}
\newcommand\ScLl{{\text{Sc,L}}}
\newcommand\ScRl{{\text{Sc,R}}}
\newcommand\Sccl{{\text{Scc}}}
\newcommand\sus{{\text{sus}}}
\newcommand\ssl{{\text{ss}}}
\newcommand\XXb{X^2_\bl}
\newcommand\XXbt{\Xt^2_\bl}
\newcommand\XXsc{X^2_\scl}
\newcommand\XXsct{\Xt^2_\scl}
\newcommand\XXSc{X^2_\Scl}
\newcommand\XXSct{\Xt^2_\Scl}
\newcommand\XXScL{X^2_\ScLl}
\newcommand\XXScR{X^2_\ScRl}
\newcommand\MMsc{M^2_\scl}
\newcommand\Deltab{\Delta_\bl}
\newcommand\Deltasc{\Delta_\scl}
\newcommand\DeltaSc{\Delta_\Scl}
\newcommand\DeltaScL{\Delta_\ScLl}
\newcommand\DeltaScR{\Delta_\ScRl}
\newcommand\prs{\sigma}
\newcommand\Nsc{N_\scl}
\newcommand\Nscp{N_{\scl,p}}
\newcommand\Nff{N_{\ff}}
\newcommand\Nffz{N_{\ff,0}}
\newcommand\Nffzp{N_{\ff,0,p}}
\newcommand\Nffl{N_{\ff,l}}
\newcommand\Nffml{N_{\ff,-l}}
\newcommand\Nmf{N_{\mf}}
\newcommand\Nmfz{N_{\mf,0}}
\newcommand\Nmfl{N_{\mf,l}}
\newcommand\Nmfml{N_{\mf,-l}}
\newcommand\ffb{\operatorname{bf}}
\newcommand\Ffb{\operatorname{bf'}}
\newcommand\ffsc{\operatorname{sf}}
\newcommand\ffSc{\operatorname{sf_C}}
\newcommand\Ffsc{\operatorname{sf'}}
\newcommand\rff{\rho_{\ff}}
\newcommand\rmf{\rho_{\mf}}
\newcommand\rffb{\rho_{\ffb}}
\newcommand\rffsc{\rho_{\ffsc}}
\newcommand\rFfsc{\rho_{\Ffsc}}
\newcommand\rffSc{\rho_{\ffSc}}
\newcommand\rinf{\rho_{\infty}}
\newcommand\CL{C_L}
\newcommand\CR{C_R}
\newcommand\betab{\beta_\bl}
\newcommand\betasc{\beta_\scl}
\newcommand\betaSc{\beta_\Scl}
\newcommand\BetaSc{\bar\beta_\Scl}
\newcommand\betaScL{\beta_\ScLl}
\newcommand\betaScR{\beta_\ScRl}
\newcommand\ScT{{}^\Scl T^*}
\newcommand\SccT{{}^\Scl \bar T^*}
\newcommand\ScS{{}^\Scl S^*}
\newcommand\Sb{{}^\bl S}
\newcommand\Tb{{}^\bl T}
\newcommand\Tsc{{}^\scl T}
\newcommand\TSc{{}^\Scl T}
\newcommand\CSc{C_\Scl}
\newcommand\Lambdasc{{}^\scl\Lambda}
\newcommand\XXXb{X^3_\bl}
\newcommand\XXXsc{X^3_\scl}
\newcommand\XXXSc{X^3_\Scl}
\newcommand\XXXScO{X^3_{\Scl,O}}
\newcommand\XXXScF{X^3_{\Scl,F}}
\newcommand\XXXScS{X^3_{\Scl,S}}
\newcommand\XXXScC{X^3_{\Scl,C}}
\newcommand\KDsc{\operatorname{KD^{\half}_\scl}}
\newcommand\KDSc{\operatorname{KD^{\half}_\Scl}}
\newcommand\KDScEF{\operatorname{KD^{E,F}_\Scl}}
\newcommand\Oh{\operatorname{\Omega^{\half}}}
\newcommand\WFSc{\WF_\Scl}
\newcommand\WFtSc{\WF_{\text 3sc}}
\newcommand\WFScmf{\WF_{\Scl,\mf}}
\newcommand\WFScff{\WF_{\Scl,\ff}}
\newcommand\WFScs{\WF_{\Scl,\prs}}
\newcommand\WFScp{\WF'_\Scl}
\newcommand\WFScmfp{\WF'_{\Scl,\mf}}
\newcommand\WFScffp{\WF'_{\Scl,\ff}}
\newcommand\WFScsp{\WF'_{\Scl,\prs}}
\newcommand\Diffscc{\Diff_\sccl}
\newcommand\DiffSc{\Diff_\Scl}
\newcommand\DiffScc{\Diff_\Sccl}
\newcommand\DiffscI{\Diff_{\scl,\text{I}}}
\newcommand\VscI{\Vf_{\scl,\text{I}}}
\newcommand\DiffsV{\operatorname{Diff}_{\sus(V)}}
\newcommand\DiffsVsc{\operatorname{Diff}_{\sus(V),\scl}}
\newcommand\DiffsVCsc{\operatorname{Diff}_{\sus(V)-C,\scl}}   
\newcommand\Psisc{\Psop_\scl}
\newcommand\Psiscc{\Psop_\sccl}
\newcommand\Psiss{\Psop_\ssl}
\newcommand\Psisch{\Psop_{\scl,h}}
\newcommand\Psiscch{\Psop_{\sccl,h}}
\newcommand\PsiSc{\Psop_\Scl}
\newcommand\PsiScph{\Psop_{\Scl,\phi}}
\newcommand\PsiScra{\Psop_{\Scl,\rho^\sharp_a}}
\newcommand\PsiScc{\Psop_\Sccl}
\newcommand\PsiSccml{\Psop^{m,l}_\Sccl}
\newcommand\PsiScxx{\Psop^{*,*}_\Scl}
\newcommand\PsiScml{\Psop^{m,l}_\Scl}
\newcommand\PsiScmz{\Psop^{m,0}_\Scl}
\newcommand\PsiScmmz{\Psop^{-m,0}_\Scl}
\newcommand\PsiSckz{\Psop^{k,0}_\Scl}
\newcommand\PsiScmmml{\Psop^{-m,-l}_\Scl}
\newcommand\Psiscmkk{\Psop^{-k,k}_\scl}
\newcommand\Psiscmmmkk{\Psop^{-m-k,k}_\scl}
\newcommand\Psiscmoo{\Psop^{-1,1}_\scl}
\newcommand\Psiscmz{\Psop^{m,0}_\scl}
\newcommand\Psiscmmz{\Psop^{-m,0}_\scl}
\newcommand\PsiSckmkl{\Psop^{km,kl}_\Scl}
\newcommand\PsiScmplp{\Psop^{m',l'}_\Scl}
\newcommand\PsiScmmpllp{\Psop^{m+m',l+l'}_\Scl}
\newcommand\Psiscml{\Psop^{m,l}_\scl}
\newcommand\PsiScid{\Psop^{0,0}_\Scl}
\newcommand\PsiSczo{\Psop^{0,1}_\Scl}
\newcommand\PsiScmii{\Psop^{-\infty,\infty}_\Scl}
\newcommand\PsiScmiz{\Psop^{-\infty,0}_\Scl}
\newcommand\PsiScmoo{\Psop^{-1,1}_\Scl}
\newcommand\PsisCid{\Psop^{0,0}_{\scl-C}}
\newcommand\PsisC{\Psop_{\scl-C}}
\newcommand\Psiinf{\Psop_{\infty}}
\newcommand\Psiinfid{\Psop_{\infty}^0}
\newcommand\PsiFinf{\Psop_{\infty-\Fr}}
\newcommand\PsisVscml{\Psop^{m,l}_{\sus(V),\scl}}
\newcommand\PsisVsc{\Psop_{\sus(V),\scl}}
\newcommand\PsisVpsc{\Psop_{\sus(V_p),\scl}}
\newcommand\PsisVCSc{\Psop_{\sus(V)-C,\scl}}
\newcommand\SFinf{S_{\infty-\Fr}}
\newcommand\YsVC{Y^2_{\sus(V)-C,\scl}}
\newcommand\ffYsc{\ffsc_{\sus(V)}}
\newcommand\SXC{S(X;C)}
\newcommand\Ios{I_{\text{os}}}
\newcommand\pbL{\pi^2_{\bl,{\text L}}}
\newcommand\pbR{\pi^2_{\bl,{\text R}}}
\newcommand\pscL{\pi^2_{\scl,{\text L}}}
\newcommand\pscR{\pi^2_{\scl,{\text R}}}
\newcommand\PbO{\pi^3_{\bl,{\text O}}}
\newcommand\PscO{\pi^3_{\scl,{\text O}}}
\newcommand\PScO{\pi^3_{\Scl,{\text O}}}
\newcommand\PScF{\pi^3_{\Scl,{\text F}}}
\newcommand\PScC{\pi^3_{\Scl,{\text C}}}
\newcommand\PScS{\pi^3_{\Scl,{\text S}}}
\newcommand\pScL{\pi^2_{\Scl,{\text L}}}
\newcommand\pScR{\pi^2_{\Scl,{\text R}}}
\newcommand\CLF{\CL^F}
\newcommand\CLO{\CL^O}
\newcommand\CLS{\CL^S}
\newcommand\CLC{\CL^C}
\newcommand\DeltaYb{\Delta_{\bl,Y}}
\newcommand\DeltaYsc{\Delta_{\sus-\scl}}
\newcommand\diag{\operatorname{diag}}
\newcommand\Vf{{\mathcal V}}
\newcommand\Vb{{\mathcal V}_{\bl}}
\newcommand\Vsc{{\mathcal V}_{\scl}}
\newcommand\VSc{{\mathcal V}_{\Scl}}
\newcommand\VfI{\Vf_{\text{I}}}
\newcommand\VfIq{\Vf_{\text{I},q}}
\newcommand\scH{{}^\scl H}
\newcommand\scHg{\scH_g}
\newcommand\Hss{H_\ssl}
\newcommand\xh{\hat x}
\newcommand\yh{\hat y}
\newcommand\sh{\hat s}
\newcommand\rh{\hat r}
\newcommand\Yh{\hat Y}
\newcommand\Zh{\hat Z}
\newcommand\Yb{\bar Y}
\newcommand\hb{\bar h}
\newcommand\xih{\hat\xi}
\newcommand\etah{\hat\eta}
\newcommand\muh{\hat\mu}
\newcommand\mub{\bar\mu}
\newcommand\nub{\bar\nu}
\newcommand\mubh{\widehat{\bar\mu}}
\newcommand\yb{\bar y}
\newcommand\zb{\bar z}
\newcommand\ub{\bar u}
\newcommand\Qb{\bar Q}
\newcommand\Wbp{{\bar W}^\perp}
\newcommand\Wp{W^\perp}
\newcommand\Kt{\tilde K}
\newcommand\Wt{\tilde W}
\newcommand\Ut{\tilde U}
\newcommand\yt{\tilde y}
\newcommand\ut{\tilde u}
\newcommand\vt{\tilde v}
\newcommand\ft{\tilde f}
\newcommand\htil{\tilde h}
\newcommand\St{\tilde S}
\newcommand\Pt{\tilde P}
\newcommand\Rt{\tilde R}
\newcommand\qt{\tilde q}
\newcommand\Qt{\tilde Q}
\newcommand\Xb{\bar X}
\newcommand\lambdat{\tilde\lambda}
\newcommand\betat{\tilde\beta}
\newcommand\Phit{\tilde\Phi}
\newcommand\epst{\tilde\epsilon}
\newcommand\ep{\epsilon}
\newcommand\bt{\tilde b}
\newcommand\Xt{\tilde X}
\newcommand\Mt{\tilde M}
\newcommand\At{\tilde A}
\newcommand\Et{\tilde E}
\newcommand\Ht{\tilde H}
\newcommand\at{\tilde a}
\newcommand\Ct{\tilde C}
\newcommand\pih{\hat\pi}
\newcommand\Rh{\hat R}
\newcommand\Ah{\hat A}
\newcommand\Bh{\hat B}
\newcommand\Ch{\hat C}
\newcommand\Gh{\hat G}
\newcommand\Hh{\hat H}
\newcommand\Qh{\hat Q}
\newcommand\Ph{\hat P}
\newcommand\Nh{\hat N}
\newcommand\Sh{\hat S}
\newcommand\Gcal{{\mathcal G}}
\newcommand\GcalC{{\mathcal G}_C}
\newcommand\Jcal{{\mathcal J}}
\newcommand\JcalC{{\mathcal J}_C}
\setcounter{secnumdepth}{3}
\newtheorem{lemma}{Lemma}[section]
\newtheorem{prop}[lemma]{Proposition}
\newtheorem{thm}[lemma]{Theorem}
\newtheorem{cor}[lemma]{Corollary}
\newtheorem{result}[lemma]{Result}
\newtheorem*{thm*}{Theorem}
\newtheorem*{prop*}{Proposition}
\newtheorem*{cor*}{Corollary}
\newtheorem*{conj*}{Conjecture}
\numberwithin{equation}{section}
\theoremstyle{remark}
\newtheorem{rem}[lemma]{Remark}
\newtheorem*{rem*}{Remark}
\theoremstyle{definition}
\newtheorem{Def}[lemma]{Definition}
\newtheorem*{Def*}{Definition}
\def\signature#1#2{\par\noindent#1\dotfill\null\\*
{\raggedleft #2\par}}

\renewcommand{\theenumi}{\roman{enumi}}
\renewcommand{\labelenumi}{(\theenumi)}

\title[Propagation of singularities]
{Propagation of singularities for the wave equation on manifolds
with corners}
\author[Andras Vasy]{Andr\'as Vasy}
\address{Department of Mathematics, Massachusetts Institute of Technology,
Cambridge MA 02139, U.S.A.}
\email{andras@math.mit.edu}
\date{May 22, 2004}
\thanks{This work is partially supported by NSF grant \#DMS-0201092, and
a Fellowship from the Alfred P.\ Sloan Foundation.}

\begin{abstract}
In this paper we describe the propagation of $\Cinf$ and Sobolev
singularities for the wave equation on $\Cinf$ manifolds with corners $M$
equipped with a Riemannian metric $g$. That is, for $X=M\times\Real_t$,
$P=D_t^2-\Delta_M$, and $u\in H^1_{\loc}(X)$ solving
$Pu=0$ with homogeneous Dirichlet
or Neumann boundary conditions, we show that $\WFb(u)$ is a union
of maximally extended generalized broken bicharacteristics. This
result is a $\Cinf$ counterpart of Lebeau's results for the propagation
of analytic singularities on real analytic manifolds with
appropriately stratified boundary,
\cite{Lebeau:Propagation}. Our methods rely on b-microlocal
positive commutator estimates, thus providing a new proof for the
propagation of singularities at hyperbolic points even if $M$
has a smooth boundary (and no corners).
\end{abstract}

\maketitle

\section{Introduction}
In this paper
we describe the propagation of $\Cinf$ singularities for
the wave equation on a manifold with corners $M$ equipped with
a smooth Riemannian metric $g$. Let $\Delta=\Delta_g$ be the positive
Laplacian of $g$, let $X=M\times\Real_t$, $P=D_t^2-\Delta$,
and consider the Dirichlet
boundary condition for $P$:
\begin{equation*}
Pu=0,\ u|_{\pa X}=0,
\end{equation*}
with the boundary condition meaning more precisely that
$u\in H^1_{0,\loc}(X)$. Here $H^1_0(X)$ is the completion of $\dCinf_c(X)$
(the vector space of $\Cinf$ functions of compact support on $X$, vanishing
with all derivatives at $\pa X$) with
respect to $\|u\|^2_{H^1(X)}=\|du\|_{L^2(X)}+\|u\|_{L^2(X)}$,
$L^2(X)=L^2(X,dg\,dt)$, and
$H^1_{0,\loc}(X)$ is its localized version, i.e.\ $u\in H^1_0(X)$ if
for all $\phi\in\Cinf_c(X)$, $\phi u\in H^1_0(X)$. At the end of the
introduction we also consider Neumann boundary conditions.

The statement of the propagation of singularities of solutions
has two additional ingredients: locating singularities of a distribution,
as captured by the wave front set, and describing the curves along
which they propagate, namely the bicharacteristics. Both of these
are closely related to an appropropriate notion of phase space, in which
both the wave front set and the bicharacteristics are located. On manifolds
without boundary, this phase space is the standard cotangent bundle. In the
presence of boundaries the phase space is the b-cotangent bundle, $\Tb^*X$,
(`b' stands for boundary)
which we now briefly describe following \cite{Melrose:Atiyah},
which mostly deals with the $\Cinf$
boundary case, and especially
\cite{Melrose-Piazza:Analytic}.

Thus, $\Vb(X)$ is, by definition,
the Lie algebra of $\Cinf$ vector fields on $X$
tangent to every boundary face of $X$. Now, $\Vb(X)$ is the
set of all $\Cinf$ sections of a vector bundle $\Tb X$ over $X$.
The dual bundle of $\Tb X$ is $\Tb^*X$; this is the phase space in our
setting.
Let $o$ denote the zero section of $\Tb^*X$ (as well as other related
vector bundles below). Then $\Tb^*X\setminus o$ is equipped with
an $\Real^+$-action (fiberwise multiplication) which has no fixed
points. It is often natural to take the
quotient with the $\Real^+$-action, and work on the b-cosphere bundle,
$\Sb^*X$.

The differential operator algebra generated by $\Vb(X)$ is denoted
by $\Diffb(X)$, and its microlocalization is $\Psib(X)$, the
algebra of b-, or totally characteristic, pseudodifferential operators.
For $A\in\Psib^m(X)$,
$\sigma_{b,m}(A)$ is a homogeneous degree $m$
function on $\Tb^*X\setminus o$. Since $X$ is not compact, even if
$M$ is, we always understand that $\Psib^m(X)$ stands for properly
supported ps.d.o's, so its elements define continuous maps
$\dCinf(X)\to\dCinf(X)$ as well as $\dist(X)\to\dist(X)$.

We are now ready to define the wave front set $\WFb(u)$ for
$u\in H^1_{\loc}(X)$. This measures if $u$ has additional regularity,
locally in $\Tb^*X$, relative to $H^1$.
For $u\in H^1_{\loc}(X)$, $q\in \Tb^* X\setminus o$, $m\geq 0$,
we say that $q\nin\WFbz^{1,m}(u)$
if there is $A\in\Psib^{m}(X)$ such that $\sigma_{b,m}(A)(q)\neq 0$ and
$Au\in H^1(X)$. Since compactly supported elements of
$\Psib^0(X)$ preserve $H^1_{\loc}(X)$, it follows that for
$u\in H^1_{\loc}(X)$, $\WFbz^{1,0}(u)=\emptyset$.
For any $m$, $\WFbz^{1,m}(u)$ is a conic subset of $\Tb^*X\setminus o$;
hence it is natural to identify it with a subset of $\Sb^*X$. Its
intersection with $\Tb^*_{X^\circ}X\setminus o$, which can be naturally
identified with $T^*X^\circ\setminus o$, is $\WF^{m+1}(u)$. Thus, in the
interior of $X$, $\WFbz^{1,m}(u)$ measures if $u$ is microlocally in $H^{m+1}$.
The main result of this paper, stated at the end of this section,
is that for $u\in H^1_0(X)$ with $Pu=0$,
$\WFbz^{1,m}(u)$ is a union of maximally extended generalized
broken bicharacteristics, which are defined below.
In fact, the requirement $u\in H^1_0(X)$ can be relaxed and
$m$ can be allowed to be negative, see
Definitions~\ref{Def:H1m-neg}-\ref{Def:WFb-neg}.
We also remark that
for such $u$, the $H^1(X)$-based b-wave front set, $\WFbz^{1,m}(u)$,
could be replaced
by an $L^2(X)$-based b-wave front set, see Lemma~\ref{lemma:H10-L2}.
In addition, our methods apply, a fortiori, for elliptic problems
such as $\Delta_g$ on $(M,g)$,
e.g.\ showing that $u\in H^1_{0,\loc}(M)$ and $(\Delta_g-\lambda)u=0$
imply $u\in H^{1,\infty}_{\bl,\loc}(M)$, so $u$ is conormal -- see the
end of Section~\ref{sec:elliptic}.

This propagation result is the $\Cinf$ (and Sobolev space)
analogue of Lebeau's result
\cite{Lebeau:Propagation} for analytic singularities of $u$ when
$M$ and $g$ are real analytic. Thus,
the geometry is similar in the two settings, but the analytic techniques
are rather different: Lebeau uses complex scaling and the analytic
wave front set of the extension of $u$ as $0$ to a neighborhood of $X$
(in an extension $\tilde X$
of the manifold $X$), while we use positive commutator
estimates and b-microlocalization relative to the form domain of
the Laplacian. In fact, our microlocalization techniques, especially
the positive commutator constructions, are very closely related
to the methods used in $N$-body scattering, \cite{Vasy:Propagation-Many},
to prove the propagation of singularities (meaning microlocal lack of
decay at infinity) there.
Although Lebeau allows more general singularities than
corners for $X$, provided that $X$ sits in a real analytic manifold
$\tilde X$
with $g$ extending to $\tilde X$,
we expect to generalize our results to settings where no analogous
$\Cinf$ extension is available, see the remarks at the end of
the introduction.

We now describe the setup in more detail so that our main theorem can
be stated in a precise fashion.
Let $F_i$, $i\in I$, be the closed boundary faces of $M$ (including $M$),
$\cF_i=F_i\times \Real$, $\cF_{i,\reg}$ the interior (`regular part') of
$\cF_i$. Note that for each $p\in X$, there is a unique $i$ such that
$p\in \cF_{i,\reg}$.
Now, there is a natural non-injective `inclusion'
$\iota:T^* X\to\Tb^*X$, and the range
of $\iota$ over the interior of a face $\cF_i$ lies in $T^* \cF_i$ (which
is well-defined as a subspace of $\Tb^*X$),
while its kernel is $N^*\cF_i$, the conormal bundle of $\cF_i$ in $X$. Thus,
we define the compressed b-cotangent bundle $\dot\Tb^* X$
\begin{equation*}
\dot\Tb^* X=\cup_{i\in I}T^*\cF_{i,\reg}\subset\Tb^*X.
\end{equation*}
Regarded as a map $T^* X\to \dot\Tb^* X$ (i.e.\ onto its range) we relabel
$\iota$ as a projection $\pi$. We write $o$ for the `zero section'
of $\dot\Tb^*X$ as well, so
\begin{equation*}
\dot\Tb^* X\setminus o=\cup_{i\in I}T^*\cF_{i,\reg}\setminus o,
\end{equation*}
and then $\pi$ restricts to a map
\begin{equation*}
T^*X\setminus\cup_i N^*\cF_i\to \dot\Tb^* X\setminus o.
\end{equation*}

Now, the characteristic set $\Char(P)\subset T^*X\setminus o$
of $P$ is defined by
$p^{-1}(\{0\})$, where $p\in\Cinf(T^*X\setminus o)$
is the principal symbol of $P$, which is homogeneous degree $2$ on $T^*X
\setminus o$. Notice that $\Char(P)\cap N^*\cF_i=\emptyset$ for all $i$,
i.e.\ the boundary faces are all non-characteristic for $P$.
Thus, $\pi(\Char(P))\subset\dot\Tb^*X\setminus o$.
We define the elliptic, glancing and hyperbolic sets by
\begin{equation*}\begin{split}
&\cE=\{q\in \dot \Tb^* X\setminus o:\ \pi^{-1}(q)\cap\Char(P)=\emptyset\},\\
&\cG=\{q\in \dot \Tb^* X\setminus o:\ \Card(\pi^{-1}(q)\cap\Char(P))=1\},\\
&\cH=\{q\in \dot \Tb^* X\setminus o:\ \Card(\pi^{-1}(q)\cap\Char(P))\geq 2\},
\end{split}\end{equation*}
with $\Card$ denoting the cardinality of a set;
each of these is a conic subset of $\dot\Tb^*X\setminus o$. Note that
in $T^*X^\circ$, $\pi$ is the identity map, so every point $q\in T^*X^\circ$
is either in $\cE$ or $\cG$ depending on whether $q\nin\Char(P)$ or
$q\in\Char(P)$.

We briefly describe these sets in local coordinates.
Let $p\in \pa X$, and let $\cF_i$ be the closed face of $X$ with the
smallest dimension that contains $p$, so $p\in \cF_{i,\reg}$.
Local coordinates near $p$ are given by $(x_1,\ldots,x_k,y_1,\ldots,y_l,t)$
where $\cF_i$ is defined by $x_1=\ldots=x_k=0$, and the other boundary
faces through $p$ are given by the vanishing of a subset of
the collection $x_1,\ldots,x_k$ of functions -- in particular, the $k$ boundary
hypersurfaces are given by $x_j=0$ for $j=1,\ldots,k$.

Such local coordinates on the base induce 
local coordinates on the cotangent bundle, namely
$(x,y,t,\xi,\zeta,\tau)$ on $T^*X$ near $\pi^{-1}(q)$,
$q\in T^*\cF_{i,\reg}$, and
corresponding coordinates $(y,t,\zeta,\tau)$ on a neighborhood
$\cU$ of $q$ in $T^*\cF_{i,\reg}$. The metric
function on $T^*M$ has the form
\begin{equation*}
g(x,y,\xi,\zeta)=\sum_{i,j}A_{ij}(x,y)\xi_i\xi_j
+\sum_{i,j}2C_{ij}(x,y)\xi_i\zeta_j+\sum_{i,j}B_{ij}(x,y)\zeta_i\zeta_j
\end{equation*}
with $A,B,C$ smooth. Moreover, these coordinates can be chosen (i.e.\ the
$y_j$ can be adjusted) so
that $C(0,y)=0$.
Thus,
\begin{equation*}
p|_{x=0}=\tau^2-\xi\cdot A(y)\xi-\zeta \cdot B(y)\zeta,
\end{equation*}
with $A$, $B$ positive definite matrices depending smoothly on $y$, so
\begin{equation*}\begin{split}
&\cE\cap\cU=\{(y,t,\zeta,\tau):\ \tau^2<\zeta \cdot B(y)\zeta,\ (\zeta,\tau)\neq 0\},\\
&\cG\cap\cU=\{(y,t,\zeta,\tau):\ \tau^2=\zeta \cdot B(y)\zeta,\ (\zeta,\tau)\neq 0\},\\
&\cH\cap\cU=\{(y,t,\zeta,\tau):\ \tau^2>\zeta \cdot B(y)\zeta,\ (\zeta,\tau)\neq 0\}.
\end{split}\end{equation*}

The compressed characteristic set is
\begin{equation*}
\dot\Sigma=\pi(\Char(P))=\cG\cup\cH,
\end{equation*}
and
\begin{equation*}
\hat\pi:\Char(P)\to\dot\Sigma
\end{equation*}
is the restriction of $\pi$ to $\Char(P)$.
Then $\dot\Sigma$ has the subspace topology of $\Tb^*X$, and it can
also be topologized by $\hat\pi$, i.e.\ requiring that $C\subset\dot\Sigma$
is closed (or open) if and only if $\hat\pi^{-1}(C)$ is closed (or open).
These two topologies are equivalent, though the former is simpler
in the present setting -- e.g.\ it is immediate that $\dot\Sigma$
is metrizable. Lebeau~\cite{Lebeau:Propagation} (following Melrose's
original approach in the $\Cinf$ boundary setting, see
\cite{Melrose:Microlocal})
uses the latter;
in extensions of the present work, to allow e.g.\ iterated conic
singularities, that approach will be needed. Again, an analogous
situation arises in $N$-body
scattering, though that is in many respects more complicated
if some subsystems have bound states
\cite{Vasy:Propagation-Many, Vasy:Bound-States}.

We are now ready to define generalized broken bicharacteristics,
essentially following Lebeau \cite{Lebeau:Propagation}.
We say that a function $f$ on $T^*X\setminus o$ is $\pi$-invariant if
$f(q)=f(q')$
whenever $\pi(q)=\pi(q')$. In this case $f$ induces a function $f_\pi$
on $\dot\Tb^*X$ which satisfies $f=f_\pi\circ\pi$. Moreover, if $f$
is continuous, then so is $f_\pi$. Notice that if $f=\iota^* f_0$,
$f_0\in\Cinf(\Tb^*X)$, then $f\in\Cinf(T^*X)$ is certainly $\pi$-invariant.

\begin{Def}\label{Def:gen-br-bichar}
A generalized broken bicharacteristic of $P$
is a continuous map
$\gamma:I\to\dot\Sigma$,
where $I\subset\Real$ is an interval, satisfying
the following requirements:

\renewcommand{\theenumi}{\roman{enumi}}
\renewcommand{\labelenumi}{(\theenumi)}
\begin{enumerate}
\item
If $q_0=\gamma(t_0)\in\cG$
then for all $\pi$-invariant functions
$f\in\Cinf(T^*X)$,
\begin{equation}\label{eq:cG-bich}
\frac{d}{dt}(f_\pi\circ \gamma)(t_0)=H_p f(\tilde q_0),\ \tilde q_0=\pih^{-1}(q_0).
\end{equation}

\item
If $q_0=\gamma(t_0)\in\cH\cap T^*\cF_{i,\reg}$ then
there exists $\ep>0$ such that
\begin{equation}
t\in I,\ 0<|t-t_0|<\ep\Rightarrow\gamma(t)\nin T^*\cF_{i,\reg}.
\end{equation}

\item
If $q_0=\gamma(t_0)\in\cG\cap T^*\cF_{i,\reg}$, and $\cF_i$ is a boundary
hypersurface (i.e.\ has codimension $1$), then in a neighborhood of $t_0$,
$\gamma$ is a generalized broken bicharacteristic in the sense
of Melrose-Sj\"ostrand \cite{Melrose-Sjostrand:I}, see also
\cite[Definition~24.3.7]{Hor}.
\end{enumerate}
\end{Def}

Note that for $q_0\in\cG$, $\hat\pi^{-1}(\{q_0\})$ consists of a single
point, so \eqref{eq:cG-bich} makes sense. Moreover, (iii) implies (i)
if $q_0$ is in a boundary hypersurface, but it is stronger at diffractive
points, see \cite[Section~24.3]{Hor}.
The propagation of analytic singularities, as in Lebeau's case,
does not distinguish between gliding
and diffractive points, hence (iii) can be dropped to define
what we may call analytic generalized broken bicharacteristics. It is
an interesting question whether in the $\Cinf$ setting
there are also analogous diffractive phenomena at higher codimension
boundary faces, i.e.\ whether the following theorem can be strengthened
at certain points.

Our main result is:

\begin{thm*}(See Corollary~\ref{cor:prop-sing}.)
Suppose that $Pu=0$, $u\in H^1_{0,\loc}(X)$.
Then $\WFbz^{1,\infty}(u)\subset\dot\Sigma$, and it
is a union of maximally extended
generalized broken bicharacteristics of $P$ in $\dot\Sigma$.
\end{thm*}

A more precise version of this theorem, with microlocal assumptions on
$Pu$, is stated in Theorem~\ref{thm:prop-sing}. In particular, one
can allow $Pu\in\Cinf(X)$, which immediately implies that the
theorem holds for solutions of the wave equation with inhomogeneous
$\Cinf$ Dirichlet boundary conditions that match across the boundary
hyperfaces, see Remark~\ref{rem:inhom-Dir}.
In addition, this
theorem generalizes to the wave operator with Neumann boundary
conditions, which need to be interpreted in terms of the quadratic
form of $P$ (i.e.\ the Dirichlet form). That is, if $u\in H^1_{\loc}(X)$
satisfies
\begin{equation*}
\langle d_M u,d_Mv\rangle_X-\langle \pa_t u, \pa_t v\rangle_X=0
\end{equation*}
for all $v\in H^1_c(X)$, then $\WFbz^{1,\infty}(u)\subset\dot\Sigma$, and it
is a union of maximally extended
generalized broken bicharacteristics of $P$ in $\dot\Sigma$. In fact, the
proof of the theorem
for Dirichlet boundary conditions also utilizes the quadratic form
of $P$. It is slightly simpler in presentation only to the extent that
one has more flexibility to integrate by parts, etc., but in the
end the proof for Neumann boundary conditions simply requires a slightly
less conceptual (in terms of the traditions of microlocal analysis)
reorganization,
e.g.\ not using commutators $[P,A]$ directly, but commuting $A$ through the
exterior derivative $d_M$ and $\pa_t$ directly.

It is expected that these results will generalize to iterated edge-type
structures (under suitable hypotheses), whose simplest example
is given by conic points, recently analyzed by Melrose and Wunsch
\cite{Melrose-Wunsch:Propagation}, extending the product cone analysis
of Cheeger and Taylor \cite{Cheeger-Taylor:Diffraction}.

To make it clear what the main theorem states, we remark that
the propagation statement means that if $u$ solves $Pu=0$ (with, say,
Dirichlet boundary condition), and
$q\in\Tb^*_{\pa X}X\setminus o$ is such that $u$ has no singularities
on bicharacteristics entering $q$ (say, from the past), then
we conclude that $u$ has no
singularities at $q$, in the sense that $q\nin\WFbz^{1,\infty}(u)$, i.e.\ we only gain
b-derivatives (or totally
characteristic derivatives) microlocally. In particular, even if
$\WFbz^{1,\infty}(u)$ is empty, we can only conclude that $u$ is conormal to the
boundary, in the precise sense that
$V_1\ldots V_k u\in H^1_{\loc}(X)$ for any $V_1,\ldots,V_k
\in\Vb(X)$, and not that $u\in H^k_{\loc}(X)$ for all $k$. Indeed, the latter cannot
be expected to hold, as can be seen by considering e.g.\ the wave
equation (or even elliptic equations) in 2-dimensional conic sectors.

This already illustrates that from a technical point of view a major
challange is to combine two differential (and pseudodifferential)
algebras: $\Diff(X)$ and $\Diffb(X)$ (or $\Psib(X)$). The wave operator
$P$ lies in $\Diff(X)$, but microlocalization needs to take place in
$\Psib(X)$: if $\Psi(\tilde X)$ is the algebra of usual pseudodifferential
operators on an extension $\tilde X$ of $X$, its elements do not even
act on $\Cinf(X)$: see \cite[Section~18.2]{Hor} when $X$ has a smooth
boundary (and no corners). In addition, one needs
an algebra whose elements $A$
respect the boundary conditions, so e.g.\ $Au|_{\pa X}$
depends only on $u|_{\pa X}$ -- this is exactly the origin of the algebra of
totally
characteristic pseudodifferential operators, denoted by $\Psib(X)$,
in the $\Cinf$ boundary setting \cite{Melrose:Transformation}.
The interaction of these two algebras also explains why we prove even
microlocal elliptic regularity via the quadratic form of $P$
(the Dirichlet form), rather than by standard arguments, valid if
one studies microlocal elliptic regularity for an element of an algebra
(such as $\Psib(X)$)
with respect to the same algebra.

The ideas of the positive commutator estimates, in particular the construction
of the commutants, are very similar to those arising in the proof
of the propagation of singularities in $N$-body scattering in previous
works of the author -- the wave
equation corresponds to the relatively simple scenario there when no
proper subsystems have bound states \cite{Vasy:Propagation-Many}. Indeed,
the author has indicated many times in lectures that there is a close
connection between these two problems, and it is a pleasure to finally
spell out in detail
how the $N$-body methods can be adapted to the present setting.

The organization of the paper is as follows.
In Section~\ref{sec:b-calc} we recall basic facts about $\Psib(X)$ and
analyze its commutation properties with
$\Diff(X)$. In Section~\ref{sec:microlocal} we describe the
mapping properties of $\Psib(X)$ on $H^1(X)$-based spaces. We also define
and discuss the b-wave front set based on $H^1(X)$ there.
The following section is devoted to the elliptic estimates for the wave
equation. These are obtained from the microlocal positivity of
the Dirichlet form, which implies in particular that in this region commutators
are negligible
for our purposes. In Section~\ref{sec:bichar} we describe basic properties
of bicharacteristics, mostly relying on Lebeau's work
\cite{Lebeau:Propagation}. In Sections~\ref{sec:hyperbolic} and
\ref{sec:glancing}, we prove propagation estimates at hyperbolic,
resp.\ glancing, points, by positive commutator arguments. Similar
arguments were used by Melrose and Sj\"ostrand \cite{Melrose-Sjostrand:I}
for the analysis of propagation at glancing points for manifolds
with smooth boundaries, but the use of such arguments for hyperbolic
points is new even in the smooth boundary setting. (The usual arguments
utilize parametrices for microlocal Cauchy problems.)
In Section~\ref{sec:prop-sing} these results are combined to prove
our main theorems. The arguments presented there are very close
to those of Melrose, Sj\"ostrand and Lebeau.

Since the changes for Neumann
boundary conditions are
minor, and the arguments for Dirichlet boundary conditions can be stated
in a form closer to those found in classical microlocal analysis
(essentially, in the Neumann case one has to pay a price for integrating
by parts, so one needs to present the proofs in an appropriately rearranged,
and less transparent, form)
the proofs in the body of the paper are primarily written for
Dirichlet boundary conditions, and the required changes are pointed out
at the end of the various sections.

In addition, the hypotheses of the propagation of singularities theorem
can be relaxed to $u\in H^{1,m}_{\bl,0,\loc}(X)$,
$m\leq 0$, defined in Definition~\ref{Def:H1m-neg}.
Since this simply requires replacing the $H^1(X)$ norms by
the $H^{1,m}_{\bl}$ norms (which are only locally well defined), we
suppress this point except in the statement of the final result, to avoid
overburdening the notation. No changes are required in the argument to
deal with this more general case. See Remark~\ref{rem:sing-soln} for more
details.

To give the reader a guide as to what the real novelty is,
Sections~\ref{sec:b-calc}-\ref{sec:microlocal} should be considered as
variations on a well-developed theme. While some of the features
of microlocal analysis, especially
wave front sets, is not discussed on manifolds with corners elsewhere, the
modifications needed are essentially trivial (cf.\ \cite[Chapter~18]{Hor}).
A slight novelty is using
$H^1(X)$ as the point of reference for the b-wave front sets (rather
than simply weighted $L^2$ spaces), which is very useful later in the
paper, but again only demands minimal changes to standard arguments.
The discussions of bicharacteristics in Section~\ref{sec:bichar}
essentially quotes Lebeau's paper \cite[Section~III]{Lebeau:Propagation}.
Moreover, given the results of
Sections~\ref{sec:elliptic}, \ref{sec:hyperbolic} and
\ref{sec:glancing}, the proof of propagation
of singularities in Section~\ref{sec:prop-sing} is standard, essentially
due to Melrose and Sj\"ostrand \cite[Section~3]{Melrose-Sjostrand:II}.
Indeed, as presented by Lebeau
\cite[Proposition~VII.1]{Lebeau:Propagation}, basically
no changes are necessary at all in this proof.

The novelty is thus the
use of the Dirichlet form (hence the $H^1$-based wave front set)
for the proof of both the elliptic and hyperbolic/glancing
estimates, and the systematic used of positive commutator estimates
in the hyperbolic/glancing regions. This approach is quite robust, hence
significant extensions of the results
can be expected, as was already indicated.

I would like to thank Richard Melrose for his interest in this project,
for reading, and thereby improving, parts of the paper, and
for numerous helpful and stimulating
discussions, especially
for the wave equation on forms. While this topic did not become a part
of the paper, it did play a role in the presentation of the arguments here.
I am also grateful to Jared Wunsch for helpful discussions and his willingness
to read large parts of the manuscript at the early stages, when the
background material was still mostly absent; his help significantly
improved the presentation here. I would also like to thank Rafe Mazzeo for
his continuing interest in this project and for his patience
when I tried to explain him the main ideas in the early days of this
project.

\section{Interaction of $\Diff(X)$ with the b-calculus}
\label{sec:b-calc}
One of the main technical issues in proving our main theorem is that
unless $\pa X=\emptyset$,
the wave operator $P$ is {\em not} a b-differential operator:
$P\nin\Diffb^2(X)$. In this section we describe the basic
properties of how $\Diff^k(X)$, which includes $P$ for $k=2$, interacts
with $\Psib(X)$. We first recall though that for $p\in\cF_{i,\reg}$,
local coordinates in $\Tb^*X$ over a neighborhood of $p$ are
given by $(x,y,t,\sigma,\zeta,\tau)$ with $\sigma_j=x_j\xi_j$.
Thus, the map $\iota$ in local coordinates is $(x,y,t,\xi,\zeta,\tau)
\mapsto(x,y,t,x\xi,\zeta,\tau)$, where by $x\xi$ we mean the
vector $(x_1\xi_1,\ldots,x_k\xi_k)$.

In fact, in this section $y$ and $t$ play a completely analogous role,
hence there is no need to distinguish them at all. The difference
will only arise when we start studying the wave operator $P$
in Section~\ref{sec:elliptic}.
Thus, we let $\bar y=(y,t)$ and $\bar\zeta=(\zeta,\tau)$
here to simplify the notation.

We briefly recall basic properties of the set of `classical' (one-step
polyhomogeneous, in the sense that the full symbols are such on the
fibers of $\Tb^*X$) pseudodifferential operators
$\Psib(X)=\cup_m \Psib^m(X)$ and the set of standard (conormal)
b-pseudodifferential operators, $\Psibc(X)=\cup_m \Psibc^m(X)$. The
difference between these two classes is in terms of the behavior
of their (full) symbols at fiber-infinity of $\Tb^*X$: elements of
$\Psibc(X)$ have full symbols that
satisfy the usual symbol estimates, while elements of
$\Psib(X)$ have in addition an asymptotic expansion in terms of homogeneous
functions, so $\Psib^m(X)\subset\Psibc^m(X)$.
Conceptually, these are best defined
via the Schwartz kernel of $A\in\Psibc^m(X)$ in terms of a certain
blow-up $X^2_\bl$ of $X\times X$, see \cite{Melrose-Piazza:Analytic}
-- the Schwartz kernel is conormal to the lift $\diag_\bl$
of the diagonal of $X^2$
to $X^2_\bl$ with infinite order vanishing on all boundary
faces of $X^2_\bl$ which are disjoint from $\diag_\bl$. Modulo
$\Psib^{-\infty}(X)$, however, the explicit quantization map we give below
describes $\Psibc^m(X)$ and $\Psib^m(X)$. Here
$\Psibc^{-\infty}(X)=\Psib^{-\infty}(X)=\cap_m\Psibc^m(X)=\cap_m\Psib^m(X)$
is the ideal of smoothing operators. The topology of $\Psibc(X)$ is given
in terms the conormal seminorms of the Schwartz kernel $K$
of its elements; these seminorms can be stated in terms of
the Besov space norms of $L_1 L_2\ldots L_k K$ as $k$ runs over non-negative
integers, and the $L_j$ over first order differential operators tangential
to $\diag_\bl$, see \cite[Definition~18.2.6]{Hor}. Recall in particular
that these seminorms are (locally) equivalent to the
$\Cinf$ seminorms away from the lifted diagonal $\diag_\bl$.

There is a
principal symbol map $\sigma_{b,m}:\Psibc^m(X)\to
S^m(\Tb^*X)/S^{m-1}(\Tb^*X)$; here, for a vector bundle $E$ over $X$,
$S^k(E)$
denotes the set of symbols of order $k$ on $E$ (i.e.\ these are symbols
in the fibers of $E$, smoothly varying over $X$). Its restriction
to $\Psib^m(X)$ can be re-interpreted as a map $\sigma_{b,m}:\Psib^m(X)\to
\Cinf(\Tb^*X\setminus o)$ with values in homogeneous functions of degree $m$;
the range can of course also be identified with $\Cinf(\Sb^*X)$ if $m=0$
(and with sections of a line bundle over $\Sb^*X$ in general).
There is a short exact sequence
\begin{equation*}
0\longrightarrow \Psib^{m-1}(X)\longrightarrow \Psib^{m}(X)\longrightarrow
S^m(\Tb^*X)/S^{m-1}(\Tb^*X)\longrightarrow 0
\end{equation*}
as usual; the last non-trivial map is $\sigma_{b,m}$. There are also
quantization maps (which depend on various choices) $q=q_m:S^m(\Tb^*X)\to
\Psibc^m(X)$, which restrict to $q:S^m_{\cl}(\Tb^*X)\to
\Psib^m(X)$, $\cl$ denoting classical symbols, and $\sigma_{b,m}\circ q_m$
is the quotient map $S^m\to S^m/S^{m-1}$. For instance, over a local
coordinate chart $U$ as above, with $a$ supported in $\Tb^*_K X$,
$K\subset U$ compact, we may take, with $n=\dim X$,
\begin{equation}\begin{split}\label{eq:b-quantize}
&q(a)u(x,y)\\
&\qquad=(2\pi)^{-n}\int e^{i(x-x')\cdot \xi+(\yb-\yb')\cdot\bar\zeta}
\phi(\frac{x-x'}{x})
a(x,y,x\xi,\bar\zeta) u(x',\yb')\,dx'\,d\bar y'\,d\xi\,d\zeta,
\end{split}\end{equation}
understood as an oscillatory integral,
where $\phi\in\Cinf_c((-1/2,1/2)^k)$ is identically $1$ near $0$
and $\frac{x-x'}{x}=(\frac{x_1-x_1'}{x_1},\ldots,\frac{x_k-x_k'}{x_k})$,
and the integral in $x'$ is over $[0,\infty)^k$.
Here the role of $\phi$ is to ensure the infinite order vanishing
at the boundary hypersurfaces of $X^2_\bl$ disjoint from $\diag_\bl$;
it is irrelevant as far as the behavior of Schwartz kernels near the
diagonal is concerned (it is identically $1$ there).
This can
be extended to a global map via a partition of unity, as usual.
Locally, for $q(a)$, $\supp a\subset \Tb^*_K X$ as above, the conormal
seminorms of the Schwartz kernel of $q(a)$ (i.e.\ the Besov space
norms described above) can be bounded in terms of the symbol
seminorms of $a$, see the beginning of \cite[Section~18.2]{Hor}, and
conversely.
Moreover, any $A\in\Psibc(X)$ with {\em properly supported Schwartz kernel}
defines continuous linear maps
$A:\dCinf(X)\to\dCinf(X)$, $A:\Cinf(X)\to\Cinf(X)$.

\begin{rem}
We often do not state it below, but in general most pseudodifferential
operators have compact support in this paper. Sometimes
we use properly supported ps.d.o's, only for not having to state precise
support conditions; these are always composed with compactly supported
ps.d.o's or applied to compactly supported distributions, so effectively
they can be treated as compactly supported. See also Remark~\ref{rem:localize}.
\end{rem}

With $\tilde g$ being any $\Cinf$ Riemannian metric on $X$, and $K\subset X$
compact,
any $A\in\Psibc^0(X)$ with Schwartz kernel supported in $K\times K$
defines a bounded operator on $L^2(X)=L^2(X,d\tilde g)$, with norm bounded by
a seminorm of $A$ in $\Psibc^0(X)$.
Indeed, this is true for
$A\in\Psib^{-\infty}(X)$ with compact support, as follows from
the Schwartz lemma and the explicit description of the Schwartz kernel of $A$
on $X^2_\bl$. The standard square root argument then shows the
boundedness for $A\in\Psibc^0(X)$, with norm bounded by a seminorm of $A$
in $\Psibc^0(X)$ -- see \cite[Equation~(2.16)]{Melrose-Piazza:Analytic}.
In fact, we get more from the argument: letting $a=\sigma_{b,0}(A)$, there
exists $A'\in\Psib^{-1}(X)$ such that for all $v\in L^2(X)$,
\begin{equation*}
\|Av\|\leq 2\sup|a|\,\|v\|+\|A'v\|.
\end{equation*}
(The factor $2$ of course can be improved, as can the order of $A'$.)
This estimate will play an
important role in our propagation estimates -- it will take the place
of constructing a square root of the commutator, which would be
difficult here as we will commute $P$ with an element of $\Psib(X)$,
so the commutator will not lie in $\Psib(X)$.
We remark here that it is more usual to take a `b-density'
in place of $d\tilde g$, i.e.\ a globally non-vanishing section
of $\Omega_\bl^1 X=\Omega_\bl X$, which thus takes the form
$(x_1\ldots x_k)^{-1}\,d\tilde g$
locally near a codimension $k$ corner,
to define an $L^2$-space,
namely $L^2_\bl(X)=L^2(X,\frac{d\tilde g}
{x_1\ldots x_k})$; then $L^2(X)=x_1^{-1/2}\ldots x_k^{-1/2}L^2_\bl(X)$
appears as a weighted space. Elements of $\Psibc^0(X)$ are bounded on
both $L^2$ spaces, in the manner stated above.
The two boundedness results are very closely related, for if $A\in\Psibc^0(X)$,
then so is $x_j^{\lambda}A x_j^{-\lambda}$, $\lambda\in\Cx$.

There is an operator wave front set associated to $\Psibc(X)$ as well:
for $A\in\Psibc^m(X)$, $\WFb'(A)$ is a conic subset of $\Tb^*X\setminus o$,
and has the interpretation that $A$ is `in $\Psibc^{-\infty}(X)$'
outside $\WFb'(A)$. (We caution the reader that unlike the previous
material, as well as the rest of the background in the next three paragraphs,
$\WFb'$ is not discussed in \cite{Melrose-Piazza:Analytic}.
This discussion, however, is standard; see e.g.~\cite[Section~18.1]{Hor},
esp.\ after Definition~18.1.25, in the boundariless case, and
\cite[Section~18.3]{Hor} for the case of a $\Cinf$ boundary, where one
simply says that the operator is order $-\infty$ on certain open cones,
see e.g.\ the proof of Theorem~18.3.27 there.)
In particular, if $\WFb'(A)=\emptyset$, then
$A\in\Psib^{-\infty}(X)$. For instance, if $A=q(a)$, $a\in S^m(\Tb^* X)$,
$q$ as in \eqref{eq:b-quantize}, $\WFb'(A)$ is
defined by the requirement that if $p\nin \WFb'(A)$ then $p$ has a conic
neighborhood $U$
in $\Tb^*X\setminus o$ such that $A=q(a)$,
$a$ is rapidly decreasing in $U$, i.e.\
$|a(x,\yb,\sigma,\bar\zeta)|\leq C_N(1+|\sigma|+|\bar\zeta|)^{-N}$ for all $N$.
Thus, $\WFb'(A)$ is a closed conic subset of $\Tb^*X\setminus o$.
Moreover, if $K\subset\Sb^*X$ is compact, and $U$ is a neighborhood of $K$,
there exists $A\in\Psib^0(X)$ such that $A$ is the identity on $K$ and vanishes
outside $U$, i.e.\ $\WFb'(A)\subset U$, $\WFb'(\Id-A)\cap K=\emptyset$
-- we can construct $a$ to be homogeneous degree zero outside a neighborhood
of $o$, such that this homogeneous function regarded as a function
on $\Sb^*X$ (and still denoted by $a$) satisfies $a\equiv 1$ near $K$,
$\supp a\subset U$, and then let $A=q(a)$. (This roughly says that
$\Psib(X)$ can be used to localize in $\Sb^*X$, i.e.\ to b-microlocalize.)

$\Psibc(X)$ forms a filtered
$*$-algebra, so $A_j\in\Psibc^{m_j}(X)$, $j=1,2$, implies $A_1A_2\in
\Psibc^{m_1+m_2}(X)$, and $A_j^*\in\Psibc^{m_j}(X)$ with
\begin{equation*}
\sigma_{b,m_1+m_2}(A_1A_2)=\sigma_{b,m_1}(A_1)\sigma_{b,m_2}(A_2),
\ \sigma_{b,m_j}(A_j^*)=\overline{\sigma_{b,m_j}(A)}.
\end{equation*}
Here the formal adjoint is defined with respect to $L^2(X)$, the $L^2$-space
of any $\Cinf$ Riemannian metric on $X$; the same statements hold with
respect to $L^2_\bl(X)$ as well, since conjugation by $x_1\ldots x_k$
preserves $\Psibc^m(X)$ (as well as $\Psib^m(X)$), as already remarked
for $m=0$.
Moreover, $[A_1,A_2]\in\Psibc^{m_1+m_2-1}(X)$ with
\begin{equation*}
\sigma_{b,m_1+m_2-1}([A_1,A_2])=\frac{1}{i}\{a_1,a_2\},
\ a_j=\sigma_{b,m_j}(A_j);
\end{equation*}
$\{\cdot,\cdot\}$ is the Poisson bracket lifted from $T^*X$ via the
identification of $T^*X^\circ$ with $\Tb^*_{X^\circ}X$.
If $A_j\in\Psib^{m_j}(X)$, then $A_1 A_2\in \Psib^{m_1+m_2}(X)$,
$A_j^*\in\Psib^{m_j}(X)$, and $[A_1,A_2]\in \Psib^{m_1+m_2-1}(X)$.
In addition, operator composition satisfies
\begin{equation*}
\WFb'(A_1 A_2)
\subset\WFb'(A_1)\cap\WFb'(A_2).
\end{equation*}

If $A\in\Psibc^m(A)$ is elliptic, i.e.\ $\sigma_{b,m}(A)$ is invertible
as a symbol (with inverse in
$S^{-m}(\Tb^*X\setminus o)/S^{-m-1}(\Tb^*X\setminus o)$), then there
is a parametrix $G\in\Psibc^{-m}(X)$ for $A$, i.e.\ $GA-\Id, AG-\Id
\in\Psibc^{-\infty}(X)$. This construction microlocalizes, so
if $\sigma_{b,m}(A)$ is elliptic at $q\in\Tb^*X\setminus o$, i.e.\ 
$\sigma_{b,m}(A)$ is invertible
as a symbol in an open cone around $q$,
then there is a {\em microlocal parametrix} $G\in\Psibc^{-m}(X)$ for $A$
at $q$, so $q\nin\WFb'(GA-\Id)$, $q\nin\WFb'(AG-\Id)$, so $GA$, $AG$
are microlocally the identity operator near $q$. More generally,
if $K\subset\Sb^*X$ is compact, and $\sigma_{b,m}(A)$ is elliptic on $K$
then there is $G\in\Psibc^{-m}(X)$ such that $K\cap \WFb'(GA-\Id)=\emptyset$,
$K\cap\WFb'(AG-\Id)=\emptyset$. For $A\in\Psib^m(X)$, $\sigma_{b,m}(A)$
can be regarded as a homogeneous degree $m$ function on $\Tb^*X\setminus o$,
and ellipticity at $q$ means that $\sigma_{b,m}(A)(q)\neq 0$. For such $A$,
one can take $G\in\Psib^{-m}(X)$ in all the cases described above.

The other important ingredient, which however rarely appears in the
following discussion, although when it appears it is crucial, is
the notion of the indicial operator. This captures the mapping
properties of $A\in\Psib(X)$ in
terms of gaining any decay at $\pa X$. It plays
a role here as $P\nin\Diffb(X)$, so even if we do not expect to gain any
decay for solutions $u$ of $Pu=0$ say, we need to understand the commutation
properties of $\Diffb(X)$ with $\Psib(X)$, which will in turn follow from
properties of the indicial operator.
There is an indicial operator
map (which can also be considered as a non-commutative analogue of the
principal symbol), denoted by $\hat N_i$, for each boundary
face $\cF_i$, $i\in I$, and $\hat N_i$ maps $\Psibc^m(X)$ to a family
of b-pseudodifferential operators on $\cF_i$. For us, only the indicial
operators associated to boundary hypersurfaces $H_j$ will be important;
in this case the family is parameterized by $\sigma_j$,
the b-dual variable of $x_j$. It is characterized by the property that
if $f\in\Cinf(H_j)$ and $u\in\Cinf(X)$ is any extension of $f$,
i.e.\ $u|_{H_j}=f$, then
\begin{equation*}
\hat N_j(A)(\sigma_j)f=(x_j^{-i\sigma_j}Ax_j^{i\sigma_j}u)|_{H_j},
\end{equation*}
where $x_j^{-i\sigma_j}Ax_j^{i\sigma_j}\in\Psibc^m(X)$, hence
$x_j^{-i\sigma_j}Ax_j^{i\sigma_j}u\in\Cinf(X)$, and the
right hand side does not depend on the choice of $u$.
(In this formulation, we need
to fix $x_j$, at least mod $x_j^2\Cinf(X)$, to fix $\hat N_j(A)$. Note that
the radial vector field, $x_j D_{x_j}$, is independent of this choice
of $x_j$, at least modulo $x_j\Vb(X)$.)
If $A\in\Psibc^m(X)$ and $\hat N_i(A)=0$, then in fact $A\in\Cinf_{\cF_i}(X)
\Psibc^m(X)$, where $\Cinf_{\cF_i}(X)$ is the ideal of $\Cinf(X)$ consisting
of functions that vanish at $\cF_i$. In particular, for a boundary hypersurface
$H_j$ defined by $x_j$, if $A\in\Psibc^m(X)$ and $\hat N_j(A)=0$,
then $A=x_jA'$
with $A'\in\Psibc^m(X)$.
The indicial operators satisfy $\hat N_i(AB)=\hat N_i(A)\hat N_i(B)$.
The indicial family of $x_jD_{x_j}$ at $H_j$ is multiplication by $\sigma_j$,
while the indicial family of $x_k D_{x_k}$,
$k\neq j$, is $x_kD_{x_k}$ and that of $D_{\bar y_k}$ is $D_{\bar y_k}$.
In particular,
$\hat N_j([x_jD_{x_j},A])=[\hat N_j(x_j D_{x_j}),\hat N_j(A)]=0$, so
\begin{equation}\label{eq:basic-comm-Psib-Diff}
[x_jD_{x_j},A]\in x_j\Psibc^{m}(X),
\end{equation}
which plays a role below.
All of the above statements also hold with $\Psibc(X)$ replaced by
$\Psib(X)$.

The key point in analyzing smooth vector fields on $X$, and thereby
differential operators such as $P$ is that while $D_{x_j}\nin\Vb(X)$,
for any $A\in\Psib^{m}(X)$ there is an operator $\tilde A\in\Psib^m(X)$
such that
\begin{equation}\label{eq:V-Psib-1}
D_{x_j}A-\tilde A D_{x_j}\in\Psib^m(X),
\end{equation}
and analogously for $\Psib^m(X)$ replaced by $\Psibc^m(X)$. Indeed,
\begin{equation*}
D_{x_j}A=x_j^{-1}(x_j D_{x_j}) A=x_j^{-1}[x_jD_{x_j},A]+x_j^{-1}Ax_j D_{x_j}.
\end{equation*}
By \eqref{eq:basic-comm-Psib-Diff}, applied for $\Psib$ rather than
$\Psibc$,
\begin{equation*}
x_j^{-1}[x_jD_{x_j},A]\in\Psib^m(X).
\end{equation*}
Thus,
we may take $\tilde A=x_j^{-1}Ax_j$,
proving \eqref{eq:V-Psib-1}.
We also have, more trivially, that
\begin{equation}\label{eq:V-Psib-2}
D_{\bar y_j}A-\tilde A D_{\bar y_j}\in\Psib^m(X),\ \tilde A\in\Psib^m(X),
\ \sigma_{b,m}(A)=\sigma_{b,m}(\tilde A).
\end{equation}
Since $\sigma_{b,m}(A)=\sigma_{b,m}(x_j^{-1}Ax_j)$, we deduce the following
lemma.

\begin{lemma}\label{lemma:V-Psib-comm}
Suppose $V\in\Vf(X)$, $A\in\Psib^m(X)$. Then $[V,A]=\sum A_j V_j+B$
with $A_j\in\Psib^{m-1}(X)$, $V_j\in\Vf(X)$, $B\in\Psib^m(X)$.

Similarly, $[V,A]=\sum V_j A'_j +B'$
with $A'_j\in\Psib^{m-1}(X)$, $V_j\in\Vf(X)$, $B'\in\Psib^m(X)$.

Analogous results hold with $\Psib(X)$ replaced by $\Psibc(X)$.
\end{lemma}

\begin{proof}
It suffices to prove this for the coordinate vector fields, and indeed
just for the $D_{x_j}$. Then with the notation of \eqref{eq:V-Psib-1},
\begin{equation*}
D_{x_j}A-AD_{x_j}=(\tilde A-A)D_{x_j}+B,
\end{equation*}
and $\sigma_{b,m}(\tilde A)=\sigma_{b,m}(A)$, so $\tilde A-A\in\Psib^{m-1}(X)$,
proving the claim.
\end{proof}

More generally, we make the definition:

\begin{Def}
$\Diff^k\Psib^s(X)$ is the vector space of operators of the form
\begin{equation}\label{eq:Diff-Psib-def}
\sum_j P_j A_j,\ P_j\in\Diff^k(X),\ A_j\in\Psib^s(X),
\end{equation}
where the sum is locally finite in $X$.
\end{Def}

\begin{rem}
Since any point $q\in\Tb^*X\setminus o$ has a conic neighborhood $U$ in
$\Tb^*X\setminus o$ on which a vector field $V\in\Vb(X)$ is elliptic,
i.e.\ $\sigma_{b,1}(V)\neq 0$ on $U$, we can always write
$A_j\in\Psib^{s+k-k_j}(X)$ with $\WFb'(A)\subset U$, $k_j\leq k$,
as $A_j=Q_jA'_j+R_j$ with
$Q_j\in\Diffb^{k-k_j}(X)$, $A'_j\in\Psib^{s}(X)$, $R_j\in\Psib^{-\infty}(X)$.
Thus, any operator which is given by a locally finite sum of the form
\begin{equation*}
\sum_j P_j A_j,\ P_j\in\Diff^{k_j}(X),\ A_j\in\Psib^{s+k-k_j}(X),
\end{equation*}
can in fact be written in the form \eqref{eq:Diff-Psib-def}.
\end{rem}

\begin{lemma}
$\Diff^*\Psib^*(X)$ is filtered algebra with respect to operator
composition, with $B_j\in\Diff^{k_j}\Psib^{s_j}(X)$,
$j=1,2$, implying $B_1 B_2\in\Diff^{k_1+k_2}\Psib^{s_1+s_2}(X)$. Moreover,
with $B_1,B_2$ as above,
\begin{equation*}
[B_1,B_2]\in\Diff^{k_1+k_2}\Psib^{s_1+s_2-1}(X).
\end{equation*}
\end{lemma}

\begin{proof}
To prove that $\Diff^*\Psib^*(X)$ is an algebra,
we only need to prove that if $A\in\Psib^s(X)$, $P\in\Diff^k(X)$, then
$AP\in\Diff^k(X)\Psib^s(X)$. Writing $P$ as a sum of products of
vector fields in $\Vf(X)$, the claim follows from
Lemma~\ref{lemma:V-Psib-comm}.

Writing $B_j=V_{j,1}\ldots V_{j,k_1}A_j$, $A_j\in\Psib^{s_j}(X)$,
$V_{j,i}\in\Vf(X)$, and expanding the commutator $[B_1,B_2]$, one gets a finite
sum, each of which is a product of the factors
$V_{j,1},\ldots V_{j,k_1},A_j$ with two factors (one with $j=1$ and one with
$j=2$) removed and replaced by a commutator. In view of the first
part of the lemma, it suffices to note that
\begin{equation*}\begin{split}
&[V_{1,i},V_{2,i'}]\in\Vf(X),\ \Diff^{k_1+k_2-1}\Psib^{s_1+s_2}(X)
\subset\Diff^{k_1+k_2}\Psib^{s_1+s_2-1}(X),\\
&[A_1,A_2]\in\Psib^{s_1+s_2-1}(X)\\
&[V_{j,i},A_{3-j}]\in\Diff^1\Psib^{s_{3-j}-1}(X),
\end{split}\end{equation*}
where the last statement is a consequence of Lemma~\ref{lemma:V-Psib-comm},
taking into account that $\Psib^m(X)\subset\Diff^1\Psib^{m-1}(X)$.
\end{proof}

Although it is possible to define the principal symbol on $\Diff^k\Psib^s(X)$,
for technical reasons we will not use this in the proofs. Still, the
behavior of the principal symbol motivates the positive commutator
constructions at $\cG\cup\cH$, so we proceed to define it here. Thus,
using $\iota:T^*X\to\Tb^*X$, we can pull pack $\sigma_{b,s}(A)$, $A\in
\Psib^s(X)$, to $T^*X$, and define:

\begin{Def}
Suppose $B=\sum P_j A_j\in\Diff^k\Psib^s(X)$,
$P_j\in\Diff^k(X)$, $A_j\in\Psib^s(X)$. The principal symbol of $B$
is the $\Cinf$ homogeneous degree $k+s$ function on $T^*X\setminus o$ defined
by
\begin{equation}\label{eq:def-prs}
\sigma_{k+s}(B)=\sum \sigma_k(P_j)\iota^*\sigma_{b,s}(A_j).
\end{equation}
\end{Def}

\begin{lemma}
$\sigma_{k+s}(B)$ is independent of all choices.
\end{lemma}

\begin{proof}
Away from $\pa X$, $B$ is a pseudodifferential operator of order $k+s$,
and $\sigma_{k+s}(B)$ is its invariantly defined symbol. Since
the right hand side of \eqref{eq:def-prs} is continuous up to $\pa X$, and
is independent of all choices in $T^*X^\circ$, it
is independent of all choices in $T^*X$.
\end{proof}

We are now ready to compute the principal symbol of the commutator of
$A\in\Psib^m(X)$ with $D_{x_j}$.

\begin{lemma}\label{lemma:comm-symbol}
Let $\pa_{x_j}$, $\pa_{\sigma_j}$ denote local coordinate vector fields
on $\Tb^*X$ in the coordinates $(x,\bar y,\sigma,\bar\zeta)$.
For $A\in\Psib^m(X)$ with Schwartz kernel supported in the coordinate patch,
$a=\sigma_{b,m}(A)\in\Cinf(\Tb^* X\setminus o)$,
we have $[D_{x_j},A]=A_1D_{x_j}+A_0\in\Diff^1\Psib^{m-1}(X)$
with $A_0\in\Psib^m(X)$, $A_1\in\Psib^{m-1}(X)$ and
\begin{equation}\label{eq:sigma-A_j}
\sigma_{b,m-1}(A_1)
=\frac{1}{i}\pa_{\sigma_j}a,\ \sigma_{b,m}(A_0)=\frac{1}{i}\pa_{x_j}a.
\end{equation}
This result also holds with $\Psib(X)$ replaced by $\Psibc(X)$ everywhere.
\end{lemma}

\begin{rem}
Notice that $\sigma_m([D_{x_j},A])=\frac{1}{i}\{\xi_j,\iota^*a\}=\frac{1}{i}
\pa_{x_j}|_{\xi}$,
$\{.,.\}$ denoting
the Poisson bracket on $T^*X$ and $\pa_{x_j}|_{\xi}$ denoting the
appropriate
coordinate vector field on $T^*X$, i.e.\ where $\xi$ is held fixed
(rather than $\sigma$),
since both sides are continuous functions
on $T^*X\setminus o$ which agree on $T^*X^\circ\setminus o$. A simple
calculation shows that the lemma is consistent with this result.
The statement of the lemma would follow from this observation if
we showed that the kernel of $\sigma_{m}$ on $\Diff^1\Psib^{m-1}(X)$ is
$\Diff^1\Psib^{m-2}(X)$ -- the proof given below avoids this point
by reducing the calculation to $\Psib(X)$.
\end{rem}

\begin{proof}
The lemma follows from
\begin{equation*}
D_{x_j}A-AD_{x_j}=x_j^{-1}[x_jD_{x_j},A]+x_j^{-1}[A,x_j]D_{x_j}.
\end{equation*}
Indeed, letting
\begin{equation}\label{eq:A_j-def}
A_0=x_j^{-1}[x_jD_{x_j},A]\in\Psib^m(X),
\ A_1=x_j^{-1}[A,x_j]\in\Psib^{m-1}(X),
\end{equation}
the principal symbols can be calculated
in the b-calculus. Since they are given by the standard Poisson bracket
in $T^*X^\circ$, hence in $\Tb^*_{X^\circ}X$, by continuity the same
calculation gives a valid result in $\Tb^*X$.
As $\pa_{\xi_j}=x_j\pa_{\sigma_j}$, $\pa_{x_j}|_{\xi}
=\pa_{x_j}|_{\sigma}+\xi_j\pa_{\sigma_j}$,
we see that for $b=\sigma_j$ or $b=x_j$,
the Poisson bracket $\{b,a\}$
is given by
\begin{equation*}\begin{split}
&x_j(\pa_{\sigma_j}b)(\pa_{x_j}|_{\sigma}a+\xi_j\pa_{\sigma_j}a)
-x_j(\pa_{\sigma_j}a)(\pa_{x_j}|_{\sigma}b+\xi_j\pa_{\sigma_j}b)\\
&\qquad=x_j(\pa_{\sigma_j}b)\pa_{x_j}|_{\sigma}a
-x_j(\pa_{\sigma_j}a)\pa_{x_j}|_{\sigma}b
\end{split}\end{equation*}
so we get
\begin{equation*}
\{\sigma_j,a\}=x_j\pa_{x_j}|_{\sigma}a,\ \{x_j,a\}=-x_j\pa_{\sigma_j}a,
\end{equation*}
so \eqref{eq:sigma-A_j} follows from \eqref{eq:A_j-def}.
\end{proof}

\section{Function spaces and microlocalization}\label{sec:microlocal}

We now turn to action of $\Psib(X)$
on function spaces related to differential operators in $\Diff(X)$,
and in particular $H^1(X)$ which corresponds to first order differential
operators, such as the exterior derivative $d$. We first recall that
$\Cinf_c(X)$ is the space of $\Cinf$ functions of compact support on $X$
(which may thus be non-zero at $\pa X$), while
$\dCinf_c(X)$ is the subspace of $\Cinf_c(X)$ consisting of functions
which vanish to infinite order at $\pa X$. Although we will mostly consider
local results, and any $\Cinf$ Riemannian metric can be used to define
$L^2_{\loc}(X)$, $L^2_c(X)$ (as different choices give the same space),
it is convenient to fix a global Riemmanian metric, $\tilde g=g+dt^2$,
on $X$, where $g$ is the metric on $M$. With this choice, $L^2(X)$ is
well-defined as a Hilbert space. For $u\in\Cinf_c(X)$, we let
\begin{equation*}
\|u\|^2_{H^1(X)}=\|du\|^2_{L^2(X)}+\|u\|^2_{L^2(X)}.
\end{equation*}
We then let $H^1(X)$ be the completion of $\Cinf_c(X)$ with respect to
the $H^1(X)$ norm. Then we define $H^1_0(X)$ as the
closure of $\dCinf_c(X)$ inside $H^1(X)$.

\begin{rem}\label{rem:Sobolev}
We recall alternative viewpoints of these Sobolev spaces.
Good references for the $\Cinf$ boundary case (and no corners) include
\cite[Appendix~B.2]{Hor} and
\cite[Section~4.4]{Taylor:Partial-I}; only minor modifications are
needed to deal with the corners for the special cases we discuss below.

We can define $H^1(X^\circ)$ as the subspace
of $L^2(X)$ consisting of functions $u$
such that $du$, defined as the distributional derivative of
$u$ in $X^\circ$, lying in $L^2(X,\Lambda^1 X)$; we then equip it with
the above norm -- this is locally equivalent to saying that
$V u\in L^2_{\loc}(X)$ for all $\Cinf$ vector fields $V$ on $X$,
where $Vu$ refers
to the distributional derivative of $u$ on $X^\circ$.

In fact, $H^1(X^\circ)=H^1(X)$,
since $H^1(X^\circ)$ is complete with respect to the $H^1$ norm and
$\Cinf_c(X)$ is easily seen to be dense in it.
For instance, locally, if $X$ is given by $x_j\geq 0$, $j=1,\ldots,k$,
and $u$ is supported in such a coordinate chart, one
can take $u_s(x,\bar y)=u(x_1+s,\ldots,x_k+s,\bar y)$ for $s>0$, and see that
$u_s|_X\to u$ in $H^1_c(X^\circ)$. Then a standard regularization argument
on $\Real^n$, $n=\dim X$, gives the claimed density
of $\Cinf_c(X)$ in $H^1_c(X^\circ)$. Thus, $H^1(X^\circ)=H^1(X)$ indeed,
which shows in particular that $H^1(X)\subset L^2(X)$. (Note that
$\|u\|_{L^2(X)}\leq\|u\|_{H^1(X)}$ only guarantees that there is a continuous
`inclusion' $H^1(X)\hookrightarrow L^2(X)$, not that it is injective,
although that can be proved easily by a direct argument, cf.\ the
Friedrichs extension method for operators, see
e.g.\ \cite[Theorem~X.23]{Reed-Simon}.)

If $\tilde X$ is a manifold without boundary, and $X$ is embedded into it,
one can also extend elements of $H^1(X)$ to elements $H^1_{\loc}(\tilde X)$
exactly as in the $\Cinf$ boundary case (or simply locally
extending in $x_1$ first,
then in $x_2$, etc., and using the $\Cinf$ boundary result), see
\cite[Section~4.4]{Taylor:Partial-I}. Thus, with the notation of
\cite[Appendix~B.2]{Hor}, $H^1_{\loc}(X)=\bar H^1_{\loc}(X^\circ)$.
As is clear from
the completion definition, $H^1_{0,\loc}(X)$ can be identified with the subset
of $H^1_{\loc}(\tilde X)$ consisting of functions supported in $X$.
Thus, $H^1_{0,\loc}(X)=\dot H^1_{\loc}(X)$ with the notation of
\cite[Appendix~B.2]{Hor}.

All of the above discussion can be easily modified for $H^m$ in place of $H^1$,
$m\geq 0$ an integer.
\end{rem}

We are now ready to state the action on Sobolev spaces. These results would
be valid, with similar proofs, if we replace $H^1(X)$ by $H^m(X)$, $m\geq 0$
integer. We also refer to \cite[Theorem~18.3.13]{Hor} for further
extensions when $X$ has a $\Cinf$
boundary (and no corners).

\begin{lemma}\label{lemma:H10-bd}
Any $A\in\Psibc^0(X)$ with compact support defines a continuous linear
maps $A:H^1(X)\to H^1(X)$, $A:H^1_0(X)\to H^1_0(X)$,
with norms bounded by a seminorm of $A$
in $\Psibc^0(X)$.

Moreover, for any $K\subset X$ compact, any $A\in\Psibc^0(X)$
with proper support
defines a continuous map from the subspace of $H^1(X)$ (resp.\ $H^1_0(X)$)
consisting of distributions supported in $K$ to $H^1_c(X)$
(resp.\ $H^1_{0,c}(X)$).
\end{lemma}

\begin{rem}
Note that all smooth vector fields $V$ of compact support
define a continuous
operator $H^1(X)\to L^2(X)$, so in particular $V\in\Vb(X)$ do so.
Now, any $A\in\Psibc^1(X)$ can be written as $\sum (D_{x_j}x_j)A_j
+\sum D_{\bar y_j}A'_j+A''$ with $A_j,A'_j,A''\in\Psibc^0(X)$ by writing
$\sigma_{b,1}(A)=\sum \sigma_j a_j+\sum \bar\zeta_j a'_j$, and taking
$A_j,A'_j$ with principal symbol $a_j,a'_j$.
Therefore the lemma implies that
any $A\in\Psibc^1(X)$ defines a continuous linear operator
$H^1(X)\to L^2(X)$, and in particular restricts to a map $H^1_0(X)\to L^2(X)$.
\end{rem}

\begin{proof}
For $A\in\Psibc^0(X)$, by \eqref{eq:V-Psib-1}
$D_{x_j} Au=\tilde A D_{x_j}u +Bu$, with $\tilde A\in\Psibc^0(X)$,
$B\in\Psibc^0(X)$ the seminorms of both in $\Psibc^0(X)$ bounded by
seminorms of $A$ in $\Psibc^0(X)$, so by the first
half of the proof
\begin{equation*}
\|D_{x_j} Au\|_{L^2(X)}
\leq \|\tilde A\|_{\bop(L^2(X),L^2(X))}\|D_{x_j} u\|_{L^2(X)}
+\|B\|_{\bop(L^2(X),L^2(X))}\|u\|_{L^2(X)}.
\end{equation*}
Since there is an analogous formula for $D_{x_j}$ replaced by $D_{\bar y_j}$,
we deduce that for some $C>0$, depending only on a seminorm of
$A$ in $\Psibc^0(X)$,
\begin{equation*}
\|d_X Au\|_{L^2(X)}\leq C(\|d_X u\|_{L^2(X)}+\|u\|_{L^2(X)}).
\end{equation*}

Thus, $A\in\Psibc^0(X)$ extends to a continuous linear map from the
completion of $\Cinf_c(X)$ with respect to the $H^1(X)$ norm
to itself, i.e.\ from $H^1(X)$ to itself as claimed. As it maps
$\dCinf_c(X)\to\dCinf_c(X)$, it also maps the $H^1$-closure of
$\dCinf(X)$ to itself, i.e.\ it defines a continuous linear
map $H^1_0(X)\to H^1_0(X)$, finishing the proof of the first half of the
lemma.

For the second half, we only need to note that $Au=A\phi u$ if $\phi\equiv 1$
near $K$ and has compact support; now $A\phi$ has compact support so
the first half of the lemma is applicable.
\end{proof}

Note that $H^1(X)\subset L^2(X)\subset\dist(X)$, with $\dist(X)$ denoting
the dual space of $\dCinf_c(X)$, i.e.\ the space of extendible distributions.
Since for any $m$, $A\in\Psibc^m(X)$ maps $\dist(X)\to\dist(X)$, we could
view $A$ already defined as a map $H^1(X)\to\dist(X)$; then the above
lemma is a continuity result for $m=0$.

We let $H^{-1}(X)$ be the dual of $H^1_0(X)$ and $\dot H^{-1}(X)$ be
the dual of $H^1(X)$, with respect to an extension of the sesquilinear form
$\langle u,v\rangle=
\int_X u\,\overline v\,d\tilde g$, i.e.\ the $L^2$ inner product.
As $H^1_0(X)$ is a closed subspace of $H^1(X)$, $H^{-1}(X)$ is the quotient
of $\dot H^{-1}(X)$ by the annihilator of $H^1_0(X)$.
In terms of the identification of the $H^1$
spaces in the penultimate
paragraph of Remark~\ref{rem:Sobolev}, $H^{-1}_{\loc}(X)
=\bar H^{-1}_{\loc}(X^\circ)$ in the notation of \cite[Appendix~B.2]{Hor},
i.e.\ its elements are the restrictions to $X^\circ$ of elements
of $H^{-1}_{\loc}(\tilde X)$. Analogously, $\dot H^{-1}_{\loc}(X)$ consists
of those elements of $H^{-1}_{\loc}(\tilde X)$ which are supported in $X$.

Any $V\in\Diff^1(X)$ of compact support
defines a continuous map $L^2(X)\to H^{-1}(X)$
via $\langle Vu,v\rangle=\langle u,V^* v\rangle$ for $u\in L^2(X)$,
$v\in H^1_0(X)$; this is the same map as induced by extending $V$ to
an element $\tilde V$ of $\Diff^1(\tilde X)$, extending $u$ to $\tilde X$,
say as $0$, and letting $Vu=\tilde V\tilde u|_{X^\circ}$. Thus,
any $P\in\Diff^2(X)$ of compact support
defines continuous maps $H^1(X)\to H^{-1}(X)$,
and in particular $H^1_0(X)\to H^{-1}(X)$, since we can write $P=\sum V_j W_j$
with $V_j,W_j\in\Diff^1(X)$. Similarly, any $P\in\Diff^2(X)$ defines
continuous maps $H^1_{\loc}(X)\to H^{-1}_{\loc}(X)$,
and in particular $H^1_{0,\loc}(X)\to H^{-1}_{\loc}(X)$. Thus,
for $P=\Delta_{\tilde g}+1$, $\langle u,v\rangle_{H^1(X)}=\langle u,Pv\rangle$
if $u\in H^1_0(X)$ and $v\in H^1(X)$. Similarly, for $P=D_t^2-\Delta_g$,
$\langle D_t u,D_t v\rangle-\langle d_M u,d_M v\rangle
=\langle u,Pv\rangle$, if $u\in H^1_0(X)$ and $v\in H^1(X)$.

We also remark that as $H^1(X)$ and $H^1_0(X)$ are Hilbert spaces, their
duals are naturally identified with themselves via the inner product.
Thus, if $f$ is a continuous linear functional on $H^1_0(X)$, then
there is a $v\in H^1_0(X)$ such that $f(u)=\langle u,v\rangle
+\langle du,dv\rangle$. Thus, regarding $H^1_0(X)$ as a subspace
of $H^1(\tilde X)$, for an extension $\tilde X$ of $X$, as in
Remark~\ref{rem:Sobolev}, we deduce that
$f(u)=\langle u,(\Delta_{\tilde g}+1)v\rangle$, so the identification of
$H^{-1}(X)$ with $H^1_0(X)$ (regarded as its own dual) is given by
$H^1_0(X)\ni v\mapsto (\Delta_{\tilde g}+1)v\in H^{-1}(X)$.

Since $\Psibc^0(X)$ is closed
under taking adjoints, the following result is an immediate consequence
of Lemma~\ref{lemma:H10-bd}.

\begin{cor}\label{cor:H-1-bd}
Any $A\in\Psibc^0(X)$ with compact support defines a continuous linear
maps $A:H^{-1}(X)\to H^{-1}(X)$, $A:\dot H^{-1}(X)\to \dot H^{-1}(X)$,
with norm bounded by a seminorm of $A$
in $\Psibc^0(X)$.
\end{cor}

We now define subspaces of $H^1(X)$ which possess additional regularity
with respect to $\Psib(X)$.

\begin{Def}
For $m\geq 0$, we define $H^{1,m}_{b,c}(X)$
as the subspace of $H^1(X)$
consisting of $u\in H^1(X)$ with $\supp u$ compact and
$Au\in H^1(X)$ for some (hence any, as shown below)
$A\in\Psib^m(X)$ (with compact support) which is elliptic over $\supp u$,
i.e.\ $A$ such that such that
$\sigma_{b,m}(A)(q)\neq 0$ for any
$q\in\Tb^*_{\supp u}X\setminus o$.

We let $H^{1,m}_{b,\loc}(X)$ be the subspace of $H^1_{\loc}(X)$
consisting of $u\in H^1_{\loc}(X)$ such that for any $\phi\in\Cinf_c(X)$,
$\phi u\in H^{1,m}_{b,c}(X)$.

We also let $H^{1,m}_{b,0,c}(X)=H^{1,m}_{b,c}(X)\cap H^1_0(X)$, and
similarly for the local space $H^{1,m}_{b,0,\loc}(X)$.
\end{Def}

\begin{rem}
The definition is independent of the choice of $A$, as can be seen
by taking a parametrix $G\in\Psib^{-m}(X)$ for $A$ in a neighborhood of
$\supp u$, so $GA-\Id=E\in\Psib^0(X)$, and $\WFb'(E)\cap
\Tb^*_{\supp u}X\setminus o=\emptyset$. Indeed, let $\rho\in\Cinf_c(X)$ be
identically $1$ near $\supp u$, $\WFb'(E)\cap\Tb^*_{\supp\rho}X=\emptyset$.
Then any $A'$ with
the properties of $A$ can be written as $A'=A'GA-A'E\rho-A'E(1-\rho)$,
$A'G,A'E\rho\in\Psib^0(X)$, while $(1-\rho)u=0$,
so by Lemma~\ref{lemma:H10-bd}, $A'u\in H^1(X)$ provided that $u,Au\in H^1(X)$.
\end{rem}

It is useful to note that if $Au\in H^1(X)$ and $u\in H^1_0(X)$, then
in fact $Au\in H^1_0(X)$:

\begin{lemma}\label{lemma:A-H1-H10}
Suppose that $u\in H^1_0(X)$, $A\in\Psib^m(X)$ and $Au\in H^1(X)$. Then
$Au\in H^1_0(X)$.
\end{lemma}

\begin{proof}
Suppose that $u\in H^1_0(X)$, $A\in\Psib^m(X)$ and $Au\in H^1(X)$.
Let $\Lambda_r$, $r\in(0,1]$,
be a uniformly bounded family in $\Psibc^0(X)$ with
$\Lambda_r\in\Psib^{-\infty}(X)$ for $r>0$, $\Lambda_r\to\Id$ in
$\Psib^{\ep}(X)$, $\ep>0$, as $r\to 0$.

Then, for $r>0$, $\Lambda_r A\in\Psib^{-\infty}(X)$, so $u\in H^1_0(X)$
implies that $\Lambda_r Au\in H^1_0(X)$ by Lemma~\ref{lemma:H10-bd}.
As $Au\in H^1(X)$, and $\Lambda_r$ is uniformly bounded as a family
of operators on $H^1(X)$, we deduce that $\Lambda_r Au$ is uniformly
bounded in $H^1(X)$. Thus, there is a weakly convergent sequence
$\Lambda_{r_j}Au$, with $r_j\to 0$, in
$H^1_0(X)$, as the latter is a closed subspace of $H^1(X)$; let $v$ be
the limit. But
$\Lambda_r A u\to Au$ in $\dist(X)$ as $r\to 0$, since $\Lambda_rA\to A$
in $\Psibc^{m+\ep}(X)$.
As $\Lambda_{r_j} A u\to v$ in $\dist(X)$
as well, $Au=v\in H^1_0(X)$ as claimed.
\end{proof}

The following wave front set microlocalizes $H^{1,m}_{b,\loc}(X)$.

\begin{Def}
Suppose $u\in H^1_{\loc}(X)$, $m\geq 0$.
We say that $q\in\Tb^*X\setminus o$ is not in $\WFbz^{1,m}(u)$ if
there exists $A\in\Psib^m(X)$ such that $\sigma_{b,m}(A)(q)\neq 0$
and $Au\in H^1(X)$.

For $m=\infty$, we say that $q\in\Tb^*X\setminus o$ is not in $\WFbz^{1,m}(u)$
if there exists $A\in\Psib^0(X)$ such that $\sigma_{b,0}(A)(q)\neq 0$
and $LAu\in H^1(X)$ for all $L\in\Diffb(X)$, i.e.\ if
$Au\in H^{1,\infty}_b(X)$.
\end{Def}

We note that, by the preceeding lemma, if $u\in H^1_{0,\loc}(X)$ then
$Au\in H^1_{0,\loc}(X)$, etc. (here $A\in\Psib^m(X)$).
Moreover, in the $m$ infinite case we may equally allow $L\in
\Psib(X)$, and we can also rewrite the finite $m$ definition analogously,
i.e.\ to state that there exists $A\in\Psib^0(X)$ such that
$\sigma_{b,0}(A)(q)\neq 0$
and $LAu\in H^1(X)$ for all $L\in\Psib^m(X)$ -- this follows immediately
from the next lemma. Although we do not need this here, so we
do not comment on it any more,
we could also allow $A\in\Psibc^m(X)$ in the
definition, provided we replace $\sigma_{b,m}(A)(q)\neq 0$
by the assumption that $A$ is elliptic at $q$ -- this follows
from the next results.

The following lemma shows that the action of elements of $\Psib(X)$ is
indeed microlocal.

\begin{lemma}\label{lemma:WFb-mic}
Suppose that $u\in H^1_{\loc}(X)$, $B\in\Psibc^k(X)$. Then $\WFbz^{1,m-k}(Bu)
\subset\WFbz^{1,m}(u)\cap\WFb'(B)$.
\end{lemma}

\begin{proof}
We assume that $m$ is finite; the proof for $m$ infinite is similar.

Suppose $q\nin\WFb'(B)$. As $\WFb'(B)$ is closed, there is a neighborhood
$U$ of $q$ such that $U\cap\WFb'(B)=\emptyset$. Let $A\in\Psib^{m-k}(X)$
satisfy $\WFb'(A)\subset U$, $\sigma_{b,m-k}(A)(q)\neq 0$. Then
$AB\in\Psib^{-\infty}(X)\subset \Psib^0(X)$, so $ABu\in H^1(X)$ by
Lemma~\ref{lemma:H10-bd}. Thus, $q\nin\WFbz^{1,m-k}(Bu)$ by definition
of the wave front set.

On the other hand, suppose that $q\nin\WFbz^{1,m}(u)$. Then there
is some $A\in\Psib^{m}(X)$ such that $Au\in H^1(X)$ and
$\sigma_{b,m}(A)(q)\neq 0$. Let $G\in\Psib^{-m}(X)$ be a microlocal
parametrix for $A$, so $GA=\Id+E$ with $E\in\Psib^0(X)$, $q\nin
\WFb'(E)$. Let $C\in\Psib^{m-k}(X)$ be such that $\WFb'(C)\cap\WFb'(E)
=\emptyset$ and $\sigma_{b,m-k}(C)(q)\neq 0$. Then
$CBE\in\Psib^{-\infty}(X)$, so $CBEu\in H^1(X)$ by Lemma~\ref{lemma:H10-bd}.
On the other hand, $CBG\in\Psibc^{0}(X)$ and $Au\in H^1(X)$, so
$CBGAu\in H^1(X)$ also by Lemma~\ref{lemma:H10-bd}. We thus deduce
that $CBu=CBGAu-CBEu\in H^1(X)$, so $q\nin\WFbz^{1,m-k}(u)$.
\end{proof}

We will need a quantitative version of this lemma
giving actual estimates,
but first we state the precise sense in which
this wave front set provides a refined version of the conormality of $u$.

\begin{lemma}\label{lemma:microloc-to-loc}
Suppose $u\in H^1_{\loc}(X)$, $m\geq 0$, $p\in X$. If $\Sb^*_p X\cap
\WFbz^{1,m}(u)=\emptyset$, then in a neighborhood of $p$,
$u$ lies in $H^{1,m}_b(X)$, i.e.\ there is $\phi\in\Cinf_c(X)$ with
$\phi\equiv 1$ near $p$ such that $\phi u\in H^{1,m}_b(X)$.
\end{lemma}

\begin{proof}
We assume that $m$ is finite; the proof for $m$ infinite is similar.

For each $q\in \Sb^*_pX$ there is $A_q\in\Psib^m(X)$ such that
$\sigma_{b,m}(A_q)(q)\neq 0$ and $A_q u\in H^1(X)$. Let $U_q$ be the
set on which $\sigma_{b,m}(A_q)\neq 0$; then $U_q$ is an open set
containing $q$. Thus, $\{U_q:\ q\in\Sb^*_pX\}$ is an open cover
of the compact set $\Sb^*_pX$. Let $U_{q_j}$, $j=1,\ldots,r$ be
a finite subcover. Then $A_0=\sum A_{q_j}^*A_{q_j}$ is elliptic
on $\Sb^*_pX$ since $\sigma_{b,2m}(A_0)=\sum |\sigma_{b,m}(A_{q_j})|^2$,
with each summand non-negative, and at any $q\in\Sb^*_pX$ at least
one term is nonzero (namely one for which $q\in U_{q_j}$).
Finally, we renormalize $A_0$ to make its order the same as that of $A$:
this is achieved by taking any $Q\in\Psib^{-m}(X)$ which is elliptic
on $\Sb^*_pX$, and letting $A=QA_0\in\Psib^{m}(X)$. Thus, $A$ is elliptic
on $\Sb^*_pX$, and
$Au\in H^1(X)$ as this holds for each summand $(QA_{q_j}^*)(A_{q_j}u)$,
for $QA_{q_j}^*\in\Psib^0(X)$ and $A_{q_j}u\in H^1(X)$. Here we used
Lemma~\ref{lemma:H10-bd}.

Let $G\in\Psib^{-m}(X)$ be a microlocal parametrix for $A$, so $GA=\Id+E$ and
$\WFb'(E)\cap\Sb^*_pX=\emptyset$. Thus, $p$ has a neighborhood $O$ in
$X$ such that $\WFb'(E)\cap \Sb^*_OX=\emptyset$. Let $\phi\in\Cinf_c(X)$
be supported in $O$, identically $1$ near $p$, and let $T\in\Psib^m(X)$ be
elliptic on $\Sb^*_{\supp\phi}X$. Then
$T\phi u=T\phi GAu-T\phi Eu$. Since $\WFb'(E)\cap\WFb'(\phi)=\emptyset$, we see
that $T\phi E\in\Psib^{-\infty}(X)$, and thus the
last term is in $H^1(X)$ by Lemma~\ref{lemma:H10-bd}. On the other hand,
the first term is in $H^1(X)$ since $Au\in H^1(X)$ and
$T\phi G\in\Psib^0(X)$. Thus, $\phi u\in
H^{1,m}_b(X)$ as claimed.
\end{proof}

\begin{cor}\label{cor:WF-to-H1}
If $u\in H^1_{\loc}(X)$
and $\WFbz^{1,m}(u)=\emptyset$, then $u\in H^{1,m}_{b,\loc}(X)$.

In particular, if $u\in H^1_{\loc}(X)$
and $\WFbz^{1,m}(u)=\emptyset$ for all $m$,
then $u\in H^{1,\infty}_{b,\loc}(X)$, i.e.\ $u$ is conormal in the sense
that $Au\in H^1_{\loc}(X)$ for all $A\in\Diffb(X)$ (or indeed $A\in\Psib(X)$).
\end{cor}

For the quantitative version of Lemma~\ref{lemma:WFb-mic}
we need a notion of the operator wave front
set that is uniform in a family of operators:

\begin{Def}
Suppose that $\cB$ is a bounded subset of $\Psibc^k(X)$, and $q\in\Sb^*X$.
We say that $q\nin\WFb'(\cB)$ if there is some $A\in\Psib(X)$
which is elliptic at $q$ such that $\{AB:\ B\in\cB\}$ is a bounded
subset of $\Psib^{-\infty}(X)$.
\end{Def}

Note that the wave front set of a family $\cB$
is only defined for bounded families. It can
be described directly in terms of quantization of (full) symbols,
much like the operator wave front set of a single operator. All
standard properties of the operator wave front set also hold for a family;
e.g.\ if $E\in\Psib(X)$ with $\WFb'(E)\cap\WFb'(\cB)=\emptyset$ then
$\{BE:\ B\in\cB\}$ is bounded in $\Psib^{-\infty}(X)$.

A quantitative version of Lemma~\ref{lemma:WFb-mic} is the following result.

\begin{lemma}\label{lemma:WFb-mic-q}
Suppose that $K\subset\Sb^*X$ is compact, and $U$ a neighborhood of $K$
in $\Sb^*X$. Let $\tilde K\subset X$ compact,
and $\tilde U$ is a neighborhood of $\tilde K$ in $X$ with compact closure.
Let $Q\in\Psib^k(X)$ be elliptic on $K$ with $\WFb'(Q)
\subset U$, with Schwartz kernel supported in $\tilde K\times\tilde K$.
Let $\cB$ be a bounded subset of $\Psibc^k(X)$
with $\WFb'(\cB)\subset K$ and
Schwartz kernel supported in $\tilde K\times\tilde K$.
Then there is a constant $C>0$ such that for $B\in\cB$, $u\in H^1_{\loc}(X)$
with $\WFbz^{1,k}(u)\cap U=\emptyset$,
\begin{equation*}
\|Bu\|_{H^1(X)}\leq C(\|u\|_{H^1(\tilde U)}+\|Qu\|_{H^1(X)}).
\end{equation*}
\end{lemma}

\begin{proof}
Let $\phi\in\Cinf_c(\tilde U)$
be identically $1$ near $\tilde K$.
We may replace $u$ by $\phi u$ in the estimate since $B\phi=B$, $Q\phi=Q$;
then $\|\phi u\|_{H^1(\tilde U)}=\|\phi u\|_{H^1(X)}$.

By Lemma~\ref{lemma:WFb-mic} and Lemma~\ref{lemma:microloc-to-loc},
all terms in the estimate are finite, since e.g.\ $\WFb'(Q)\cap\WFbz^{1,k}(u)
=\emptyset$ so $\WFbz^{1,0}(u)=\emptyset$, so $Qu\in H^{1,0}_{b,\loc}(X)
=H^1_{\loc}(X)$,
and indeed $Qu\in H^{1}_{c}(X)$, as the Schwartz kernel of $Q$ has compact
support.

Let $G$ be a microlocal parametrix for $Q$, so $GQ=\Id+E$ with
$E\in\Psib^0(X)$, $\WFb'(E)\cap K=\emptyset$. Thus,
$Bu=BGQu-BEu$. Now, $BE\in\Psib^{-\infty}(X)$ since $\WFb'(E)\cap K=\emptyset$
and $\WFb'(B)\subset K$, and it lies in a bounded subset of
$\Psib^{-\infty}(X)$ for $B\in\cB$. Thus,
$\|BEu\|_{H^1(X)}\leq C_1\|u\|_{H^1(X)}$ by Lemma~\ref{lemma:H10-bd}.
On the other hand, $BG\in\Psib^0(X)$ and indeed in a bounded subset
of $\Psibc^0(X)$ for $B\in\cB$, so Lemma~\ref{lemma:H10-bd}
also gives that for some $C_2>0$ (independent of $B\in\cB$),
$\|BGQu\|_{H^1(X)}\leq C_2\|Qu\|_{H^1(X)}$. Combining these proves the
lemma.
\end{proof}

We can similarly microlocalize $H^{-1}_{\loc}(X)$:

\begin{Def}\label{Def:WFbd}
Suppose $u\in H^{-1}_{\loc}(X)$, $m\geq 0$.
We say that $q\in\Tb^*X\setminus o$ is not in $\WFbd^{-1,m}(u)$ if
there exists $A\in\Psib^m(X)$ such that $\sigma_{b,m}(A)(q)\neq 0$
and $Au\in H^{-1}(X)$.
\end{Def}

Then the analogues of Lemma~\ref{lemma:WFb-mic}-\ref{lemma:WFb-mic-q}
remain valid with $H^1(X)$ replaced by $H^{-1}(X)$ and
$\WFbz^{1,\cdot}$ replaced by $\WFbd^{-1,\cdot}$, with analogous proofs using
Corollary~\ref{cor:H-1-bd} in place of Lemma~\ref{lemma:H10-bd}.

These results can be extended in another way, by considering Sobolev
spaces with a negative order of regularity relative to $H^1(X)$.

\begin{Def}\label{Def:H1m-neg}
Let $k$ be an integer, $m<0$, and $A\in\Psib^{-m}(X)$ be elliptic on $\Sb^*X$
with proper support.
We let $H^{k,m}_{\bl,c}(X)$ be the space of all $u\in\dist(X)$ of the
form $u=u_1+Au_2$ with $u_1,u_2\in H^k_c(X)$.
We let
\begin{equation*}
\|u\|_{H^{k,m}_{\bl,c}(X)}
=\inf\{\|u_1\|_{H^k(X)}+\|u_2\|_{H^k(X)}:\ u=u_1+Au_2\}.
\end{equation*}

We also let $H^{k,m}_{\bl,\loc}(X)$ be the space of all $u\in\dist(X)$
such that $\phi u\in H^{k,m}_{\bl,c}(X)$ for all $\phi\in\Cinf_c(X)$.

We also define $\dot H^{k,m}_{\bl,c}(X)$ and $\dot H^{k,m}_{\bl,\loc}(X)$
analogously, replacing $H^k(X)$ by $\dot H^k(X)$ throughout the
above discussion. Here, for $k\geq 0$, $\dot H^k(X)$ stands for
$H^k_0(X)$, see Remark~\ref{rem:Sobolev}, so we also write
$\dot H^{k,m}_{\bl,c}(X)=H^{k,m}_{\bl,0,c}(X)$ for $k\geq 0$.
\end{Def}

\begin{rem}
In this paper we are only concerned with the cases $k=\pm 1$. There
is no difference between these two cases for the ensuing discussion,
except for the boundary values considered in the next paragraph.
For the sake of definiteness, we will use
$k=1$ throughout the discussion. We will also not consider
$\dot H^k(X)$ explicitly for most of the discussion; there is no
difference for the treatment of these spaces either.

We also remark that we can talk about the boundary values of
$u\in H^{1,m}_{\bl,c}(X)$ at boundary hypersurfaces $H_j$ for $m<0$, although
we do not need this here.
One way to do this is
to define, for $u=u_1+Au_2$, $u|_{H_j}=u_1|_{H_j}+\hat N_j(A)(0)(u_2|_{H_j})$,
regarded e.g.\ as an element of $\dist(H_j)$ (note that $\hat N_j(A)(0):
\dist(H_j)\to\dist(H_j)$), and this is independent of the choices of $u_1$,
$u_2$ and $A$. Of course, for $u\in H^{1,m}_{\bl,0,c}(X)$, in the
sense just sketched, $u|_{H_j}=0$
for all $j$.
It is straightforward to see that for
$u\in H^{1,m}_{\bl,c}$ with $u|_{H_j}=0$ for all $j$, there exist
$u_1,u_2\in H^1_{0,c}(X)$ with $u=u_1+Au_2$, so $u\in H^{1,m}_{\bl,0,c}(X)$.

We also remark that
Lemma~\ref{lemma:A-H1-H10} still holds if one only assumes
$u\in H^{1,m}_{\bl,0,c}(X)$.
\end{rem}

First note that given any $K\subset X$ compact there is another
$K'\subset X$ compact such that $u\in H^{1,m}_{\bl,c}(X)$ with
$\supp u\subset K$ can be written as $u=u_1+Au_2$ with $u_1,u_2\in H^1_c(X)$
both supported in $K'$. Indeed, let $\phi\in\Cinf_c(X)$ be identically $1$
on a neighborhood of $K$, and let $G\in\Psib^{m}(X)$ be a properly
supported parametrix for $A$, so $AG=\Id+E$, $E\in\Psib^{-\infty}(X)$,
$E$ also properly supported. By definition, if $u\in H^{1,m}_{\bl,c}(X)$
then there are $u_1',u_2'\in H^1_c(X)$ with $u=u_1'+Au_2'$, and as
$\phi\equiv 1$ on a neighborhood of $\supp u$, $\phi u=u$.
Thus,
\begin{equation*}\begin{split}
&u=\phi u=\phi u_1'-E\phi Au_2'+AG\phi Au_2'=u_1+u_2,\\
&\qquad u_1=\phi u_1'-E\phi A u_2',\ u_2=G\phi Au_2',
\end{split}\end{equation*}
so $u_1,u_2\in H^1_c(X)$
as $E\phi A,G\phi A\in\Psib^0(X)$, and $\supp u_j$, $j=1,2$,
is bounded in terms of $\supp\phi$, $\supp E$ and $\supp G$. Namely,
\begin{equation*}\begin{split}
&\supp u_j\subset K',\\
& K'=\supp\phi\cup\pi_L(\supp E\cap\pi_R^{-1}(\supp\phi))
\cup\pi_L(\supp G\cap\pi_R^{-1}(\supp\phi)),
\end{split}\end{equation*}
where $\pi_L,\pi_R:X\times X\to X$ are the projections to the left
and right factors; $K'$ is compact as $E$ and $G$ are
properly supported, so $\supp E\cap\pi_R^{-1}(\supp\phi)$,
$\supp G\cap\pi_R^{-1}(\supp\phi)$ are compact.
Note also that, by Lemma~\ref{lemma:H10-bd},
$\|u_1\|_{H^1(X)}+\|u_2\|_{H^1(X)}\leq C(\|u_1'\|_{H^1(X)}+\|u_2'\|_{H^1(X)})$.
Since this holds for any $u_1'$, $u_2'$ with $u=u_1'+Au_2'$, we deduce
that with this $K'$, if we restrict $\supp u_j\subset K'$, and take
$\inf$ just over these $u_j$, we get an equivalent norm on the
subspace of $H^1_c(X)$ consisting of elements supported in $K$.

In fact, as $\supp G$, $\supp E$ can be made to lie in any neighborhood
of the diagonal in $X\times X$, and $\supp\phi$ can be made to lie
in any neighborhood of $K$, this argument shows that given any $K$ compact
and any $U$ open with $K\subset U$, $\supp u_j$ may be assumed to lie
in $K'=\overline{U}$, with the resulting norm equivalent to the $H^1_c(X)$
norm of the definition (with the equivalence constant of course depending
on $U$!).

Moreover, Definition~\ref{Def:H1m-neg}
is independent of the choice of $A$. Indeed, if
$A'\in\Psib^{-m}(X)$ is elliptic and has proper support, then it
has a parametrix $G'\in\Psib^{m}(X)$ with $E'=A'G'-\Id\in\Psib^{-\infty}(X)$,
all with proper support. Then $u=u_1+Au_2=u_1-E'Au_2+A'G'Au_2$, and
$u_1'=u_1-E'Au_2\in H^1_c(X)$ since $E'A\in\Psib^{-\infty}(X)$, and
$u_2'=G'A u_2\in H^1_c(X)$ since $G'A\in\Psib^0(X)$.
Moreover, if we fix $K\subset X$ compact, then for $u$ with $\supp u\subset K$,
the norms $\|u\|_{H^{1,m}_{\bl,c}(X)}$ are equivalent for different
choices of $A$ -- this follows from Lemma~\ref{lemma:H10-bd} and the
preceeding remark that we may take the support of $u_1$, $u_2$ lie
in a compact set depending on $K$ only.

Note also that for $F\in\Psibc^{m}(X)$ with compactly supported
Schwartz kernel, $F:H^{1,m}_{\bl,c}(X)\to H^1(X)$ is continuous. Indeed,
$Fu=Fu_1+FAu_2\in H^1_c(X)$ by Lemma~\ref{lemma:H10-bd}
since $F,FA\in\Psibc^{0}(X)$ and $u_1,u_2\in H^1_c(X)$, and this
also gives a bound for $\|Fu\|_{H^1(X)}$ in terms of
$\|u\|_{H^{1,m}_{\bl,c}(X)}$ and a seminorm of $F$ in $\Psibc^{m}(X)$.
In particular, $\Psib^{-\infty}(X)$
maps $H^{1,m}_{\bl,c}(X)\to H^1(X)$, and indeed into the conormal space
$H^{1,\infty}_{\bl,c}(X)$.

Since any $A\in\Psib^m(X)$ defines a map $A:\dist(X)\to\dist(X)$, our
definition of the wave front set makes sense for $m<0$ as well; it
is independent of $s$ if we take $u\in H^{1,s}_{\loc}(X)$ since
the action of $\Psib(X)$ is well-defined on the larger space
$\dist(X)$ already.

\begin{Def}\label{Def:WFb-neg}
Suppose $u\in H^{1,s}_{\loc}(X)$ for some $s\leq
0$, and suppose that $m\in\Real$.
We say that $q\in\Tb^*X\setminus o$ is not in $\WFbz^{1,m}(u)$ if
there exists $A\in\Psib^m(X)$ such that $\sigma_{b,m}(A)(q)\neq 0$
and $Au\in H^1(X)$.

For $m=\infty$, we say that $q\in\Tb^*X\setminus o$ is not in $\WFbz^{1,m}(u)$
if there exists $A\in\Psib^0(X)$ such that $\sigma_{b,0}(A)(q)\neq 0$
and $LAu\in H^1(X)$ for all $L\in\Diffb(X)$, i.e.\ if
$Au\in H^{1,\infty}_b(X)$.
\end{Def}

Again, the analogues of Lemma~\ref{lemma:WFb-mic}-\ref{lemma:WFb-mic-q}
remain valid with $H^1(X)$ replaced by $H^{1,s}_{\bl,c}(X)$ for some $s$, and
$m$ allowed to be negative in $\WFbz^{1,m}(u)$. In particular,
Lemma~\ref{lemma:WFb-mic-q} takes the form:

\begin{lemma}\label{lemma:WFb-mic-q-p}
Suppose that $K\subset\Sb^*X$ is compact, and $U$ a neighborhood of $K$
in $\Sb^*X$. Let $\tilde K\subset X$ compact,
and $\tilde U$ is a neighborhood of $\tilde K$ in $X$ with compact closure.
Let $Q\in\Psib^k(X)$ be elliptic on $K$ with $\WFb'(Q)
\subset U$, with Schwartz kernel supported in $\tilde K\times\tilde K$.
Let $\cB$ be a bounded subset of $\Psibc^k(X)$
with $\WFb'(\cB)\subset K$ and
Schwartz kernel supported in $\tilde K\times\tilde K$.
Then for any $s<0$
there is a constant $C>0$ such that for $B\in\cB$, $u\in H^{1,s}_{\bl,\loc}(X)$
with $\WFbz^{1,k}(u)\cap U=\emptyset$,
\begin{equation*}
\|Bu\|_{H^1(X)}\leq C(\|u\|_{H^{1,s}_{\bl}(\tilde U)}+\|Qu\|_{H^1(X)}),
\end{equation*}
where $\|u\|_{H^{1,s}_{\bl}(\tilde U)}$ stands for
$\|\phi u\|_{H^{1,s}_{\bl,c}(X)}$ for some fixed $\phi\in\Cinf_c(X)$
with $\supp\phi\subset \tilde U$, $\phi\equiv 1$ on a neighborhood
of $\tilde K$.
\end{lemma}

Finally, connecting $H^{k,m}_{\bl,\loc}(X)$ for $k=\pm 1$, we remark
that any $P\in\Diff^2_\bl(X)$ defines a continuous
linear map $P:H^{1,m}_{\bl,\loc}(X)\to H^{-1,m}_{\bl,\loc}(X)$,
as discussed before the statement of Corollary~\ref{cor:H-1-bd};
now we need to use \eqref{eq:V-Psib-1} as well to deduce this.

\section{The elliptic set}\label{sec:elliptic}
We first prove an estimate that microlocally
controls the Dirichlet form for microlocalized solutions
$Pu=0$, $u\in H^1_0(X)$, in terms of a {\em lower order} microlocal
information and a global bound in $H^1_0(X)$. In fact, as it does not
require much additional effort, we consider microlocal solutions,
i.e.\ we make assumptions on $\WFbd^{-1,\infty}(Pu)$, or
indeed $\WFbd^{-1,s}(Pu)$.

\begin{rem}\label{rem:localize}
Since $X$ is non-compact and our results are microlocal, we may always
fix a compact set $\tilde K\subset X$ and assume that all ps.d.o's
have Schwartz kernel supported in $\tilde K\times\tilde K$. We
also let $\tilde U$ be a neighborhood of $\tilde K$ in $X$ such that
$\tilde U$ has compact closure, and use the $H^1(\tilde U)$ norm
in place of the $H^1(X)$ norm to accommodate $u\in H^1_{0,\loc}(X)$.
Below we use the notation $\|.\|_{H^1_{\loc}(X)}$ for
$\|.\|_{H^1(\tilde U)}$ to avoid having to specify $\tilde U$.
We also use $\|.\|_{H^{-1}_{\loc}(X)}$ for
$\|.\|_{H^1(\tilde U)}$.
\end{rem}

We give two versions of the Dirichlet estimates: the first one
suffices for most purposes, but it does not give the optimal
estimates in terms of the order $m$ in $\WFbd^{-1,m}(Pu)$. The second one
takes care of this issue.

\begin{lemma}\label{lemma:Dirichlet-form}
Suppose that $K\subset\Sb^*X$ is compact, $U\subset\Sb^*X$ is open,
$K\subset U$.
Suppose that $\calA=\{A_r:\ r\in(0,1]\}$ be a bounded
family of ps.d.o's in $\Psibc^s(X)$ with $\WFb'(\calA)\subset K$, and
with $A_r\in\Psib^{s-1}(X)$ for $r\in (0,1]$.
Then there are $G\in\Psib^{s-1/2}(X)$, $\tilde G\in\Psib^{s+1/2}(X)$
with $\WFb'(G),\WFb'(\tilde G)\subset U$
and $C_0>0$ such that for $r\in(0,1]$, $u\in H^1_{0,\loc}(X)$
with $\WFbz^{1,s-1/2}(u)
\cap U=\emptyset$, $\WFbd^{-1,s+1/2}(Pu)\cap U=\emptyset$, we have
\begin{equation*}\begin{split}
&|\int_X \left(|d_M A_r u|^2-|D_t A_r u|^2\right)|\\
&\qquad\leq
C_0(\|u\|^2_{H^1_{\loc}(X)}+\|Gu\|^2_{H^1(X)}+\|Pu\|^2_{H^{-1}_{\loc}(X)}
+\|\tilde G Pu\|^2_{H^{-1}(X)}).
\end{split}\end{equation*}
In particular, if the assumption on $Pu$ is strengthened to $Pu=0$, we have
\begin{equation*}
|\int_X \left(|d_M A_r u|^2-|D_t A_r u|^2\right)|\leq
C_0(\|u\|^2_{H^1_{\loc}(X)}+\|Gu\|^2_{H^1(X)}).
\end{equation*}
The meaning of $\|u\|^2_{H^1_{\loc}(X)}$ and
$\|Pu\|^2_{H^{-1}_{\loc}(X)}$ is stated above in
Remark~\ref{rem:localize}.
\end{lemma}

\begin{rem}
The point of this lemma is $G$ is $1/2$ order lower ($s-1/2$ vs. $s$)
than the {\em family} $\calA$. We will later take a limit, $r\to 0$,
which gives control of the Dirichlet form evaluated on $A_0u$,
$A_0\in\Psibc^s(X)$, in terms of lower order information.

The role of $A_r$, $r>0$, is to regularize such an argument, i.e.\ to make
sure various terms in a formal computation, in which one uses $A_0$
directly, actually make sense.
\end{rem}

\begin{proof}
Then for $r\in(0,1]$, $A_r u\in H^1_0(X)$, so
\begin{equation*}
\int_X (|d_M A_r u|^2-|D_t A_r u|^2)=-\int_X PA_r u \,\overline{A_r u}.
\end{equation*}
Here the right hand side is the pairing of $H^{-1}(X)$ with $H^1_0(X)$.
Writing $PA_r=A_rP+[P,A_r]$, we see that
the right hand side can be estimated by
\begin{equation}\label{eq:Dirichlet-form-8}
|\int_X A_rPu \,\overline{A_r u}|+|\int_X [P,A_r] u \,\overline{A_r u}|.
\end{equation}
The lemma is thus proved if we show that the first term of
\eqref{eq:Dirichlet-form-8} is bounded by
\begin{equation}\label{eq:Dirichlet-bd}
C_0'(\|u\|^2_{H^1_{\loc}(X)}+\|Gu\|^2_{H^1(X)}+\|Pu\|^2_{H^{-1}_{\loc}(X)}
+\|\tilde G Pu\|^2_{H^{-1}(X)}),
\end{equation}
the second term
is bounded by $C_0''(\|u\|^2_{H^1_{\loc}(X)}+\|Gu\|^2_{H^1(X)})$.

The first term is straightforward to estimate.
Let $\Lambda\in\Psib^{-1/2}(X)$ be elliptic with
$\Lambda^-\in\Psib^{1/2}(X)$ a parametrix, so
\begin{equation*}
E=\Lambda\Lambda^--\Id,E'=\Lambda^-\Lambda-\Id\in\Psib^{-\infty}(X).
\end{equation*}
Then
\begin{equation*}\begin{split}
\int_X A_rPu \,\overline{A_r u}&=\int_X (\Lambda\Lambda^--E)A_r Pu
\,\overline{A_r Pu}\\
&=\int_X \Lambda^-A_r Pu\,\overline{\Lambda^* A_r Pu}
-\int_X A_r Pu\,\overline{E^*A_r u}.
\end{split}\end{equation*}
Since $\Lambda^- A_r$ is uniformly bounded in $\Psibc^{s+1/2}(X)$,
and $\Lambda^* A_r$ is uniformly bounded in $\Psibc^{s-1/2}(X)$,
$\int_X \Lambda^-A_r Pu\,\overline{\Lambda^* A_r Pu}$
is uniformly bounded, with a bound like \eqref{eq:Dirichlet-bd}
using Cauchy-Schwartz and Lemma~\ref{lemma:WFb-mic-q}.
Indeed, by Lemma~\ref{lemma:WFb-mic-q},
choosing any $G\in\Psib^{s-1/2}(X)$ which is elliptic on
$K$, there is a constant $C_1>0$ such that
\begin{equation*}
\|\Lambda^*A_r u\|^2_{H^1(X)}\leq
C_1(\|u\|^2_{H^1_{\loc}(X)}+\|Gu\|^2_{H^1(X)}).
\end{equation*}
Similarly, by Lemma~\ref{lemma:WFb-mic-q} and the remark following
Definition~\ref{Def:WFbd},
choosing any $\tilde G\in\Psib^{s+1/2}(X)$ which is elliptic on
$K$, there is a constant $C_1'>0$ such that
$\|\Lambda^- A_rPu\|^2_{H^{-1}(X)}\leq
C_1'(\|Pu\|^2_{H^{-1}_{\loc}(X)}+\|\tilde GPu\|^2_{H^{-1}(X)})$.
Combining these gives, with $C_0'=C_1+C_1'$,
\begin{equation*}\begin{split}
|\int_X &\Lambda^-A_r Pu\,\overline{\Lambda^* A_r Pu}|\leq\|\Lambda^-A_r Pu\|
\,\|\Lambda^*A_r u\|\leq \|\Lambda^-A_r Pu\|^2+\|\Lambda^*A_r u\|^2\\
&\leq C_0'(\|u\|^2_{H^1_{\loc}(X)}+\|Gu\|^2_{H^1(X)}
+\|Pu\|^2_{H^{-1}_{\loc}(X)}+\|\tilde GPu\|^2_{H^{-1}(X)}),
\end{split}\end{equation*}
as desired.

A similar argument, using that
$A_r$ is uniformly bounded in $\Psibc^{s+1/2}(X)$ (in fact in
$\Psibc^{s}(X)$),
and $E^* A_r$ is uniformly bounded in $\Psibc^{s-1/2}(X)$
(in fact in $\Psibc^{-\infty}(X)$), shows that
$\int_X A_r Pu\,\overline{E^* A_r Pu}$
is uniformly bounded.

Now we turn to the second term in \eqref{eq:Dirichlet-form-8}.
Using \eqref{eq:V-Psib-1},
\begin{equation*}
[P,A_r]=\sum_{i,j}D_{x_i}D_{x_j}B_{ij,r}+\sum D_{x_j} B_{j,r}+B_r,
\end{equation*}
$B_r\in\Psib^{s-1}(X)$,
$B_{j,r}\in\Psib^{s-2}(X)$, $B_{ij,r}\in\Psib^{s-3}(X)$,
uniformly bounded in $\Psibc^{s+1}(X)$, resp.\ $\Psibc^{s}(X)$,
resp.\ $\Psibc^{s-1}(X)$.
With $\Lambda\in\Psib^{-1/2}(X)$ as above,
we can write further
\begin{equation*}
\Lambda^{-}[P,A_r]=\sum_{i,j}D_{x_i}D_{x_j}\Lambda^{-}B'_{ij,r}
+\sum D_{x_j}\Lambda^{-} B'_{j,r}+\Lambda^{-}B'_r,
\end{equation*}
with $B'_{ij,r}$, $B'_{j,r}$, $B'_r$ having the same properties as
the $B_{ij,r}$, etc., listed above. Thus,
\begin{equation}\begin{split}\label{eq:ell-16}
\int_X [P,A_r] u \,\overline{A_r u}
&=\sum_{ij}\int_X D_{x_i}D_{x_j}\Lambda^{-}B'_{ij,r}u
\,\overline{\Lambda^*A_r u}
-\sum_{ij}\int_X D_{x_i}D_{x_j} E B'_{ij,r}u
\,\overline{A_r u}\\
&\quad+\sum_j \int_X D_{x_j}\Lambda^{-} B'_{j,r} u\,\overline{\Lambda^*A_r u}
-\sum_j \int_X D_{x_j}E B'_{j,r} u\,\overline{A_r u}\\
&\quad+\int_X \Lambda^{-}B'_r u\,\overline{\Lambda^*A_r u}
-\int_X E B'_r u\,\overline{A_r u}\\
&=\sum_{ij}\int_X D_{x_j}\Lambda^{-}B'_{ij,r}u
\,\overline{D_{x_i}^t\Lambda^*A_r u}
-\sum_{ij}\int_X D_{x_j} E B'_{ij,r}u
\,\overline{D_{x_i}^t A_r u}\\
&\quad+\sum_j \int_X D_{x_j}\Lambda^{-} B'_{j,r} u\,\overline{\Lambda^*A_r u}
-\sum_j \int_X D_{x_j} E B'_{j,r} u\,\overline{A_r u}\\
&\quad
+\int_X \Lambda^{-}B'_r u\,\overline{\Lambda^*A_r u}
-\int_X E B'_r u\,\overline{A_r u},
\end{split}\end{equation}
where $D_{x_i}^t$ is the formal adjoint of $D_{x_i}$ with respect to $dg$, and
where in the last step we used that
\begin{equation*}
\Lambda^{-}B'_{ij,r}u,\Lambda^*A_r u, EB'_{ij,r}u,A_ru\in H^1_0(X).
\end{equation*}
Note that
$D_{x_i}^t=J ^{-1}D_{x_i}J$ if $dg=J dx_1\ldots dx_k\, dy_1\ldots dy_l$ is
the Riemannian density, so $D_{x_i}^t=D_{x_i}+b$, $b\in\Cinf(X)$.
Thus,
\begin{equation*}\begin{split}
|\int_X D_{x_j}\Lambda^{-}B'_{ij,r}u
\,\overline{D_{x_i}^t\Lambda^*A_r u}|
&\leq \|D_{x_j}\Lambda^{-}B'_{ij,r}u\|_{L^2(X)}
\|D_{x_i}\Lambda^*A_r u\|_{L^2(X)}\\
&\qquad+\|D_{x_j}\Lambda^{-}B'_{ij,r}u\|_{L^2(X)}
\|\Lambda^*A_r u\|_{L^2(X)},
\end{split}\end{equation*}
and both factors in both terms
are uniformly bounded for $r\in(0,1]$ since $\Lambda^* A_r$,
$\Lambda^{-}B'_{ij,r}$
are uniformly bounded in $\Psibc^{s-1/2}(X)$ with a uniform wave front bound
disjoint from $\WFbz^{1,s-1/2}(u)$.
Indeed, as noted above, by Lemma~\ref{lemma:WFb-mic-q},
choosing any $G\in\Psib^{s-1/2}(X)$ which is elliptic on
$K$, there is a constant $C_1>0$ such that
the right hand side is bounded by
$C_1(\|u\|^2_{H^1_{\loc}(X)}+\|Gu\|^2_{H^1(X)})$.
Similar estimates apply to the other
terms on the right hand side of \eqref{eq:ell-16}, showing that
$\int_X [P,A_r] u \,\overline{A_r u}$ is uniformly bounded for $r\in(0,1]$,
indeed is bounded by $C_0(\|u\|^2_{H^1_{\loc}(X)}+\|Gu\|^2_{H^1(X)})$,
proving the lemma.
\end{proof}

The lemma which allows more precise estimates is the following.

\begin{lemma}\label{lemma:Dirichlet-form-2}
Suppose that $K\subset\Sb^*X$ is compact, $U\subset\Sb^*X$ is open,
$K\subset U$.
Suppose that $\calA=\{A_r:\ r\in(0,1]\}$ be a bounded
family of ps.d.o's in $\Psibc^s(X)$ with $\WFb'(\calA)\subset K$, and
with $A_r\in\Psib^{s-1}(X)$ for $r\in (0,1]$.
Then there are $G\in\Psib^{s-1/2}(X)$, $\tilde G\in\Psib^{s}(X)$
with $\WFb'(G),\WFb'(\tilde G)\subset U$
and $C_0>0$ such that for $\ep>0$, $r\in(0,1]$, $u\in H^1_{0,\loc}(X)$
with $\WFbz^{1,s-1/2}(u)
\cap U=\emptyset$, $\WFbd^{-1,s}(Pu)\cap U=\emptyset$, we have
\begin{equation*}\begin{split}
|\int_X \left(|d_M A_r u|^2-|D_t A_r u|^2\right)|\leq
\ep&\|d_X A_r u\|^2_{L^2(X)}
+C_0(\|u\|^2_{H^1_{\loc}(X)}+\|Gu\|^2_{H^1(X)}\\
&+\ep^{-1}\|Pu\|^2_{H^{-1}_{\loc}(X)}
+\ep^{-1}\|\tilde G Pu\|^2_{H^{-1}(X)}).
\end{split}\end{equation*}
\end{lemma}

\begin{rem}
The point of this lemma is that on the one hand
the new term $\ep\|d_X A_r u\|^2$ can be absorbed
in the left hand side in the elliptic region, hence is negligible,
on the other hand, there is a gain in the
order of $\tilde G$ ($s$, versus $s+1/2$ in the previous lemma).
\end{rem}

\begin{proof}
We only need to modify the previous proof slightly. Thus, we need to
estimate the term $|\int_X A_rPu \,\overline{A_r u}|$ in
\eqref{eq:Dirichlet-form-8} differently, namely
\begin{equation*}
|\int_X A_rPu \,\overline{A_r u}|\leq \|A_r Pu\|_{H^{-1}(X)}
\|A_r u\|_{H^1(X)}\leq \ep\|A_r u\|^2_{H^1(X)}+\ep^{-1}
\|A_r Pu\|^2_{H^{-1}(X)}.
\end{equation*}
Now the lemma follows by using
Lemma~\ref{lemma:WFb-mic-q} and the remark following
Definition~\ref{Def:WFbd}, namely
choosing any $\tilde G\in\Psib^{s}(X)$ which is elliptic on
$K$, there is a constant $C_1'>0$ such that
$\|A_rPu\|^2_{H^{-1}(X)}\leq
C_1'(\|Pu\|^2_{H^{-1}_{\loc}(X)}+\|\tilde GPu\|^2_{H^{-1}(X)})$,
and finishing the proof exactly as for Lemma~\ref{lemma:Dirichlet-form}.
\end{proof}

Using the microlocal positivity of the Dirichlet form, we now
prove the elliptic estimates.

\begin{prop}\label{prop:elliptic}(Microlocal elliptic regularity.)
If $u\in H^1_{0,\loc}(X)$ then
\begin{equation*}
\WFbz^{1,m}(u)\subset \WFbd^{-1,m}(Pu)\cup\dot\Tb^*X,\Mand
\WFbz^{1,m}(u)\cap\cE\subset \WFbd^{-1,m}(Pu).
\end{equation*}

In particular, if $Pu=0$, $u\in H^1_{0,\loc}(X)$ then
\begin{equation*}
\WFbz^{1,\infty}(u)\subset \dot\Tb^*X,\Mand \WFbz^{1,\infty}(u)\cap\cE=\emptyset.
\end{equation*}
\end{prop}

\begin{proof}
We first prove a slightly weaker result in which
$\WFbd^{-1,m}(Pu)$ is replaced by $\WFbd^{-1,m+1/2}(Pu)$ -- we rely on
Lemma~\ref{lemma:Dirichlet-form}. We then prove the
original statement using Lemma~\ref{lemma:Dirichlet-form-2}.

Suppose that either $q\in \Tb^*X\setminus\dot\Tb^*X$ or
$q\in\cE$. We may assume iteratively that $q\nin\WFbz^{1,s-1/2}(u)$;
we need to prove then that $q\nin\WFbz^{1,s}(u)$ (note that the inductive
hypothesis holds for $s=1/2$ since $u\in H^1_{\loc}(X)$).
Let $A\in\Psib^{s}(X)$ be such that
$\WFb'(A)\cap \WFbz^{1,s-1/2}(u)=\emptyset$, $\WFb'(A)\cap\WFbz^{1,s+1/2}(Pu)
=\emptyset$, and have $\WFb'(A)$ in a small
conic neighborhood $U$ of $q$ so that for a suitable $C>0$ or $\ep>0$, in $U$
\begin{enumerate}
\item
$\tau^2<C \sum_j \sigma_j^2$ if $q\in \Tb^*X\setminus\dot\Tb^*X$,
\item
$|\sigma_j|<\ep(\tau^2+|\zeta|^2)^{1/2}$ for all $j$, and
$\frac{|\zeta|}{|\tau|}>1+\ep$, if $q\in\cE$.
\end{enumerate}
Let $\Lambda_r\in\Psib^{-2}(X)$ for $r>0$, such that $\calL=\{\Lambda_r:
\ r\in(0,1]\}$ is a bounded family in $\Psib^0(X)$, and $\Lambda_r\to\Id$
as $r\to 0$ in $\Psib^{\tilde\ep}(X)$, $\tilde\ep>0$,
e.g.\ the symbol of $\Lambda_r$ could be taken as
$(1+r(\tau^2+|\zeta|^2+|\sigma|^2))^{-1}$. Let $A_r=\Lambda_r A$.
Let $a$ be the symbol of $A$, and
let $A_r$ have symbol $(1+r(\tau^2+|\zeta|^2+|\sigma|^2))^{-1}a$, $r>0$,
so $A_r\in
\Psib^{s-2}(X)$ for $r>0$, and $A_r$ is uniformly bounded in
$\Psibc^{s}(X)$, $A_r\to A$ in $\Psibc^{s+\tilde\ep}(X)$.

By Lemma~\ref{lemma:Dirichlet-form},
\begin{equation*}
\int_X \left(|d_M A_r u|^2-|D_t A_r u|^2\right)
\end{equation*}
is uniformly bounded for $r\in(0,1]$.
On the other hand,
\begin{equation*}\begin{split}
\int_X |d_M A_r u|^2
=\int_X \sum A_{ij} D_{x_i}A_r u \,\overline{D_{x_j} A_r u}
&+\int_X \sum B_{ij} D_{y_i}A_r u \,\overline{D_{y_j} A_r u}\\
&+\int_X \sum C_{ij} D_{x_i}A_r u \,\overline{D_{y_j} A_r u}.
\end{split}\end{equation*}
Using that $A_{ij}(x,y)=A_{ij}(0,y)+\sum x_k A'_{ijk}(x,y)$, we see that
if $A_r$ is supported in $x_k<\delta$ for all $k$,
\begin{equation}\label{eq:elliptic-28}
|\int_X \sum x_k A'_{ijk} D_{x_i}A_r u \,\overline{D_{x_j} A_r u}|
\leq C\delta\sum_{i',j'}\|D_{x_{i'}}A_r u\|\,\|D_{x_{j'}}A_r u\|,
\end{equation}
with analogous estimates for $B_{ij}(x,y)-B_{ij}(0,y)$ and
for $C_{ij}(x,y)$. Moreover, as the matrix $A_{ij}$ is positive definite,
for some $c>0$,
\begin{equation*}
c\int_X \sum_j |D_{x_j}A_r u|^2\leq
\frac{1}{2}\int_X \sum_{ij} A_{ij} D_{x_i}A_r u \overline{D_{x_j} A_r u}.
\end{equation*}
Thus,
there exists $\tilde C>0$ and $\delta_0>0$ such that if $\delta<\delta_0$
and $A$ is supported
in $|x|<\delta$ then
\begin{equation}\begin{split}\label{eq:ell-32}
&c\int_X \sum_j|D_{x_j} A_r u|^2
+\int_X((1-\tilde C\delta)\sum_j |D_{y_j} A_r u|_h^2-|D_t A_r u|^2)\\
&\qquad\leq \int_X (|d_M A_r u|^2-|D_t A_r u|^2),
\end{split}\end{equation}
where we used the notation
\begin{equation*}
|D_{y_j} A_r u|_h^2=\sum_{ij} B_{ij}(0,y) D_{y_i}A_r u\,\overline{D_{y_j}A_r u},\end{equation*}
i.e.\ $h$ is the metric $g$ restricted to the span of the $\pa_{y_j}$,
$j=1,\ldots,l$.

Now we distinguish the cases $q\in\cE$ and $q\in\Tb^*X\setminus\dot\Tb^*X$.
If $q\in\cE$, $A$ is supported near $\cE$,
we choose $\delta\in(0,\frac{1}{2\tilde C})$ so that
$(1-\tilde C\delta)\frac{|\zeta|^2}{\tau^2}>1+\delta$ on
a neighborhood of $\WFb'(A)$, which is possible in view of (ii) at the
beginning of the proof.
Then the second integral on the left hand side of \eqref{eq:ell-32}
can be written
as $\| BA_r u\|^2$, with the symbol of $B$ given by
$((1-\tilde C\delta)|\zeta|^2
-\tau^2)^{1/2}$ )(which is $\geq\delta\tau$), modulo a term
\begin{equation*}
\int_X F A_r u\,\overline{A_r u},\ F\in\Psib^{1}(X).
\end{equation*}
But this expression is uniformly bounded as $r\to 0$ by the
argument above. We thus deduce that
\begin{equation*}
c\int_X (\sum_j|D_{x_j} A_r u|^2)+\| BA_r u\|^2
\end{equation*}
is uniformly bounded as $r\to 0$.

If $q\in\Tb^*X\setminus\dot\Tb^*X$, and $A$ is supported in $|x|<\delta$,
\begin{equation*}
\int_X \delta^{-2}|x_j D_{x_j} A_r u|^2\leq \int_X |D_{x_j} A_r u|^2,
\end{equation*}
On the other hand, near $\Tb^*X\setminus\dot\Tb^*X$,
for $\delta>0$ sufficiently small,
\begin{equation*}
\int_X\left(\frac{c}{2\delta^{2}}
\sum_j |x_j D_{x_j}A_r u|^2-|D_t A_r u|^2\right)=\|BA_r u\|^2
+\int_X F A_r u\,\overline{A_r u},
\end{equation*}
with the symbol of $B$ given by
$(\frac{c}{2\delta^{2}}\sum \sigma_j^2-\tau^2)^{1/2}$
(which does not vanish on $U$ for
$\delta>0$ small), while
$F\in\Psib^{1}(X)$, so the second term on the right hand side
is uniformly bounded as $r\to 0$.
We thus deduce in this case that
\begin{equation*}
\frac{c}{2}\int_X (\sum_j|D_{x_j} A_r u|^2)+\| BA_r u\|^2
\end{equation*}
is uniformly bounded as $r\to 0$.

We thus conclude that
$D_{x_j} A_ru, BA_ru$ are uniformly bounded $L^2(X)$. Correspondingly
there are sequences $D_{x_j}A_{r_k}u$, $BA_{r_k}u$,
weakly convergent in $L^2(X)$, and such that $r_k\to 0$, as $k\to\infty$.
Since they converge to $D_{x_j}Au$, $BAu$, respectively, in $\dist(X)$,
we deduce that the weak limits are $D_{x_j}Au$, $BAu$, which therefore
lie in $L^2(X)$. Consequently,
$d Au\in L^2(X)$ proving the
proposition with $\WFbd^{-1,m}(Pu)$ replaced by $\WFbd^{-1,m+1/2}(Pu)$.

To obtain the optimal result, we note that due to
Lemma~\ref{lemma:Dirichlet-form-2} we still have, for any $\ep>0$,
that
\begin{equation*}\begin{split}
&\int_X \left(|d_M A_r u|^2-|D_t A_r u|^2-\ep|d_X A_r u|^2\right)\\
&\qquad=\int_X \left((1-\ep)|d_M A_r u|^2-(1+\ep)|D_t A_r u|^2\right)
\end{split}\end{equation*}
is uniformly bounded above for $r\in(0,1]$.
By arguing just as above, with $B$ as above,
for sufficiently small $\ep>0$, the right hand side gives an upper bound for
\begin{equation*}
\frac{c}{2}\int_X (\sum_j|D_{x_j} A_r u|^2)+\| BA_r u\|^2,
\end{equation*}
which is thus uniformly bounded as $r\to 0$. The proof is then finished
exactly as above.
\end{proof}

A slightly different formulation of this argument is the following. Below
$w=(x,y)$.
Consider
\begin{equation*}\begin{split}
\|d_M A_r u\|^2-&\|D_t A_r u\|^2\\
&=\int_X \sum_{i,j}g^{ij}D_{w_i} A_r u
\,\overline{D_{w_j}A_r u} \,J \,dw\,dt
-\int_X D_t A_r u\,\overline{D_t A_ru} \,J\,dw\,dt.
\end{split}\end{equation*}
We move the $A_r$ in the first factor of each term on the right hand side
by first commuting it through $g^{ij}D_{w_i}$ (or $D_t$), then
taking its adjoint with respect to $J\,dw\,dt$, and finally commuting
it through $D_{w_j}$. Each of the commutator terms can be controlled by
the inductive hypothesis as above. Modulo such terms the result is
\begin{equation}\label{eq:DN-16}
\int_X \left(\sum_{i,j}g^{ij}D_{w_i} u\,\overline{D_{w_j}A_r^* A_r u}
- D_t u\,\overline{D_t A_r^* A_r u}\right)\,J \,dw\,dt.
\end{equation}
But by definition, a solution of the wave equation $Pu=f$ satisfying the
Dirichlet boundary condition is $u\in H^1_{0,\loc}(X)$ with
\begin{equation*}
\int_X \left(\sum_{i,j}g^{ij}D_{w_i} u\,\overline{D_{w_j}v}
- D_t u\,\overline{D_t v}\right)\,J \,dw\,dt
=-\int_X f\,\overline{v}\,J \,dw\,dt
\end{equation*}
for every $v\in H^1_{0,c}(X)$. 
In particular, as $A_r^*A_r$ preserves $H^1_{0,\loc}(X)$,
this holds for $v=A_r^*A_r u$ when $A_r$ has a compactly supported Schwartz
kernel. If $f\in\dCinf(X)$, e.g.\ if $f=0$, the
right hand side now can also be estimated by the inductive hypothesis, showing
that $\|d_M A_r u\|^2-\|D_t A_r u\|^2$ is uniformly bounded as $r\to 0$.
The rest of the arguments presented above apply then, so we can conclude
that $q\nin\WFbz^{1,\infty}(u)$ as above.

This argument is immediately applicable for Neumann boundary conditions
as well. Thus, we still get \eqref{eq:DN-16} modulo terms that
can be estimated by the inductive hypothesis. 
Now, by definition, a solution of the wave equation $Pu=f$ satisfying the
Neumann boundary condition is $u\in H^1_{\loc}(X)$ with
\begin{equation}\label{eq:Neumann-solution-def}
\int_X \left(\sum_{i,j}g^{ij}D_{w_i} u\,\overline{D_{w_j}v}
- D_t u\,\overline{D_t v}\right)\,J \,dw\,dt
=-\int_X f\,\overline{v}\,J \,dw\,dt
\end{equation}
for every $v\in H^1_c(X)$. Here, for $f\in \dot H^{-1}_{\loc}(X)$, the
right hand side is the pairing of $\dot H^{-1}_{\loc}(X)$ with $H^1_c(X)$
via duality. In particular, as $A_r^*A_r$ preserves $H^1_{\loc}(X)$,
this holds for $v=A_r^*A_r u$, and the rest of the elliptic argument is
as for the Dirichlet boundary condition.

We use this opportunity to remark that our methods also immediately give
elliptic regularity for the Laplacian on $M$.

\begin{thm} (Microlocal elliptic regularity for $\Delta$.)
Suppose that $u\in H^1_{0,\loc}(M)$, and $\Delta u=f$, i.e.
\begin{equation*}
\langle d u,d v\rangle_M
=\langle f,v\rangle_M
\end{equation*}
for all $v\in H^1_{0,c}(M)$; here $\langle\cdot,\cdot\rangle_M$ is the
$L^2$ inner product on $M$. Then $\WFbz^{1,m}(u)\subset\WFbd^{-1,m}(f)$.

In particular, if $f\in H^{-1,m}_{b,\loc}(M)$ then $u\in H^{1,m}_{b,\loc}(M)$.

The same conclusions hold for Neumann boundary conditions,
i.e.\ with $H^1_0(M)$ replaced by $H^1(M)$.
\end{thm}

\begin{cor}
Suppose that $u\in H^1_{0,\loc}(M)$, and $(\Delta-\lambda) u=0$. Then
$u\in H^{1,\infty}_{\bl,\loc}(M)$. The conclusion also holds if $u$
satisfies Neumann boundary conditions.
\end{cor}

\begin{proof}
We have $\Delta u=f$ with $f=\lambda u\in H^1_{0,\loc}(M)\subset H^{-1,2}_{b,\loc}(M)$,
so $u\in H^{1,2}_{\bl,\loc}(M)$. Iterating this, using
$H^{1,m}_{\bl,\loc}\subset H^{-1,m+2}_{\bl,\loc}(M)$, completes the proof.
\end{proof}

\section{Bicharacteristics}\label{sec:bichar}
In this section we state the basic properties of generalized
broken bicharacteristics that are instrumental in proving the propagation
of singularities theorem in Section~\ref{thm:prop-sing}.
The philosophy originating from the work of Melrose and Sj\"ostrand
\cite{Melrose-Sjostrand:I, Melrose-Sjostrand:II} is that it is easier
to analyze the bicharacteristics (i.e.\ the `classical' system) precisely,
and prove only rough propagation estimates for the `quantum' system
(in this case the wave equation), essentially merely getting the
direction of the propagation correct, than to prove the precise propagation
statements directly, for many different aspects (not only the classical
geometry) interact in the latter setting. The precise propagation
statement is thus a combination of the rough propagation statements
with the detailed analysis of the bicharacteristics -- this is
the content of Section~\ref{sec:prop-sing} here.

Turning to the generalized broken bicharacteristics, these have been
described by Lebeau~\cite[Section~III]{Lebeau:Propagation} in his
setting, i.e.\ for domains $M$ in real analytic
manifolds $\tilde M$, equipped with a real analytic metric $g$,
with the boundary of $M$ admitting a stratification. However, analyticity
does {\em not} enter into the analysis of generalized broken bicharacteristic
(called `rayons' there), and manifolds with corners, by definition,
admit the desired stratification (stratified by the boundary faces),
in a $\Cinf$ sense. Thus, all of Lebeau's results on generalized broken
bicharacteristics apply in our setting, at least if one adopts his
definitions.

Our definition differs from that of Lebeau in two ways. First, at boundary
hypersurfaces (i.e.\ codimension $1$ faces),
Definition~\ref{Def:gen-br-bichar},
part (iii), demands more than Lebeau's definition
(from which (iii) is missing). Thus, our bicharacteristics are a subset
of those of Lebeau's. However, since the analysis of bicharacteristics
is local in $X$, the $\Cinf$ boundary analysis of Melrose and Sj\"ostrand
applies. As this only necessitates trivial changes, we point these
out below after the statement of the propositions of this section.

The other difference is that we defined the topology of
$\dot\Sigma$ as the subspace topology inherited from $\Tb^*X$, while
Lebeau defines it by requiring that $\hat\pi$ be continuous, so we need
to show that these are indeed the same, which we proceed to do now.

\begin{lemma}\label{lemma:topology-equiv}
Define the topology of $\dot\Sigma$ as the subspace topology of $\Tb^*X$.
Then $O\subset\dot\Sigma$ is open (resp.\ closed)
if and only if $\hat\pi^{-1}(O)$ is open (resp.\ closed).
\end{lemma}

Since the bundle inclusion map $\iota:T^*X\to\Tb^*X$ is $\Cinf$,
hence continuous,
$\hat\pi$ is automatically continuous, so it only remains to show
that if $\hat\pi^{-1}(O)$ is open, then $O$ is open, which we do below.

First, however we remark
that a basis of the subspace topology is given by
\begin{equation}\begin{split}\label{eq:top-basis-2}
B_\delta(q_0)=\{q\in\dot\Sigma:\ &|x(q)|<\delta,\ |y(q)-y_0(q)|<\delta,\ |t(q)-t(q_0)|<\delta,\\
&|\tau(q)-\tau(q_0)|<\delta,\ |\zeta(q)-\zeta(q_0)|<\delta\},
\end{split}\end{equation}
as $q_0$ and $\delta>0$ vary. Indeed, on $\dot\Sigma=\pi(\Char(P))$,
$|\sigma(q)|\leq C|x(q)|\,|\tau(q)|$ over compact subsets of $X$.
Assuming $\delta<1$, $\delta<|\tau(q_0)|/2$, as we may,
the above inequalities imply that $|\sigma(q)|<2C\delta|\tau(q_0)|$.
Given $\delta_0>0$, this set will thus be included in a $\delta_0$-ball
in $\Tb^*X$, centered at $q_0$,
provided we choose $\delta<\delta_0/2C|\tau(q_0)|$, so every neighborhood
of $q_0$ in $\dot\Sigma$ contains a set of the form \eqref{eq:top-basis-2}.

\begin{proof}[Proof of Lemma~\ref{lemma:topology-equiv}.]
We now show that if $\pih^{-1}(O)$ is open, then so is $O$.
That is, we need to show for any set $O$ with $\pih^{-1}(O)$ open, and for
any $q_0\in O\cap T^*\cF_{i,\reg}$,
there is a $\delta>0$ such that $B_\delta(q_0)\subset O$.
But $\pih^{-1}(\{q_0\})$
is the set of points $\tilde q_0=(x,y,t,\xi,\zeta,\tau)$ in $T^*X$
with $(x,y,t,\xi,\zeta,\tau)=(0,y(q_0),t(q_0),\xi,\zeta(q_0),\tau(q_0))$
and $\xi\cdot A(y(q_0))\xi=\tau(q_0)^2-|\zeta(q_0)|^2_{y(q_0)}$.
As $A$ is positive definite, the last equation implies that $\xi$ is
bounded on $\pih^{-1}(\{q_0\})$, and indeed $\pih^{-1}(\{q_0\})$ is
compact. So if $\pih^{-1}(O)$ open, then for some $\delta>0$
it contains the intersection of $\Char(P)$ with
the set
\begin{equation*}\begin{split}
\{\tilde q\in T^*X:\ &|x(\tilde q)|<\delta,\ |y(\tilde q)-y(q_0)|<\delta,
\ |t(\tilde q)-t(q_0)|<\delta,\\
&|\tau(\tilde q)-\tau(q_0)|<\delta,\ |\zeta(\tilde q)-\zeta(q_0)|<\delta,
|p(\tilde q)|<\delta\},
\end{split}\end{equation*}
i.e.\ it contains the set
\begin{equation*}\begin{split}
\tilde B_\delta(q_0)=
\{\tilde q\in \Char(P):\ &|x(\tilde q)|<\delta,\ |y(\tilde q)-y(q_0)|<\delta,
\ |t(\tilde q)-t(q_0)|<\delta,\\
&|\tau(\tilde q)-\tau(q_0)|<\delta,\ |\zeta(\tilde q)-\zeta(q_0)|<\delta\}.
\end{split}\end{equation*}
Now $\pih(\tilde B_\delta)=B_\delta(q_0)$, while $\pih(\pih^{-1}(O))=O$,
so we deduce that $B_\delta(q_0)\subset O$, and hence $O$ is open as
claimed.
\end{proof}

Being a subset of $\Tb^*X$, $\dot\Sigma$ is a separable,
locally compact metrizable space, although this follows also
directly using the topology induced by $\pih$ as in Lebeau's paper.

A stronger characterization of generalized broken bicharacteristics
at $\cH$ follows as in Lebeau's
paper. Notice that if $\gamma:I\to\dot\Sigma$ is continuous then the
conclusion of the following proposition certainly implies (i) and (ii)
((ii) follows as $x_j$ are $\pi$-invariant) of
Definition~\ref{Def:gen-br-bichar}, so the proposition
indeed provides an alternative to (i)-(ii) of
our definition. Note that (iii) is
not required for this proposition, and conversely, it does not imply (iii).
(We also remark paranthetically that there is yet another way of
phrasing (i) and (ii) in
the definition of generalized broken bicharacteristics, which is
important in $N$-body scattering in the presence of bound states, see
\cite[Definition~2.1]{Vasy:Bound-States}.)

\begin{prop}(Lebeau, \cite[Proposition~1]{Lebeau:Propagation})
\label{prop:Lebeau-2-sided}
If $\gamma$ is a generalized broken bicharacteristic, $t_0\in I$,
$q_0=\gamma(t_0)$, then there exist unique
$\tilde q_+, \tilde q_-\in\Char(P)$ satisfying
$\pi(\tilde q_\pm)=q_0$ and having
the property that if $f\in\Cinf(T^*X)$ is $\pi$-invariant
then $t\mapsto f_\pi(\gamma(t))$ is differentiable both from the left and
from the right at $t_0$ and
\begin{equation}\label{eq:Ham-23}
\left(\frac{d}{dt}\right)(f_\pi\circ\gamma)|_{t_0\pm}=H_p f(\tilde q_{\pm}).
\end{equation}
\end{prop}

\begin{cor}(Lebeau, \cite[Corollaire~2]{Lebeau:Propagation})
\label{cor:Lebeau-Lipschitz}
Suppose that $K$ is a compact subset of $\dot\Sigma$. Then there is
a constant $C>0$ such that for all generalized broken bicharacteristics
$\gamma:I\to K$, and for all $\pi$-invariant functions $f$
on a neighborhood of $\pi^{-1}(K)$ in $T^*X$, one has the uniform
Lipschitz estimate
\begin{equation*}
|f_\pi\circ\gamma(s_1)-f_\pi\circ\gamma(s_2)|\leq M\|f\|_{C^1}\,|s_1-s_2|,
\ s_1,s_2\in I.
\end{equation*}
In particular, (locally) the functions $x$, $\bar y$ and $\bar\zeta$
are Lipschitz on generalized broken bicharacteristics.
\end{cor}

We also need to analyze the uniform behavior of generalized broken
bicharacteristics. Here we quote Lebeau's results.

\begin{prop}(Lebeau, \cite[Proposition~5]{Lebeau:Propagation})
\label{prop:Lebeau-unif-limits}
Suppose that $K$ is a compact subset of $\dot\Sigma$, $\gamma_n:[a,b]\to K$ is
a sequence of generalized broken bicharacteristics which converge uniformly
to $\gamma$. Then $\gamma$ is a generalized broken bicharacteristic.
\end{prop}

\begin{proof}
By Lebeau's result, $\gamma$ is a `rayon', i.e.\ it satisfies (i)-(ii)
of Definition~\ref{Def:gen-br-bichar}. Thus, we only need
to show that it satisfies (iii) in order to prove that it is
a generalized broken bicharacteristic.
But if $\gamma(t_0)\in\cG\cap T^*\cF_{i,\reg}$,
$\cF_i$ a boundary hypersurface, then, using that the projection
of $\gamma$ to $X$ is Lipschitz by Corollary~\ref{cor:Lebeau-Lipschitz},
we see that for $\delta>0$ sufficiently
small, $\tilde\gamma_n=\gamma_n|_{[t_0-\delta,t_0+\delta]}$ lie
in $T^*X^\circ\cup T^*\cF_{i,\reg}$ for all $n$, as does
$\tilde\gamma=\gamma|_{[t_0-\delta,t_0+\delta]}$. Thus, $\tilde\gamma$
is a generalized broken bicharacteristic by the results
of \cite{Melrose-Sjostrand:II}, which implies that $\gamma$
satisfies (iii), finishing the proof.
\end{proof}

\begin{prop}(Lebeau, \cite[Proposition~6]{Lebeau:Propagation})
\label{prop:Lebeau-compactness}
Suppose that $K$ is a compact subset of $\dot\Sigma$, $[a,b]\subset\Real$
and
\begin{equation}
\calR=\{\text{generalized broken bicharacteristics}\ \gamma:[a,b]\to K\}.
\end{equation}
If $\calR$ is not empty then it is compact in the topology of uniform
convergence.
\end{prop}

\begin{proof}
$\calR$ is equicontinuous, as in Lebeau's proof (since every generalized
broken bicharacteristic is a rayon), so the proposition follows from
the theorem of Ascoli-Arzel\`a and Proposition~\ref{prop:Lebeau-unif-limits}.
\end{proof}

\begin{cor}(Lebeau, \cite[Corollaire~7]{Lebeau:Propagation})
\label{cor:Lebeau-bichar-ext}
If $\gamma:(a,b)\to\Real$ is a generalized broken bicharacteristic then
$\gamma$ extends to $[a,b]$.
\end{cor}

\section{The hyperbolic set}\label{sec:hyperbolic}
In $\cH\cup\cG$ the Dirichlet form is not positive, but
Lemma~\ref{lemma:Dirichlet-form} immediately gives the following
estimate, by simply rearranging its concluding estimate.
We do not need the sharp elliptic version,
as in Lemma~\ref{lemma:Dirichlet-form-2}, since
Lemma~\ref{lemma:Dirichlet-form} is only $1/2$ derivative weaker than
Lemma~\ref{lemma:Dirichlet-form-2}, and at $\cH\cup\cG$, $u$ loses
a whole derivative as compared to the elliptic estimates.

\begin{lemma}\label{lemma:Dt-dX}
Suppose that $K\subset\Sb^*X$ is compact, $U\subset\Sb^*X$ is open,
$K\subset U$.
Suppose that $\calA=\{A_r:\ r\in(0,1]\}$ be a bounded
family of ps.d.o's in $\Psibc^s(X)$ with $\WFb'(\calA)\subset K$, and
with $A_r\in\Psib^{s-1}(X)$ for $r\in (0,1]$.
Then there exist $B\in\Psib^{s-1/2}(X)$, $\tilde B\in\Psib^{s+1/2}(X)$
with $\WFb'(B),\WFb'(\tilde B)\subset U$
and $C_0>0$ such that for $r\in(0,1]$, $u\in H^1_{0,\loc}(X)$
with $\WFbz^{1,s-1/2}(u)\cap U=\emptyset$,
$\WFbd^{-1,s+1/2}(Pu)\cap U=\emptyset$, the following estimate holds:
\begin{equation*}\begin{split}
\|d_M A_r u\|^2\leq&\|D_t A_r u\|^2\\
&\quad+C_0(\|u\|^2_{H^1_{\loc}(X)}+\|Bu\|^2_{H^1(X)}
+\|Pu\|^2_{H^{-1}_{\loc}(X)}
+\|\tilde B Pu\|^2_{H^{-1}(X)})).
\end{split}\end{equation*}
In particular, if the assumption on $Pu$ is strengthened to $Pu=0$, we have
\begin{equation*}
\|d_M A_r u\|^2\leq\|D_t A_r u\|^2+C_0(\|u\|^2_{H^1_{\loc}(X)}+\|Bu\|^2_{H^1(X)}).
\end{equation*}
The meaning of $\|u\|^2_{H^1_{\loc}(X)}$ and
$\|Pu\|^2_{H^{-1}_{\loc}(X)}$ is stated in
Remark~\ref{rem:localize}.
\end{lemma}

This lemma roughly says that $D_{x_i}A_r u$ (and also $D_{y_i}A_r u$, but the
latter follows more directly from general
properties of the b-ps.d.o's near $\cH\cup\cG$) is bounded by
$D_t A_r u$, modulo lower order error terms. This
allows us to estimate various error terms in the positive
commutator argument below, and it shows that we only need to find
a uniform
bound on $\|D_t A_r u\|^2$ in terms of other terms on the right hand side
in order to get a bound on $\|d_M A_r u\|^2$, hence conclude that
points at which $\sigma_{b,s}(A)\neq 0$ do not lie in $\WFbz^{1,s}(u)$.
(Here $A_r\to A$ in a suitable sense.)

A related
consequence of this lemma is that for microlocal solutions of $Pu=0$,
$u\in H^1_0(X)$, $\WFbz^{1,m}(u)$ agrees
with the b-wave front set of $u$ defined with respect to the more traditional
$L^2$ space.

\begin{lemma}\label{lemma:H10-L2}
Suppose $u\in H^1_{0,\loc}(X)$, $\WFbd^{-1,\infty}(Pu)=\emptyset$. Then
\begin{equation*}
\WFbz^{1,m}(u)^c=\{q\in\Tb^*X\setminus o:
\ \exists A\in\Psib^{m+1}(X),\ \sigma_{b,m+1}(A)(q)
\neq 0,\ Au\in L^2(X)\}.
\end{equation*}

More generally, for $u\in H^1_{0,\loc}(X)$,
\begin{equation*}\begin{split}
&\WFbz^{1,m}(u)^c\cap\WFbd^{-1,\infty}(Pu)^c\\
&\ =\{q\in\WFbd^{-1,\infty}(Pu)^c:
\ \exists A\in\Psib^{m+1}(X),\ \sigma_{b,m+1}(A)(q)
\neq 0,\ Au\in L^2(X)\}.
\end{split}\end{equation*}
\end{lemma}

\begin{proof}
In $T^*X^\circ$, both sides are the standard wave front set, $\WF^{m+1}(u)$,
so it suffices to consider the case when $q$ lies over $\pa X$.

First we show that the left hand side is a subset of the right hand side,
which is the `easy direction', and does not use any condition on $Pu$.
Now, if $q\in\WFbz^{1,m}(u)^c$, then there is some $B\in\Psib^m(X)$
with $\sigma_{b,m}(B)(q)\neq 0$ and $Bu\in H^1_0(X)$. We may assume that
$B$ is supported near the projection of $q$ to $X$, so in particular
we can use local coordinates in the rest of the argument.
If $\zeta_j(q)\neq 0$, then $A=D_{y_j}B\in\Psib^{m+1}(X)$
with non-vanishing principal symbol at $q$ and $D_{y_j}Bu\in L^2(X)$ since
$Bu\in H^1_0(X)$, so $q$ indeed lies in the right hand side.
A similar argument works of $\tau(q)\neq 0$. If $\sigma_j(q)\neq 0$,
then $A=x_j D_{x_j}B\in\Psib^{m+1}(X)$
with non-vanishing principal symbol at $q$ and $D_{x_j}Bu\in L^2(X)$
since $Bu\in H^1_0(X)$, so $x_jD_{x_j}Bu\in L^2(X)$ as well -- thus, again,
$q$ lies in the right hand side. Therefore the left hand side is indeed
a subset of the right hand side.

To see the converse direction, i.e.\ that the right hand side is a subset
of the left hand side, we note that as $u\in H^1_{0,\loc}(X)$,
$\WFbz^{1,m}(u)^c
\supset((\dot\Tb^*X)^c\cup\cE)\setminus \WFbd^{-1,\infty}(Pu)$
by Proposition~\ref{prop:elliptic},
so it suffices to consider $q\in\cG\cup\cH$.
We use induction on $m$ to prove that if $q$ is in the right hand side 
then it is also in the left hand side -- with the case $m=0$ being trivial
as we are assuming $u\in H^1_{0,\loc}(X)$. In general, suppose that the
inclusion has been proved for $m$ replaced by $m-1/2$.
Suppose that $q\in\cG\cup\cH$ is in the
right hand side, so there is $A\in\Psib^{m+1}(X)$, $A$ elliptic at $q$,
$Au\in L^2(X)$, and $q\nin\WFbz^{1,m-1/2}(u)$ by the inductive hypothesis.
Note that $\tau(q)\neq 0$, i.e.\ $D_t$ is elliptic
at $q$. We may assume that $\WFb'(A)$ lies close to $q$, hence that
$\tau$ is elliptic on $\WFb'(A)$, and in addition $\WFbz^{1,m-1/2}(u)\cap
\WFb'(A)=\emptyset$. Then we can write $A=D_t B+R$,
$B\in\Psib^m(X)$ elliptic at $q$ and $R\in\Psib^{-\infty}(X)$.
Thus, (as $u\in L^2(X)$) $Ru\in L^2(X)$, so $D_t Bu\in L^2(X)$.
Taking $B_r\in\Psib^{m-1}(X)$ uniformly bounded with $B_r\to B$ in $\Psibc^{m+\ep}(X)$ ($\ep>0$),
Lemma~\ref{lemma:Dt-dX} gives that $d_M B_r u$ is uniformly bounded in $L^2$.
Since it converges to $d_M Bu$ in $\dist(X)$ on the one hand, and there must
be a weakly convergent sequence $d_M B_{r_k} u$ in $L^2(X)$, $r_k\to 0$ as
$k\to\infty$, by the uniform bound, we deduce that
$d_M Bu\in L^2(X)$
as well, so $d_X Bu\in L^2(X)$, hence $Bu\in H^1_0(X)$.
\end{proof}

After these preliminary discussions, we turn to the propagation estimate
at $q\in\cH$. As usual, the key ingredient is to find a $\Cinf$ function $f$
on $\Tb^* X$ such that, at least near $q$, $H_p \iota^* f$ has a fixed
sign. We usually drop the pull-back $\iota^*$ below; recall that
$\iota:T^*X\to\Tb^*X$ is the `inclusion', and $\pi$ is $\iota$, considered
as a map onto $\dot\Tb^*X$. In our setting, we can take $f=\eta$
where $\eta=-\frac{x\cdot\xi}{|\tau|}=-\frac{\sum\sigma_j}{|\tau|}$.
Indeed, the Hamilton vector field $H_p$ of $p$ is given by
\begin{equation}\begin{split}\label{eq:H_p-form}
H_p=2\tau\pa_t-H_g&=2\tau\pa_t-2 A\xi\cdot \pa_x-2B\zeta\cdot \pa_y
-2\sum C_{ij}\zeta_j\pa_{x_i}-2\sum C_{ij}\xi_i\pa_{y_j}\\
&\qquad+2\sum (\pa_{x_k} A_{ij})\xi_i\xi_j\pa_{\xi_k}
+2\sum (\pa_{x_k} C_{ij})\xi_i\zeta_j\pa_{\xi_k}\\
&\qquad\qquad\qquad
+2\sum (\pa_{x_k} B_{ij})\zeta_i\zeta_j\pa_{\xi_k}\\
&\qquad+2\sum (\pa_{y_k} A_{ij})\xi_i\xi_j\pa_{\zeta_k}
+2\sum (\pa_{y_k} C_{ij})\xi_i\zeta_j\pa_{\zeta_k}\\
&\qquad\qquad\qquad
+2\sum (\pa_{y_k} B_{ij})\zeta_i\zeta_j\pa_{\zeta_k}.
\end{split}\end{equation}
Thus,
\begin{equation*}\begin{split}
|\tau| H_p\eta=2\xi\cdot A\xi+2\sum C_{ij}\xi_i\zeta_j
&-2\sum (\pa_{x_k} A_{ij})\xi_i\xi_j x_k\\
&-2\sum (\pa_{x_k} C_{ij})\xi_i\zeta_j x_k
-2\sum (\pa_{x_k} B_{ij})\zeta_i\zeta_j x_k,
\end{split}\end{equation*}
so at $x=0$, where $C$ vanishes,
\begin{equation}\label{eq:H_p-eta-at-0}
|\tau| H_p\eta=2\xi\cdot A\xi=2\tau^2-2\zeta\cdot B\zeta-2p=2\tau^2
-2|\zeta|_y^2-2p.
\end{equation}
Thus, $H_p\eta>0$ at $\pi^{-1}(\cH)\cap\Char(P)=\pih^{-1}(H)$.

We only state the following propagation result for propagation in the
forward direction along the generalized broken bicharacteristics.
A similar result holds in the backward direction, i.e.\ if we replace
$\eta(\xi)<0$ by $\eta(\xi)>0$
in \eqref{eq:prop-9a};
the proof in this case only requires changes in some signs
in the argument given below. The construction of a positive commutator below
closely mirrors that of \cite{Vasy:Propagation-Many} in the $N$-body setting.

\begin{prop}\label{prop:normal-prop}
Let $q_0=(y_0,t_0,\zeta_0,\tau_0)\in\cH\cap T^* \cF_{\reg}$
and let $\eta=-\frac{x\cdot\xi}{|\tau|}$ be the $\pi$-invariant function defined
in the local coordinates discussed above, and suppose that $u\in
H^1_{0,\loc}(X)$, $q_0\nin\WFbd^{-1,\infty}(Pu)$.
If there exists a conic neighborhood $U$ of $q_0$ in $\dot\Tb^*X$ such that
\begin{equation}\begin{split}\label{eq:prop-9a}
q\in U\Mand \eta(q)<0\Rightarrow q\nin\WFbz^{1,\infty}(u)
\end{split}\end{equation}
then $q_0\nin\WFbz^{1,\infty}(u)$.

In fact, if the wave front set assumptions are relaxed to
$q_0\nin\WFbd^{-1,s+1}(Pu)$ and the existence of a conic neighborhood $U$ of
$q_0$ in $\dot\Tb^*X$ such that
\begin{equation}\begin{split}\label{eq:prop-9a-s}
q\in U\Mand \eta(q)<0\Rightarrow q\nin\WFbz^{1,s}(u),
\end{split}\end{equation}
then we can still conclude that $q_0\nin\WFbz^{1,s}(u)$.
\end{prop}

\begin{rem}
Note that $\eta(q)<0$ implies $x\neq 0$, so $q\nin T^*\cF$.
\end{rem}

\begin{rem}
We recall that every conic neighborhood $U$ of
$q_0=(y_0,t_0,\zeta_0,\tau_0)\in\cH\cap T^* \cF_{\reg}$ in $\dot\Sigma$
contains an open set of the form
\begin{equation}\label{eq:prop-rem-9b}
\{q:\ |x(q)|^2+|y(q)-y_0|^2+|t(q)-t_0|^2+|\hat\zeta(q)
-\hat\zeta_0|^2<\delta\},
\end{equation}
$\hat\zeta=\frac{\zeta}{\tau}$.
Note also that
\eqref{eq:prop-9a} implies the same statement with $U$ replaced by
any smaller neighborhood of $q_0$; in particular, for
the set \eqref{eq:prop-rem-9b}, provided that $\delta$ is sufficiently small.
We can also assume that $\WFbd^{-1,\infty}(Pu)\cap U
=\emptyset$.
\end{rem}

\begin{proof}
As in Proposition~\ref{prop:elliptic} we use an inductive argument to
show that $q_0\nin\WFbz^{1,s}(u)$, provided that $q_0\nin\WFbz^{1,s-1/2}(u)$;
again the inductive hypothesis holds for $s=1/2$ since $u\in H^1_{\loc}(X)$.
Because of Lemma~\ref{lemma:Dt-dX}, we only need to show that for some
$B\in\Psib^{s+1}(X)$ with $\sigma_{b,s+1}(B)(q_0)\neq 0$, $Bu\in L^2(X)$.

Below we fix a small neighborhood $U_0$ of $q_0$ such that
$U_0$ is inside a coordinate neighborhood of $q_0$ and $\WFbd^{-1,\infty}(Pu)\cap
U_0=\emptyset$.

The key is to construct an operator $A$ with $\WFb'(u)\cap U$ and
$i[A^*A,P]$ positive, modulo
terms that we can estimate either by the a priori assumptions, namely those
on $Pu$ and those on $\WFb(u)$, summarized in
\eqref{eq:prop-9a} above. Thus, we do not need to make the commutator
positive in $\eta<0$, and also `away from $\Char(P)$', although the
latter is a moral statement as the locus of the microlocalization is
$\Tb^*X\setminus o$, not $T^*X\setminus o$. Our $A$ will in fact be formally
self-adjoint modulo lower order operators, and we only take $A^*A$ to avoid
having to comment on the subprincipal terms.

The main technical problem below is that $P$ does not lie in $\Psib(X)$,
so we cannot simply use the symbol calculus on $\Psib(X)$ -- we need
to write out various expressions semi-explicitly as elements
of $\Diff\Psib(X)$. On the other hand, while $\Psib(X)$ is the
locus of the microlocalization, at the level of the symbol calculus
one can rely on standard ps.d.o's on an extension $\tilde X$ of $X$,
i.e.\ work with symbols on $T^*X$. This has the advantage that $p$ is
a symbol on $T^*X$, as is the pull-back of symbols on $\Tb^*X$ via
$\pi$, so one can calculate their Poisson bracket, etc. However, it
is not trivial to make this into a technically useful computation,
since we need to control various expression in $\Diff\Psib(X)$.
In order to make the argument more digestable, we start with a symbol
construction, and do a formal commutator computation in $\Psi(\tilde X)$
(in fact, we will ignore that we need an extension $\tilde X$ here
and write `$\Psi(X)$' at times)
to show why the constructed symbol {\em should} be useful, and then
give the actual proof.

We construct the symbol of $A$ in a few steps. The two main ingredients are
a homogeneous degree zero
function that is increasing along the Hamilton flow, which will be $\eta$,
and a homogeneous degree zero
function $\omega$ on a conic neighborhood of $q_0$ in
$\Tb^*X\setminus o$ that roughly measures the square of the distance
from $q_0$ in $\dot\Tb^*X$. Note that $\omega$ can also be regarded
as a function on a subset of $\Sb^*X$, if desired.
Thus, we let
\begin{equation}\label{eq:prop-omega-def}
\omega(q)=|x(q)|^2+|y(q)-y_0|^2+|t(q)-t_0|^2+|\hat\zeta(q)
-\hat\zeta_0|^2,
\end{equation}
$|.|$ denoting the Euclidean norm, and $\hat\zeta=\frac{\zeta}{\tau}$ as
above.
Then $\omega$ vanishes quadratically at $q_0$, in fact is a sum of
squares, so $|d\omega|\leq C'_1\omega^{1/2}$, and in particular
\begin{equation}\label{eq:omega-est-a0}
|\tau^{-1}H_p\omega|\leq C_1''\omega^{1/2}.
\end{equation}
Were we merely using the symbol calculus for $\Psib(X)$ or `$\Psi(X)$',
this is all that would matter. Since this is not the case,
we need that more explicitly,
\begin{equation}\begin{split}\label{eq:omega-est-a}
&\tau^{-1}H_p\omega=f_0+\sum_i f_i\tau^{-1}\xi_i
+\sum_{i,j}f_{ij}\tau^{-2}\xi_i\xi_j,\\
&\qquad\qquad f_i,f_{ij}\in\Cinf(\Tb^*X),\ |f_i|,|f_{ij}|\leq C_1\omega^{1/2},
\end{split}\end{equation}
$f_i$, $f_{ij}$ homogeneous of degree $0$, which follows from
\eqref{eq:H_p-form}.

Next, we use the variable $\eta=-\frac{x\cdot\xi}{|\tau|}$ to measure
propagation. Since
\begin{equation*}
\eta=-\frac{x\cdot\xi}{|\tau|}=-\sum_j\sigma_j|\tau|^{-1},
\end{equation*}
$\eta$ is a homogeneous degree zero
$\Cinf$ function on a conic neighborhood of $q_0$ in $\Tb^*X\setminus o$,
hence it (or more precisely its pullback by $\pi$)
is a $\Cinf$, $\pi$-invariant function on $T^*X$.
This function indeed measures the flow along bicharacteristics near $q_0$
since at points $\tilde q_0$ in $\pih^{-1}(\{q_0\})$, where thus $p=0$,
\begin{equation}
|\tau| H_p\eta(\tilde q_0)=\tau_0^2-|\zeta_0|^2_{y_0}=c_0\tau_0^2>0,
\end{equation}
due to \eqref{eq:H_p-eta-at-0}, where we used that $q_0\in\cH$. Again,
if we could use `$\Psi(X)$', all we would need is that
$|\tau| H_p\eta>c_0\tau^2/2>0$
on $U_0$, which is automatic if the neighborhood $U_0$
is small enough.
Now, however, we need the more explicit expression
\begin{equation*}\begin{split}
|\tau|^{-1}H_p\eta=&\tau^{-2}(2\tau^2-2|\zeta|^2-2p)+g_0+\sum_i \xi_i \tau^{-1}
g_i+\sum_{i,j}g_{ij}\tau^{-2}\xi_i\xi_j,\\
&\ g_i,g_{ij}\in\Cinf(\Tb^*X),
\ |g_i|, |g_{ij}|\leq C_1\omega^{1/2},
\end{split}\end{equation*}
$g_i$, $g_{ij}$ homogeneous of degree $0$, which again follows from
\eqref{eq:H_p-form}.

We are now ready to define the symbol $a$ of $A$.
For $\ep>0$, $\delta>0$, with other restrictions to be imposed later on,
let
\begin{equation}
\phi=\eta+\frac{1}{\ep^2\delta}\omega,
\end{equation}
so $\phi$ is a homogeneous degree zero
$\Cinf$ function on a conic neighborhood of $q_0$ in $\Tb^*X\setminus o$
-- we can again regard it as a $\pi$-invariant function on $T^*X\setminus o$.
(Here $\ep^{-2}$ plays the role of $\beta$ in the analogous -- normal --
propagation estimate of \cite{Vasy:Propagation-Many}.)

Let $\chi_0\in\Cinf(\Real)$ be equal to $0$ on $(-\infty,0]$ and
$\chi_0(t)=\exp(-1/t)$ for $t>0$. Thus, $\chi_0'(t)=t^{-2}\chi_0(t)$.
Let $\chi_1\in\Cinf(\Real)$ be $0$
on $(-\infty,0]$, $1$ on $[1,\infty)$, with $\chi_1'\geq 0$ satisfying
$\chi_1'\in\Cinf_c((0,1))$. Finally, let $\chi_2\in\Cinf_c(\Real)$
be supported in $[-2c_1,2c_1]$, identically $1$ on $[-c_1,c_1]$,
where $c_1$ is such that if
$|\sigma|^2/\tau^2<c_1/2$ in $\dot\Sigma\cap U_0$. Thus,
$\chi_2(|\sigma|^2/\tau^2)$ is a cutoff in $|\sigma|/|\tau|$, with its
support properties ensuring that $d\chi_2(|\sigma|^2/\tau^2)$ is supported in
$|\sigma|^2/\tau^2\in [c_1,2c_1]$ hence outside $\dot\Sigma$ -- it should
be thought of as a factor that microlocalizes near the characteristic
set but effectively commutes with $P$.
Then, for $A_0>0$ large, to be determined, let
\begin{equation}\label{eq:prop-22}
a=\chi_0(A_0^{-1}(2-\phi/\delta))\chi_1(\eta/
\delta+2)\chi_2(|\sigma|^2/\tau^2);
\end{equation}
so $a$ is a homogeneous degree zero $\Cinf$ function on a conic neighborhood
of $q_0$ in $\Tb^*X$.
Indeed, as we see momentarily, for any $\ep>0$,
$a$ has compact
support inside this neighborhood (regarded as a subset of $\Sb^*X$, i.e.
quotienting out by the $\Real^+$-action) for $\delta$ sufficiently small,
so in fact it is globally well-defined.
In fact, on $\supp a$ we have $\phi\leq 2\delta$ and $\eta\geq-2\delta$.
Since $\omega\geq 0$, the first of these inequalities implies that
$\eta\leq 2\delta$, so on $\supp a$
\begin{equation}
|\eta|\leq 2\delta.
\end{equation}
Hence,
\begin{equation}\label{eq:omega-delta-est}
\omega\leq \ep^2\delta(2\delta-\eta)\leq4\delta^2\ep^2.
\end{equation}
In view of \eqref{eq:prop-omega-def} and \eqref{eq:prop-rem-9b},
this shows that for any $\ep>0$, $a$ is supported in $U$, provided
$\delta>0$ is sufficiently small.
The role that $A_0$ large
plays is that it increases the size of the first derivatives of $a$ relative
to the size of $a$, hence it allows us to give a bound for
$a$ in terms of a small multiple of
its derivative along the Hamilton vector field. This is crucial as we
need to deal with weight factors, such as $|\tau|^{s+1/2}$ in the
next paragraph, if the weight factors do not commute with $P$.
In this case, they can be arranged to commute (at least microlocally, which
suffices), so we could eliminate $A_0$, but its presence is helpful if
one is to weaken the assumptions on the structure of $P$.

This is the point where the technical argument needs significantly more
details than the motivational one.
So we start with the motivation.
Thus, using \eqref{eq:omega-est-a0}, \eqref{eq:omega-delta-est},
\begin{equation*}
|\tau|^{-1} H_p\phi=H_p\eta+\frac{1}{\ep^2\delta}H_p\omega\geq c_0/2
-\frac{1}{\ep^2\delta}C_1''\omega^{1/2}\geq c_0/2-2C_1''\ep^{-1}\geq
c_0/4>0
\end{equation*}
provided that $\ep>\frac{8C_1''}{c_0}$, i.e.\ that
$\ep$ is not too small. We fix some such $\ep$ for the rest of
the arguments below, and then we will take $\delta>0$ sufficiently
small. With this,
\begin{equation*}
H_p a^2=-b^2+e,\ b=|\tau|^{1/2} (2|\tau|^{-1} H_p\phi)^{1/2}
(A_0\delta)^{-1/2}(\chi_0\chi_0')^{1/2}\chi_1\chi_2,
\end{equation*}
with $e$ arising from the derivative of $\chi_1\chi_2$. Here
$\chi_0$ stands for $\chi_0(A_0^{-1}(2-\frac{\phi}{\delta}))$, etc.
Since $\eta<0$
on $d\chi_1$ while $d\chi_2$ is disjoint from the characteristic set,
both being regions disjoint from $\WFb(u)$,
$i[A^*A,P]$ is positive modulo terms that we can a priori control,
so the standard positive commutator argument gives an estimate
for $Bu$, where $B$ has symbol $b$. Replacing $a$ by $a|\tau|^{s+1/2}$,
we still have a positive commutator (in this case $\tau$, or rather $D_t$,
actually commutes
with $P$, but in any case we could use $A_0$ to bound the additional
commutator term), which now gives (with the new $B$) that $Bu\in L^2(X)$,
which means in particular that $q_0\nin\WFbz^{1,s}(u)$.

This argument is of course {\em very} imprecise.
The technically correct version is the following. First,
\begin{equation}\begin{split}\label{eq:H_p-phi}
|\tau|^{-1} H_p\phi&=|\tau|^{-1}H_p\eta+\frac{1}{\ep^2\delta}|\tau|^{-1}H_p\omega\\
&=-2p\tau^{-2}
+\tau^{-2}(2\tau^2-2|\zeta|^2_y)
+g_0+\sum_i \tau^{-1}\xi_i g_i+\sum_{ij}\tau^{-2}\xi_i\xi_j g_{ij}\\
&\qquad+\frac{1}{\ep^2\delta}
(f_0+\sum \xi_i \tau^{-1}f_i+\sum\tau^{-2}\xi_i\xi_j g_{ij})
\end{split}\end{equation}
Let $\tilde B\in\Psib^{1/2}(X)$ with
\begin{equation}\label{eq:tilde-b}
\tilde b=\sigma_{b,1/2}(\tilde B)=|\tau|^{1/2}
(A_0\delta)^{-1/2}(\chi_0\chi_0')^{1/2}\chi_1\chi_2\in\Cinf(\Tb^*X\setminus o),
\end{equation}
and let $A\in\Psib^0(X)$ with $\sigma_{b,0}(A)=a$. Again,
$\chi_0$ stands for $\chi_0(A_0^{-1}(2-\frac{\phi}{\delta}))$, etc.
Also, let $C\in\Psib^0(X)$
have symbol $\sigma_{b,0}(C)=|\tau|^{-1}(2\tau^2-2|\zeta|^2_y)^{1/2}\psi$ where
$\psi\in S^0(\Tb^*X)$ is identically $1$ on $U$ considered as a subset
of $\Tb^*X$.
Then an explicit calculation using Lemma~\ref{lemma:comm-symbol}
and $P=D_t^2-\Delta$,
\begin{equation*}
\Delta=\sum_{i,j}A_{ij}(x,y)D_{x_i}D_{x_j}
+\sum_{i,j}2C_{ij}(x,y)D_{x_i}D_{y_j}+\sum_{i,j}B_{ij}(x,y)D_{y_i}D_{y_j}
+P_1,
\end{equation*}
$P_1\in\Diff^1(X)$, gives,
in accordance with \eqref{eq:H_p-phi},
\begin{equation}\begin{split}\label{eq:P-comm}
&i[A^*A,P]\\
&\quad=R' P+\tilde B^*(C^*C+R_0+\sum_i D_{x_i} R_i
+\sum_{ij} D_{x_i}R_{ij}D_{x_j})\tilde B+R''+E+E'
\end{split}\end{equation}
with
\begin{equation*}\begin{split}
&R_0\in\Psib^0(X),\ R_i\in\Psib^{-1}(X),\ R_{ij}\in\Psib^{-2}(X),\\
&R'\in\Psib^{-1}(X),\ R''\in\Diff^2\Psib^{-2}(X),
\ E,E'\in\Diff^2\Psib^{-1}(X),
\end{split}\end{equation*}
with $\WFb'(E)\subset\eta^{-1}((-\infty,-\delta])
\cap U$, $\WFb'(E')\cap\dot\Sigma=\emptyset$
($E$ arises from the commutator of $P$ with an operator with symbol
$\chi_1(\eta/\delta+2)$, while $E'$ from the commutator of $P$ with an operator
with symbol $\chi_2(|\sigma|^2/\tau^2)$)
and with
$r_0=\sigma_{b,0}(R_0)$, $r_i=\sigma_{b,-1}(R_i)$,
$r_{ij}\in\sigma_{b,-2}(R_{ij})$,
\begin{equation*}
|r_0|\leq C_2(1+\frac{1}{\ep^2\delta})\omega^{1/2},
\ |\tau r_i|\leq C_2(1+\frac{1}{\ep^2\delta})\omega^{1/2},
\ |\tau^2 r_{ij}|\leq C_2(1+\frac{1}{\ep^2\delta})\omega^{1/2},
\end{equation*}
and $\supp r_j$ lying in $\omega\leq 9\delta^2\ep^2$. Thus,
\begin{equation*}
|r_0|\leq 3C_2(\delta\ep+\ep^{-1}),
\ |\tau r_i|\leq 3C_2(\delta\ep+\ep^{-1}),
\ |\tau^2 r_{ij}|\leq 3C_2(\delta\ep+\ep^{-1}).
\end{equation*}

Having calculated the commutator, we proceed to estimate the `error terms'
$R_0$, $R_i$, $R_{ij}$ as operators. We start with $R_0$.
As follows from the standard square root construction to prove the
boundedness of ps.d.o's on $L^2$, there exists $R_0'\in\Psib^{-1}(X)$ such that
\begin{equation*}
\|R_0v\|\leq2\sup|r_0|\,\|v\|+\|R_0'v\|
\end{equation*}
for all $v\in L^2(X)$. Here $\|\cdot\|$ is the $L^2(X)$-norm, as usual.
Thus, we can estimate, for any $\gamma>0$,
\begin{equation*}\begin{split}
|\langle R_0 v,v\rangle|&\leq \|R_0 v\|\,\|v\|
\leq 2\sup |r_0|\,\|v\|^2+\|R_0' v\|\,\|v\|\\
&\leq 6C_2(\delta\ep+\ep^{-1})\|v\|^2
+\gamma^{-1}\|R_0' v\|^2+\gamma \|v\|^2.
\end{split}\end{equation*}

Now we turn to $R_i$.
Let $T\in\Psib^{-1}(X)$ be elliptic (which we use to keep track
of the orders of ps.d.o's), $T^-\in\Psib^1(X)$ a parametrix, so
$T^-T=\Id+F$, $F\in\Psib^{-\infty}(X)$.
Then there exist
$R'_i\in\Psib^{-1}(X)$ such that
\begin{equation*}\begin{split}
\|R_i w\|=\|R_i (T^- T -F)w\|&\leq\|(R_i T^-)(Tw)\|+\|R_iFw\|\\
&\leq 6C_2(\delta\ep+\ep^{-1})\|Tw\|+
\|R_i' Tw\|+\|R_iFw\|
\end{split}\end{equation*}
for all $w$ with $Tw\in L^2(X)$.
Similarly, there exist $R'_{ij}\in\Psib^{-1}(X)$ such that
\begin{equation*}
\|(T^-)^*R_{ij} w\|\leq 6C_2(\delta\ep+\ep^{-1})\|Tw\|+
\|R_{ij}' Tw\|+\|(T^-)^*R_{ij}Fw\|
\end{equation*}
for all $w$ with $Tw\in L^2(X)$.
Thus,
\begin{equation*}\begin{split}
|\langle R_i D_{x_i} v,v\rangle|\leq
&6C_2(\delta\ep+\ep^{-1})\|TD_{x_i} v\|\,\|v\|\\
&\qquad+2\gamma\|v\|^2+\gamma^{-1}\|R'_i TD_{x_i} v\|^2+\gamma^{-1}\|F_i D_{x_i}v\|^2,
\end{split}\end{equation*}
and, writing $D_{x_j}v=T^-Tv-Fv$ in the right factor, and taking the
adjoint of $T^-$,
\begin{equation*}\begin{split}
|\langle R_{ij} D_{x_i} v,D_{x_j}v\rangle|\leq
&6C_2(\delta\ep+\ep^{-1})\|TD_{x_i}v\|\,\|TD_{x_j}v\|\\
&\qquad
+2\gamma\|TD_{x_j}v\|^2+\gamma^{-1}\|R'_{ij}T D_{x_i}v\|^2
+\gamma^{-1}\|F_{ij} D_{x_i}v\|^2\\
&\qquad+\|R_{ij}D_{x_i}v\|\,\|FD_{x_j}v\|,
\end{split}\end{equation*}
with $F_i,F_{ij}\in\Psib^{-\infty}(X)$.

Let $\Lambda_r$ have symbol
\begin{equation}
|\tau|^{s+1/2}(1+r|\tau|^2)^{-s},\quad r\in[0,1),
\end{equation}
so $A_r=A\Lambda_r\in\Psib^{0}(X)$ for $r>0$ and it is
uniformly bounded in $\Psibc^{s+1/2}(X)$.
In similar constructions in general, the commutator $[P,\Lambda_r]$
can be controlled by the
other terms using $A_0$, for $A_0$ large -- in the present setting
$[P,\Lambda_r]=0$.

Now, by \eqref{eq:P-comm},
\begin{equation}\begin{split}\label{eq:pos-comm}
\langle i[A_r^*A_r,P]u,u\rangle
&=\|C\tilde B\Lambda_r u\|^2
+\langle R'P\Lambda_r u,\Lambda_r u\rangle
+\langle R_0 \tilde B\Lambda_r u,\tilde B\Lambda_r u\rangle\\
&\qquad+\sum \langle R_i D_{x_i} \tilde B \Lambda_r u,\tilde B\Lambda_r u\rangle
+\sum \langle R_{ij} D_{x_i} \tilde B\Lambda_r u,
D_{x_j}\tilde B\Lambda_r u\rangle\\
&\qquad+\langle R''\Lambda_r u,\Lambda_r u\rangle
+\langle (E+E')\Lambda_r u,\Lambda_r u\rangle
\end{split}\end{equation}
On the other hand, as $A_r\in\Psib^0(X)$ for $r>0$ and $u\in H^1_0(X)$, so
$A_r^*A_ru\in H^1_0(X)$,
\begin{equation}\begin{split}\label{eq:comm-expansion}
\langle [A_r^*A_r,P]u,u\rangle
&=\langle A_r^*A_r Pu,u\rangle-\langle PA_r^*A_ru,u\rangle\\
&=\langle A_r Pu,A_ru\rangle-\langle A_ru,A_rPu\rangle
=2i\im\langle A_rPu,A_ru\rangle;
\end{split}\end{equation}
the pairing makes sense for $r>0$ since $A_r\in\Psib^{0}(X)$ then.

Assume for the moment that $\WFbd^{-1,s+3/2}(Pu)\cap U=\emptyset$ -- this
is certainly the case in our setup if $q_0\nin\WFbd^{-1,\infty}(Pu)$, but this
assumption is a little stronger that $q_0\nin\WFbd^{-1,s+1}(Pu)$,
which is what we need to assume for the second paragraph in the
statement of the proposition. We deal with the weakened
hypothesis $q_0\nin\WFbd^{-1,s+1}(Pu)$ at the end of the proof.
Returning to \eqref{eq:comm-expansion}, the utility of the commutator
calculation is that we have good information about $Pu$ (this is
where we use that we have a microlocal solution of the PDE!). Namely,
we estimate the right hand side as
\begin{equation}\begin{split}\label{eq:Pu-s+32}
|\langle A_r Pu,A_r u\rangle|&\leq|\langle (T^-)^*A_rPu,TA_ru\rangle|
+|\langle A_rPu,FA_r u\rangle|\\
&\leq \|(T^-)^*A_r Pu\|_{H^{-1}(X)}
\|TA_r u\|_{H^1(X)}\\
&\qquad+\|A_r Pu\|_{H^{-1}(X)}\|FA_r u\|_{H^1(X)}.
\end{split}\end{equation}
Since $(T^-)^*A_r$ is uniformly bounded in $\Psibc^{s+3/2}(X)$,
$TA_r$ is uniformly bounded in $\Psibc^{s-1/2}(X)$, both with $\WFb'$
in $U$, with $\WFbd^{-1,s+3/2}(Pu)$, resp.\ $\WFbz^{1,s-1/2}(u)$ disjoint
from them, we deduce (using Lemma~\ref{lemma:WFb-mic-q} and its
$H^{-1}$ analogue) that $|\langle (T^-)^*A_rPu,TA_ru\rangle|$
is uniformly bounded. Similarly, taking into account that
$FA_r$ is uniformly bounded
in $\Psib^{-\infty}(X)$, we see that
$|\langle A_rPu,FA_r u\rangle|$ is also uniformly
bounded, so $|\langle A_r Pu,A_r u\rangle|$ is uniformly bounded
for $r\in(0,1]$.

Thus, for some $C_3>0$ depending only on the dimension of $X$,
\begin{equation}\begin{split}\label{eq:prop-64}
\|C\tilde B\Lambda_r u\|^2\leq
&2|\langle A_rPu,A_ru\rangle|+|\langle (E+E')\Lambda_r u,\Lambda_r u\rangle|\\
&\qquad
+\left(6C_2(\delta\ep+\ep^{-1})+C_3\gamma\right)
\|\tilde B\Lambda_r u\|^2
+\gamma^{-1}\|R'_0 \tilde B\Lambda_r u\|^2\\
&\qquad+6C_2(\delta\ep+\ep^{-1})\|\tilde B\Lambda_ru\|
\sum_i \|TD_{x_i}\tilde B\Lambda_ru\|\\
&\qquad
+\gamma^{-1}\sum_i\|T R'_i D_{x_i}
\tilde B\Lambda_ru\|^2+\gamma\|\tilde B\Lambda_ru\|^2\\
&\qquad+\left(6C_2(\delta\ep+\ep^{-1})+C_3\gamma\right)\sum_i\|TD_{x_i}\tilde B\Lambda_ru\|^2\\
&\qquad+\gamma^{-1}\sum_{ij}\|R'_{ij} TD_{x_i}\tilde B\Lambda_r u\|^2\\
&\qquad+\gamma^{-1}\sum_i\|F_i D_{x_i}\tilde B\Lambda_r u\|^2
+\gamma^{-1}\sum_{ij}\|F_{ij} D_{x_i}\tilde B\Lambda_r u\|^2\\
&\qquad+\sum_{ij}
\|R_{ij}D_{x_i}\tilde B\Lambda_r u\|\,\|FD_{x_j}\tilde B\Lambda_r u\|.
\end{split}\end{equation}
All terms but the ones involving $C_2$ or $\gamma$ (not $\gamma^{-1}$)
remain bounded as $r\to 0$.
The $C_2$ and $\gamma$ terms can be estimated by
writing $TD_{x_i}=D_{x_i}T'_i+T_i''$ for some $T'_i,T''_i\in\Psib^{-1}(X)$, and
using Lemma~\ref{lemma:Dt-dX} where necessary, to conclude
that there exist $\gamma>0$, $\ep>0$,
$\delta_0>0$ and $C_4>0$, $C_5>0$ such that for $\delta\in(0,\delta_0)$,
\begin{equation*}\begin{split}
C_4\|\tilde B\Lambda_r u\|^2\leq &2|\im\langle A_rPu,A_ru\rangle|
+|\langle (E+E')\Lambda_r u,\Lambda_r u\rangle|\\
&\qquad+\gamma^{-1}\|R'_0 \tilde B\Lambda_r u\|^2
+C_5\gamma^{-1}\|d_X T^2\tilde B\Lambda_r u\|^2.
\end{split}\end{equation*}
Letting $r\to 0$ now keeps the right hand side
bounded, proving that $\|\tilde B\Lambda_r u\|$ is uniformly bounded as
$r\to 0$, hence $\tilde B\Lambda_0 u\in L^2(X)$ (cf.\ the proof
of Proposition~\ref{prop:elliptic}).
In view of Lemma~\ref{lemma:Dt-dX}
this proves that $q_0\nin\WFbz^{1,s}(u)$, and hence proves the
first statement of the proposition.

In fact,
recalling that we needed $q_0\nin\WFbd^{-1,s+3/2}(Pu)$ for
the uniform boundedness in \eqref{eq:Pu-s+32}, this proves a slightly
weaker version of the second statement of the proposition with
$\WFbd^{-1,s+1}(Pu)$ replaced by $\WFbd^{-1,s+3/2}(Pu)$.
For the more precise statement we modify \eqref{eq:Pu-s+32} -- this
is the only term in \eqref{eq:prop-64} that needs modification
to prove the optimal statement. Let
$\tilde T\in\Psib^{-1/2}(X)$ be elliptic, $\tilde T^-\in\Psib^{1/2}(X)$
a parametrix, $\tilde F=\tilde T^-\tilde T-\Id\in\Psib^{-\infty}(X)$.
Then, similarly to \eqref{eq:Pu-s+32}, we have for any $\gamma>0$,
\begin{equation}\begin{split}\label{eq:Pu-s+1}
|\langle A_r Pu,A_r u\rangle|&\leq|\langle (\tilde T^-)^*A_rPu,\tilde
TA_ru\rangle|
+|\langle A_rPu,\tilde FA_r u\rangle|\\
&\leq \gamma^{-1}\|(\tilde T^-)^*A_r Pu\|_{H^{-1}(X)}^2+
\gamma\|\tilde TA_r u\|_{H^1(X)}^2\\
&\qquad+\|A_r Pu\|_{H^{-1}(X)}\|\tilde FA_r u\|_{H^1(X)}.
\end{split}\end{equation}
The last term on the right hand side can be estimated as before.
As $(\tilde T^-)^*A_r$ is bounded in $\Psibc^{s+1}(X)$ with $\WFb'$
disjoint from $U$, we see that $\|(\tilde T^-)^*A_r Pu\|_{H^{-1}(X)}$
is uniformly bounded. Moreover, $\|d_X\tilde TA\Lambda_r u\|^2$ can
be estimated, using Lemma~\ref{lemma:Dt-dX}, by $\|D_t\tilde TA\Lambda_r u\|^2$
modulo terms that are uniformly bounded as $r\to 0$.
The principal symbol of $D_t\tilde TA$ is $\tau\sigma_{b,-1/2}(\tilde T)a$,
with
$a=\chi_0\chi_1\chi_2$, where $\chi_0$ stands for
$\chi_0(A_0^{-1}(2-\frac{\phi}{\delta}))$, etc., while the
principal symbol $\tilde b$ of $\tilde B$ is given by \eqref{eq:tilde-b},
so we can write:
\begin{equation*}
|\tau|^{1/2}a=|\tau|^{1/2}\chi_0\chi_1\chi_2=A_0^{-1}(2-\phi/\delta)
|\tau|^{1/2}(\chi_0\chi_0')^{1/2}\chi_1\chi_2=
A_0^{-1/2}\delta^{1/2}(2-\phi/\delta)\tilde b, 
\end{equation*}
where we used that
\begin{equation*}
\chi'_0(A_0^{-1}(2-\phi/\delta))
=A_0^2(2-\phi/\delta)^{-2}\chi_0(A_0^{-1}(2-\phi/\delta))
\end{equation*}
when $2-\phi/\delta>0$, while $a$, $\tilde b$ vanish otherwise.
Correspondingly, as $|\tau|^{1/2}\sigma_{b,-1/2}(\tilde T)$ is $\Cinf$,
homogeneous degree zero, near the support of $a$ in $\Tb^*X\setminus o$,
we can write $D_t\tilde TA=G\tilde B+F$, $G\in\Psib^0(X)$, $F\in\Psib^{-1/2}(X)$.
Correspondingly, modulo terms that are bounded
as $r\to 0$, $\|D_t\tilde TA\Lambda_r u\|^2$ (hence
$\|d_X\tilde TA\Lambda_r u\|^2$)
can be estimated from above
by $C_6\|\tilde B\Lambda_r u\|^2$. Thus, modulo terms that
are bounded as $r\to 0$,
for $\gamma>0$ sufficiently small,
$\gamma\|\tilde TA_r u\|_{H^1(X)}^2$ can be absorbed
into $\|C\tilde B\Lambda_r u\|^2$. As the treatment of the other terms
on the right hand side of \eqref{eq:prop-64} requires no change,
we deduce as above
that $\tilde B\Lambda_0 u\in L^2(X)$, which
(in view of Lemma~\ref{lemma:Dt-dX})
proves that $q_0\nin\WFbz^{1,s}(u)$, completing the proof of the iterative
step.

We need to make one more remark to prove the proposition
for $\WFbz^{1,\infty}(u)$, namely
we need to show that the neighborhoods of $q_0$ which are disjoint from
$\WFbz^{1,s}(u)$ do not shrink uncontrollably to $\{q_0\}$ as $s\to\infty$.
This argument parallels to last paragraph of the proof of
\cite[Proposition~24.5.1]{Hor}.
In fact, note that above we have proved that the elliptic set of
$\tilde B=\tilde B_s$
is disjoint from $\WFbz^{1,s}(u)$.
In the next step, when we are proving $q_0\nin\WFbz^{1,s+1/2}(u)$,
we decrease $\delta>0$ slightly (by an arbitrary small amount),
thus decreasing the support
of $a=a_{s+1/2}$ in \eqref{eq:prop-22}, to make sure that $\supp a_{s+1/2}$
is a subset of the elliptic set of the union of $\tilde B_s$ with the
region $\eta<0$, and hence that $\WFbz^{1,s}(u)
\cap\supp a_{s+1/2}=\emptyset$. Each iterative step thus shrinks the
elliptic set of $\tilde B_s$ by an arbitrarily small amount, which allows
us to conclude that $q_0$ has a neighborhood $U'$ such that $\WFbz^{1,s}(u)
\cap U'=\emptyset$ for all $s$. This proves that $q_0\nin \WFbz^{1,\infty}(u)$,
and indeed that $\WFbz^{1,\infty}(u)\cap U'=\emptyset$, for if $A\in\Psib^m(X)$
with $\WFb'(A)\subset U'$ then $Au\in H^1(X)$ by Lemma~\ref{lemma:WFb-mic}
and Corollary~\ref{cor:WF-to-H1}.
\end{proof}

Again, this can be modified to allow Neumann boundary conditions. Namely,
rather than consider $[A_r^*A_r,P]$, we work directly with the quadratic form,
see \eqref{eq:Neumann-solution-def}.
Thus, writing $w=(x,y,t)$ and $\tilde g$ for the semi-Riemannian metric
$g-dt^2$, while $J\,dw$ is the volume form of $g+dt^2$, and $\langle\cdot,
\cdot\rangle$ is the corresponding inner product on $L^2(X)$,
\eqref{eq:Neumann-solution-def} shows that
\begin{equation}\begin{split}\label{eq:Neumann-comm}
&\langle A_r^*A_r u,f\rangle-\langle f,A_r^*A_r u\rangle\\
&\qquad
=\sum_{ij}\langle\tilde g^{ij} D_{w_i}u,D_{w_j}A_r^*A_r u\rangle
-\sum_{ij}\langle\tilde g^{ij} D_{w_i}A_r^*A_r u,D_{w_j}u\rangle.
\end{split}\end{equation}
Then the replacement of \eqref{eq:comm-expansion} is achieved by expanding
the right hand side:
\begin{equation}\begin{split}\label{eq:comm-expansion-Neumann}
&\sum_{ij}\langle\tilde g^{ij} D_{w_i}u,D_{w_j}A_r^*A_r u\rangle
-\sum_{ij}\langle\tilde g^{ij} D_{w_i}A_r^*A_r u,D_{w_j}u\rangle\\
&=\sum_{ij}\langle\tilde g^{ij} D_{w_i}u,[D_{w_j},A_r^*A_r] u\rangle
+\sum_{ij}\langle\tilde g^{ij} D_{w_i}u,A_r^*A_rD_{w_j} u\rangle\\
&\qquad
-\sum_{ij}\langle[\tilde g^{ij} D_{w_i},A_r^*A_r] u,D_{w_j}u\rangle
-\sum_{ij}\langle A_r^*A_r\tilde g^{ij} D_{w_i} u,D_{w_j}u\rangle\\
&=\sum_{ij}\langle\tilde g^{ij} D_{w_i}u,[D_{w_j},A_r^*A_r] u\rangle
-\sum_{ij}\langle[\tilde g^{ij} D_{w_i},A_r^*A_r] u,D_{w_j}u\rangle;
\end{split}\end{equation}
the second and fourth terms in the middle cancel as $A_r^*A_r$ is symmetric.
If there were no boundary present, i.e.\ if $\pa X=\emptyset$,
we could of course write the right
hand side as
\begin{equation*}\begin{split}
&-\sum_{ij}\langle([D_{w_j}^*,A_r^*A_r]\tilde g^{ij} D_{w_i}
+D_{w_j}^*[\tilde g^{ij} D_{w_i},A_r^*A_r])u,u\rangle\\
&\qquad=\langle[D_t^2-\Delta,A_r^*A_r]u,u\rangle,
\end{split}\end{equation*}
so formally this is indeed the same commutator as the one considered
in \eqref{eq:comm-expansion}. The actual expression,
the right hand side of \eqref{eq:comm-expansion-Neumann},
can be analyzed much as in the Dirichlet problem, using
Lemma~\ref{lemma:comm-symbol}
to compute the commutators.

To illustrate the form that \eqref{eq:Neumann-comm}
takes, replace $A_r^*A_r$ by $A^*A$ temporarily, now $\sigma_{\bl,0}(A^*A)=
a^2$. Thus, by Lemma~\ref{lemma:comm-symbol}, up to terms of similar
form with vanishing symbol at $x=0$, $y=y_0$, $t=t_0$, the right hand
side of \eqref{eq:Neumann-comm} is, $\frac{1}{i}$ times,
\begin{equation*}
\int \sum_{ij} g^{ij}D_{x_i} u\,\overline{\tilde C D_{x_j} u}\,J\,dw
+\int \sum_{ij}  g^{ij} \tilde C D_{x_i} u\,\overline{D_{x_j}u}\,J\,dw,
\end{equation*}
where the summation is only over the coordinates vanishing at the corner
(i.e.\ $x_1,\ldots,x_k$), and $\tilde C\in\Psib^{-1}(X)$ with
$\sigma_{\bl,-1}(\tilde C)=|\tau|^{-1}(A_0\delta)^{-1}\chi_0\chi_0'\chi_1^2
\chi_2^2$, cf.\ \eqref{eq:tilde-b} and the sentence afterwards.
We can subtract this from the PDE (which corresponds to restricting
to the characteristic set of $P$, or allowing the term $R'P$ in
\eqref{eq:P-comm}), considered in the form
\begin{equation*}
\int \sum_{ij} \tilde g^{ij}D_{w_i} u\,\overline{D_{w_j} \tilde Cu}\,J\,dw
+\int \sum_{ij}  \tilde g^{ij} D_{w_i} \tilde Cu\,\overline{D_{w_j}u}\,J\,dw,
\end{equation*}
plus terms involving $f$, commute the $C$ through the $D_{w_i}$, $D_{w_j}$
(the commutators are
lower order in terms of b-differential order, so we ignore them), to obtain
an expression for
\begin{equation*}
\int \sum_{ij} g^{ij}D_{\bar y_i} u\,\overline{\tilde C D_{\bar y_j} u}\,J\,dw
+\int \sum_{ij}  g^{ij} \tilde C D_{\bar y_i} u\,
\overline{D_{\bar y_j}u}\,J\,dw,
\end{equation*}
$\bar y=(y,t)$ as usual. Shifting the tangential derivatives $D_{\bar y_i}$
over and rearranging this gives (modulo lower order terms),
with $\tilde B$ as in \eqref{eq:tilde-b}, and $C$
also as there,
\begin{equation*}
\int C\tilde B u\,\overline{C\tilde B u}\,J\,dw=\|C\tilde B u\|^2.
\end{equation*}
The neglected error terms can be treated much as in the Dirichlet problem,
giving the desired
positivity estimate.

\section{Glancing points}\label{sec:glancing}

We again need a technical lemma, roughly stating that when applied to solutions
of $Pu=0$, $u\in H^1_0(X)$, microlocally near $\cG$, $D_{x_i}$
is not merely bounded by $D_t$, but it is small compared to it.
Such an estimate is natural since
$p|_{x=0}=\tau^2-|\xi|^2_y-|\zeta|^2_y$ gives
$\tau^{-2}|\xi|^2\leq C(\tau^{-2}|p|+|x|+|1-\tau^{-2}|\zeta|_y^2|)$, and
$1-\tau^{-2}|\zeta|_y^2$ is homogeneous degree zero and vanishes at $\cG$,
so the right hand size is small near $\cG$.
Below a $\delta$-neighborhood refers to a $\delta$-neighborhood with respect
to the metric associated to any Riemannian metric on the manifold $\Tb^*X$,
and we identify $\Sb^*X$ as the unit ball bundle
with respect to some fibre metric
on $\Tb^*X$.

\begin{lemma}\label{lemma:Dt-Dx}
Suppose $u\in H^1_{0,\loc}(X)$, and suppose that we are given
$K\subset\Sb^*X$ compact satisfying
\begin{equation*}
K\subset\cG\cap T^*\cF_{k,\reg}\setminus\WFbd^{-1,s+1/2}(Pu).
\end{equation*}
Then
there exist $\delta_0>0$ and $C_0>0$ with the following property.
Let $\delta<\delta_0$,
$U\subset\Sb^*X$ open in a $\delta$-neighborhood of $K$,
and $\calA=\{A_r:\ r\in(0,1]\}$ be a bounded
family of ps.d.o's in $\Psibc^s(X)$ with $\WFb'(\calA)\subset U$, and
with $A_r\in\Psib^{s-1}(X)$ for $r\in (0,1]$.

Then
there exist $B\in\Psib^{s-1/2}(X)$, $\tilde B\in\Psib^{s+1/2}(X)$
with $\WFb'(B),\WFb'(\tilde B)\subset U$ and $\tilde C_0
=\tilde C_0(\delta)>0$ such that
for all $r>0$,
\begin{equation*}\begin{split}
\sum_i\|D_{x_i} A_r u\|^2
\leq C_0\delta\|D_t A_r u\|^2+\tilde C_0(&\|u\|^2_{H^1_{\loc}(X)}
+\|Bu\|^2_{H^1(X)}\\
&+\|Pu\|^2_{H^{-1}_{\loc}(X)}
+\|\tilde B Pu\|^2_{H^{-1}(X)}).
\end{split}\end{equation*}
The meaning of $\|u\|_{H^1_{\loc}(X)}$ and
$\|Pu\|^2_{H^{-1}_{\loc}(X)}$ is stated in
Remark~\ref{rem:localize}.
\end{lemma}

\begin{rem}
As $K$ is compact, this is essentially a local result. In particular, we
may assume that $K$ is a subset of $\Tb^*X$ over
a suitable local coordinate patch. Moreover, we may assume that $\delta_0>0$
is sufficiently small so that $D_t$ is elliptic on $U$.
\end{rem}

\begin{proof}
By Lemma~\ref{lemma:Dt-dX}, applied with $K$ replaced by $\WFb'(\calA)$
in the hypothesis (note that the latter is compact), we already know that
\begin{equation}\begin{split}\label{eq:mic-ell-gl-8}
\|d_X A_r u\|^2
\leq &\|D_t A_r u\|^2\\
&\ +C'_0(\|u\|^2_{H^1_{\loc}(X)}+\|Bu\|^2_{H^1(X)}
+\|Pu\|^2_{H^{-1}_{\loc}(X)}
+\|\tilde B Pu\|^2_{H^{-1}(X)}).
\end{split}\end{equation}
for some $C'_0>0$ and for some $B$, $\tilde B$ as in the statement of the
lemma. Thus, we only need to show that if we replace the left hand side
by $\sum_i\|D_{x_i} A_r u\|^2$ (i.e.\ we drop the tangential derivatives,
at least roughly speaking), the constant in front of $\|D_t A_r u\|^2$
can be made small.

As a first step, we freeze the coefficients at $\cF_k$, i.e.\ replace
$A_{ij}(x,y)$, etc., by $A_{ij}(0,y)$. Writing $A_{ij}(x,y)=A_{ij}(0,y)+
\sum x_l A'_{ijl}(x,y)$
as in the proof of Proposition~\ref{prop:elliptic}, we deduce
that if the operators $A_r$ are supported in $|x|<\delta$,
then \eqref{eq:elliptic-28} holds, i.e.
\begin{equation*}
|\int_X \sum x_l A'_{ijl} D_{x_i}A_r u \,\overline{D_{x_j} A_r u}|
\leq C\delta\sum_{i',j'}\|D_{x_{i'}}A_r u\|\,\|D_{x_{j'}}A_r u\|,
\end{equation*}
with analogous estimates with $A_{ij}(x,y)-A_{ij}(0,y)$ replaced
by $B_{ij}(x,y)-B_{ij}(0,y)$ or $C_{ij}(x,y)$. Combined with
\eqref{eq:mic-ell-gl-8} above, this gives that
\begin{equation*}\begin{split}
\int_X&\left(\sum_{ij}A_{ij}(0,y)D_{x_i} A_r u \,\overline{D_{x_j}A_r u}+
\sum_{ij}B_{ij}(0,y)D_{y_i} A_r u \,\overline{D_{y_j}A_r u}\right)\\
&\leq (1+C_1\delta)\|D_t A_r u\|^2\\
&\qquad
+C''_0(\|u\|^2_{H^1_{\loc}(X)}+\|Bu\|^2_{H^1(X)}+\|Pu\|^2_{H^{-1}_{\loc}(X)}
+\|\tilde B Pu\|^2_{H^{-1}(X)}),
\end{split}\end{equation*}
and hence, after rearrangement, that
\begin{equation*}\begin{split}
&\int_X\sum_{ij}A_{ij}(0,y)D_{x_i} A_r u \,\overline{D_{x_j}A_r u}\\
&\qquad\leq \int_X
\left((D_t^2-\sum B_{ij}(0,y)D_{y_i}D_{y_j}) A_r u\,\overline{A_r u}\right)
+C_1\delta\|D_t A_r u\|^2\\
&\qquad\qquad
+C''_0(\|u\|^2_{H^1_{\loc}(X)}+\|Bu\|^2_{H^1(X)}+\|Pu\|^2_{H^{-1}_{\loc}(X)}
+\|\tilde B Pu\|^2_{H^{-1}(X)}).
\end{split}\end{equation*}
It thus suffices to prove that
\begin{equation}\begin{split}\label{eq:Dt-Dx-16}
&\left|\int_X
\left((D_t^2-\sum B_{ij}(0,y)D_{y_i}D_{y_j}) A_r u\,\overline{A_r u}\right)\right|\\
&\qquad\leq C_2\delta \|D_t A_r u\|^2+\tilde C_2(\delta)
(\|u\|^2_{H^1_{\loc}(X)}+\|Bu\|^2_{H^1(X)}),
\end{split}\end{equation}
which we proceed to do.

Let $\psi\in\Cinf(\Sb^*X)$
(which can thus be identified with a homogeneous degree zero
function on $\Tb^*X\setminus o$)
with $\psi\equiv 1$ near $\WFb'(\calA)$,
$\supp \psi\subset U$, $|\psi|\leq 1$, and let $G\in\Psib^0(X)$ be such that
\begin{equation}\begin{split}\label{eq:Dt-Dx-20}
&\WFb'(G)\subset U,
\ \WFb'\left(D_t GD_t-(D_t^2-\sum B_{ij}D_{y_i}D_{y_j})\right)
\cap\WFb'(\calA)=\emptyset\\
&g=\sigma_{b,0}(G)=\psi(1-\tau^{-2}\sum B_{ij}\zeta_i\zeta_j).
\end{split}\end{equation}
Such $\psi$ and $G$ exist, since $D_t$ is elliptic on $\WFb'(\calA)$.
Now,
\begin{equation*}\begin{split}
&\left|\int_X
\left((D_t GD_t-(D_t^2-\sum B_{ij}(0,y)D_{y_i}D_{y_j})) A_r u\,\overline{A_r u}\right)\right|
\leq C_2'
\|u\|^2_{H^1_{\loc}(X)}
\end{split}\end{equation*}
since $(D_t GD_t-(D_t^2-\sum B_{ij}D_{y_i}D_{y_j}))A_r$
is uniformly bounded in $\Psib^{-\infty}(X)$,
by the first line of
\eqref{eq:Dt-Dx-20}.
Moreover,
\begin{equation*}
\sup|g|\leq C_3\delta
\end{equation*}
since $|1-\tau^{-2}\sum B_{ij}\zeta_i\zeta_j|<C_3\delta$ on
a $\delta$-neighborhood of $K$. Indeed,
$1-\tau^{-2}\sum B_{ij}\zeta_i\zeta_j$ is a homogeneous degree zero
$\Cinf$ function on
a neighborhood of $K$ in $\Tb^*X$ (hence $\Cinf$ near $K$ in $\Sb^*X$) which
vanishes at $\cG\cap T^*\cF_k$. Since there exists
$G'\in\Psib^{-1}(X)$ with $\WFb'(G')\subset U$
satisfying
\begin{equation*}
\|Gv\|\leq2\sup|g|\,\|v\|+\|G'v\|
\end{equation*}
for all $v\in L^2(X)$, we deduce that $\|Gv\|\leq2 C_3 \delta\|v\|+\|G'v\|$
for all $v\in L^2(X)$. Applying this with $v=D_t A_r u$, and estimating
$\|G'v\|$ using Lemma~\ref{lemma:WFb-mic-q},
\eqref{eq:Dt-Dx-16} follows,
which in turn completes the proof of the lemma.
\end{proof}

We are now ready to state and prove the tangential propagation estimate.
First, local coordinates $(x,y,t)$ near $p\in \cF_{i,\reg}$
give a product decomposition of a neighborhood
of $p\in \cF_{i,\reg}$ in $X$ of the form $U\times V$, $U\subset[0,\infty)^k$,
$V\subset\Real^{l+1}$, hence of $T^*X$ as $T^*U\times T^*V$. We denote the
projection $T^*X\to T^*V$ by $\pi^e_i$. Explicitly,
in local coordinates $(x,y,t,\xi,\zeta,\tau)$ on $T^*X$,
\begin{equation*}
\pi^e_i(x,y,t,\xi,\zeta,\tau)=(y,t,\zeta,\tau).
\end{equation*}
With
$\pi_i:T^*_{\cF_{i,\reg}}X\to\dot\Tb^*X$ being the restriction of $\pi$
to $T^*_{\cF_{i,\reg}}X$, $\pi^e_i$ is an extension of $\pi_i$ in the sense that
$\pi^e_i|_{T^*_{\cF_{i,\reg}}X\cap (T^*U\times T^*V)}=\pi_i$.
The tangential propagation estimate is then the following:

\begin{prop}\label{prop:tgt-prop}
Let $u\in H^1_{0,\loc}(X)$. Given $K\subset\Sb^*X$ compact with
\begin{equation}
K\subset(\cG\cap T^*\cF_{i,\reg})\setminus\WFbd^{-1,\infty}(Pu),
\end{equation}
there exist constants $C_0>0$,
$\delta_0>0$ such that the following holds. If
$q_0=(y_0,t_0,\zeta_0,\tau_0)\in K$ and
for some $0<\delta<\delta_0$, $C_0\delta\leq\epsilon<1$ and
for all $\alpha=(x,y,t,\xi,\zeta,\tau)\in \Char(P)$
\begin{equation}\begin{split}\label{eq:tgt-prop-est}
\alpha\in T^*\cF_{j,\reg} &\Mand
|\pi_i^e(\alpha-\exp(-\delta H_p)(\pih^{-1}(q_0)))|
\leq\epsilon\delta\Mand
|x(\alpha)|\leq\epsilon\delta\\
&\Rightarrow
\pi_j(\alpha)\nin\WFb(u),
\end{split}\end{equation}
then $q_0\nin\WFb(u)$.
\end{prop}

\begin{rem}\label{rem:tgt-est}
In the estimate \eqref{eq:tgt-prop-est}, $H_p$ can be replaced by any
$\Cinf$ vector field which agrees
with $H_p$ at the point $\pih^{-1}(q_0)$, since flow to distance $\delta$
along a vector field only depends on the vector field evaluated
at the initial
point of the flow, up to committing an error $\calO(\delta^2)$.
In particular, it can be replaced by the vector field $W^\flat$ defined
below. Similarly,
changing the initial point of the flow by $\calO(\delta^2)$ will not
affect the endpoint up to an error $\calO(\delta^2)$.
Thus, estimate \eqref{eq:tgt-prop-est} can be further
rewritten, at the cost of changing
$C_0$ again, as
\begin{equation}\begin{split}\label{eq:tgt-prop-est-3}
\alpha\in T^*\cF_{j,\reg} &\Mand
|\pi_i^e(\exp(\delta W^\flat)(\alpha))-\xi_0|\leq\epsilon\delta\Mand
|x(\exp(\delta W^\flat)(\alpha))|\leq\epsilon\delta\\
&\Rightarrow
\pi_j(\alpha)\nin\WFb(u);
\end{split}\end{equation}
here we also interchanged the roles of the intial and final points of the flow.
\end{rem}

\begin{proof}
The proof is very similar to the previous one and now
the positive commutator
construction follows that of Melrose and Sj\"ostrand
\cite{Melrose-Sjostrand:I}, as well as \cite{Vasy:Propagation-Many}
in $N$-body scattering without bound states.
Thus, we take local coordinates as above, i.e.\ of the form
$(x,y,t)$ with the $\cF_j$ intersecting the coordinate neighborhood
defined by the vanishing of components of $x$. We can use
$t-t_0$ now to measure propagation, since $\tau^{-1}H_p(t-t_0)=2>0$.
More precisely, to allow for both signs of $\tau$ and yet keep the sign
of the derivative along $H_p$ fixed, we need to take
\begin{equation*}
\tilde\eta=(\sign\tau) (t-t_0)
\end{equation*}
as the propagation variable, so $|\tau|^{-1}H_p\tilde\eta=2$. However,
for the sake of notational simplicity and clarity, we take $\tau_0>0$,
and make all symbols below supported in $\tau>0$
-- the general setting only requires replacing $t-t_0$ by $\tilde\eta$
in \eqref{eq:glancing-phi-def} below.

Then we could
construct $\omega_0\in\Cinf(T^*\cF_i)$ (defined near $q_0$) to
measure
the squared distance from the integral curve of
\begin{equation}
W^\flat=2\tau\pa_t-H_h,\ h(y,\zeta)=\zeta\cdot B(y)\zeta
\end{equation}
through $q_0$;
this can be achieved by solving a Cauchy problem as in
\cite{Melrose-Sjostrand:I},
\cite{Vasy:Propagation-Many}.
In fact, this does not need to be done
precisely -- after all, $W^\flat$ is only an approximation to $H_p$
in the very first place.
Thus, all we need is that $\omega_0$ is the
sum of squares of $2l$ homogeneous degree zero functions $\rho_j$:
\begin{equation*}
\omega_0=\sum_{j=1}^{2l} \rho_j^2,\ W^\flat\rho_j(q_0)=0,
\end{equation*}
$d\rho_j(q_0)$, $j=1,\ldots,2l$ linearly independent at $q_0$. Since
$\dim \cF_j=l+1$, $d\rho_j(q_0)$, $j=1,\ldots,2l$, together with $dt$
($t$ is also homogeneous degree zero), span the cotangent space of
the quotient of $T^*\cF_i$ by the $\Real^+$-action, for dimensional reasons
(note that $W^\flat t(q_0)\neq 0$). In particular,
\begin{equation*}
|\tau^{-1}W^\flat\omega_0|\leq C_1'\omega_0^{1/2}(\omega_0^{1/2}+|t-t_0|)
\end{equation*}
Then we
extend $\omega_0$ to a function on $\Tb^*X$ (using the coordinates
$(x,y,t,\sigma,\zeta,\tau)$),
let
\begin{equation}
\omega=\omega_0+|x|^2.
\end{equation}
Then the `naive' estimate, playing an analogous role
to \eqref{eq:omega-est-a0} in the
hyperbolic region, is
\begin{equation}\begin{split}\label{eq:omega-est-a0-t}
|\tau^{-1}H_p\omega|&\leq \tilde C_1''\omega^{1/2}(\omega^{1/2}+|t-t_0|
+\tau^{-2}|\xi|^2)\\
&\leq C_1''\omega^{1/2}(\omega^{1/2}+|t-t_0|+\tau^{-2}|p|),
\end{split}\end{equation}
where we used that $p|_{x=0}=\tau^2-|\xi|^2_y-|\zeta|^2_y$ lets us
estimate
\begin{equation*}
\tau^{-2}|\xi|^2\leq C(\tau^{-2}|p|+|x|+\omega_0^{1/2}+|t-t_0|),
\end{equation*}
for
$1-\tau^{-2}|\zeta|_y^2$ is homogeneous degree zero and vanishes at $\cG$
(recall from the beginning of the section that this last estimate
motivates Lemma~\ref{lemma:Dt-Dx}).
Note that \eqref{eq:omega-est-a0-t}
is much more precise than \eqref{eq:omega-est-a0}:
we have a factor of $\omega^{1/2}+|t-t_0|+\tau^{-2}|p|$ in addition to
$\omega^{1/2}$ -- this is crucial since we need to get the direction
of propagation right. Again, we in fact need a more explicit version of
this:
\begin{equation}\begin{split}\label{eq:omega-est-a-t}
&\tau^{-1}H_p\omega=f_0+\sum_i f_i\tau^{-1}\xi_i
+\sum_{i,j}f_{ij}\tau^{-2}\xi_i\xi_j,\\
&\ f_i,f_{ij}\in\Cinf(\Tb^*X),\ |f_i|
\leq C_1\omega^{1/2}(\omega^{1/2}+|t-t_0|),\ |f_{ij}|\leq C_1\omega^{1/2}
\end{split}\end{equation}
$f_i$, $f_{ij}$ homogeneous of degree $0$.
Note that the estimates on $f_{ij}$ are weaker than the estimates on $f_i$.
In fact, $f_{ij}$ arises from the
$2\sum (\pa_{y_k} A_{ij})\xi_i\xi_j\pa_{\zeta_k}$ term of $H_p$ in
\eqref{eq:H_p-form} -- when applied to $\rho_j^2$, it gives a result of the
stated form. The reason for the sufficiency of this weaker estimate
is that at $\pih^{-1}(q_0)$, $\xi=0$, so the $f_{ij}$ term can be
estimated using $P$ (as will be done below),
as was already done at a formal level in
\eqref{eq:omega-est-a0-t}.

Finally, we let
\begin{equation}\label{eq:glancing-phi-def}
\phi=t-t_0+\frac{1}{\ep^2\delta}\omega,
\end{equation}
and define $a$ almost as in \eqref{eq:prop-22}, with $\eta$
replaced by $t-t_0$, namely
\begin{equation}\label{eq:prop-22t}
a=\chi_0(A_0^{-1}(2-\phi/\delta))\chi_1((t-t_0+\delta)/
\ep\delta+1)\chi_2(|\sigma|^2/\tau^2).
\end{equation}
The slight difference is in the argument of $\chi_1$, in order to
microlocalize more precisely in the `hypothesis region', i.e.\ where
$u$ is a priori assumed to have no wave front set. This is natural, since
for the hyperbolic points we only needed to prove that singularities
cannot stay at the given boundary face $\cF_{i,\reg}$, while
for glancing points we need to get the correct direction of
propagation.
We always assume $\ep<1$, so on $\supp a$ we have
\begin{equation*}
\phi\leq 2\delta\Mand t-t_0\geq-\ep\delta-\delta\geq-2\delta.
\end{equation*}
Since $\omega\geq 0$, the first of these inequalities implies that
$t-t_0\leq 2\delta$, so on $\supp a$
\begin{equation}
|t-t_0|\leq 2\delta.
\end{equation}
Hence,
\begin{equation}\label{eq:omega-delta-est-t}
\omega\leq \ep^2\delta(2\delta-(t-t_0))\leq4\delta^2\ep^2.
\end{equation}
Moreover, on $\supp d\chi_1$,
\begin{equation}\label{eq:dchi_1-supp}
t-t_0\in[-\delta-\ep\delta,-\delta],\ \omega^{1/2}\leq 2\ep\delta,
\end{equation}
so this region lies in \eqref{eq:tgt-prop-est-3} after $\ep$ and $\delta$
are both replaced by appropriate constant multiples, namely the present
$\delta$ should be replaced by $\delta/2\tau_0$.

We again start with the imprecise motivational argument.
Thus, using \eqref{eq:omega-est-a0-t}, \eqref{eq:omega-delta-est-t},
$\tau^{-1}H_p(t-t_0)=2=c_0>0$, we deduce that
at $p=0$,
\begin{equation*}\begin{split}
\tau^{-1} H_p\phi&=H_p(t-t_0)+\frac{1}{\ep^2\delta}H_p\omega\\
&\geq c_0/2
-\frac{1}{\ep^2\delta}C_1''\omega^{1/2}(\omega^{1/2}+|t-t_0|)\\
&\geq c_0/2-2C_1''(\delta+\frac{\delta}{\ep})\geq
c_0/4>0
\end{split}\end{equation*}
provided that $\delta<\frac{c_0}{16 C_1''}$,
$\frac{\ep}{\delta}>\frac{16C_1''}{c_0}$, i.e.\ that $\delta$ is small, but
$\ep/\delta$ is not too small -- roughly, $\ep$ can go to $0$ at most
proportionally to $\delta$ (with an appropriate constant) as $\delta\to 0$.
(Recall also that $\ep<1$, so there is an upper bound as well for $\ep$,
but this is of no significance as we let $\delta\to 0$. It is also
worth remembering that in the hyperbolic region, $\ep$ roughly played the
same role as here, but was bounded below by an absolute constant, rather
than by a suitable multiple of $\delta$, hence could not go to $0$ as
$\delta\to 0$.)
With this, we can proceed exactly as in the
hyperbolic region, so (recall that $\tau>0$ on $\supp a$!)
\begin{equation*}
H_p a^2=-b^2+e,\ b=\tau^{1/2} (2\tau^{-1} H_p\phi)^{1/2}
(A_0\delta)^{-1/2}(\chi_0\chi_0')^{1/2}\chi_1\chi_2,
\end{equation*}
with $e$ arising from the derivative of $\chi_1\chi_2$. Again,
$\chi_0$ stands for $\chi_0(A_0^{-1}(2-\frac{\phi}{\delta}))$, etc.
In view of \eqref{eq:dchi_1-supp} and \eqref{eq:tgt-prop-est-3} on the
one hand, and that
$d\chi_2$ is disjoint from the characteristic set on the other,
both $\supp d\chi_1$ and $\supp d\chi_2$ are disjoint from $\WFb(u)$. Thus,
$i[A^*A,P]$ is positive modulo terms that we can a priori control,
so the standard positive commutator argument gives an estimate
for $Bu$, where $B$ has symbol $b$. Replacing $a$ by $a\tau^{s+1/2}$,
we still have a positive commutator (again, $D_t$ actually commutes
with $P$, but in any case we could use $A_0$ to bound the additional
commutator term), which now gives (with the new $B$) that $Bu\in L^2(X)$,
which means in particular that $q_0\nin\WFbz^{1,s}(u)$.

The detailed proof is analogous to the hyperbolic case, with the biggest
difference being the treatment of the $f_{ij}$ term in $\tau^{-1}H_p\omega$.
First,
\begin{equation}\begin{split}\label{eq:H_p-phi-8}
\tau^{-1} H_p\phi
&=\tau^{-1}H_p(t-t_0)+\frac{1}{\ep^2\delta}\tau^{-1}H_p\omega\\
&=2+\frac{1}{\ep^2\delta}(f_0+\sum_i f_i\tau^{-1}\xi_i
+\sum_{i,j}f_{ij}\tau^{-2}\xi_i\xi_j).
\end{split}\end{equation}
Let $\tilde B\in\Psib^{1/2}(X)$ with
\begin{equation*}
\tilde b=\sigma_{b,0}(\tilde B)=\tau^{1/2}
(A_0\delta)^{-1/2}(\chi_0\chi_0')^{1/2}\chi_1\chi_2\in\Cinf(\Tb^*X\setminus o),
\end{equation*}
and let $A\in\Psib^0(X)$ with $\sigma_{b,0}(A)=a$. Again,
$\chi_0$ stands for $\chi_0(A_0^{-1}(2-\frac{\phi}{\delta}))$, etc.
Also, let $C\in\Psib^0(X)$
have symbol $\sigma_{b,0}(C)=\sqrt{2}\,\psi$ where
$\psi\in S^0(\Tb^*X)$ is identically $1$ on $U$ considered as a subset
of $\Tb^*X$.
Then an explicit calculation using Lemma~\ref{lemma:comm-symbol} gives,
in accordance with \eqref{eq:H_p-phi-8},
\begin{equation*}\begin{split}
&i[A^*A,P]\\
&\quad=R' P+\tilde B^*(C^*C+R_0+\sum_i D_{x_i} R_i
+\sum_{ij} D_{x_i}R_{ij}D_{x_j})\tilde B+R''+E+E'
\end{split}\end{equation*}
with
\begin{equation*}\begin{split}
&R_0\in\Psib^0(X),\ R_i\in\Psib^{-1}(X),\ R_{ij}\in\Psib^{-2}(X),\\
&R'\in\Psib^{-1}(X),\ R''\in\Diff^2\Psib^{-2}(X),
\ E,E'\in\Diff^2\Psib^{-1}(X),
\end{split}\end{equation*}
with $\WFb'(E)\subset\eta^{-1}((-\infty,-\delta])
\cap U$, $\WFb'(E')\cap\dot\Sigma=\emptyset$
($E$ arises from the commutator of $P$ with an operator with symbol
$\chi_1(\eta/\delta+2)$, while $E'$ from the commutator of $P$ with an operator
with symbol $\chi_2(|\sigma|^2/\tau^2)$)
and with
$r_0=\sigma_{b,0}(R_0)$, $r_i=\sigma_{b,-1}(R_i)$,
$r_{ij}\in\sigma_{b,-2}(R_{ij})$,
\begin{equation*}
|r_0|\leq \frac{C_2}{\ep^2\delta}\omega^{1/2}(|t-t_0|+\omega^{1/2}),
\ |\tau r_i|\leq \frac{C_2}{\ep^2\delta}\omega^{1/2}
(|t-t_0|+\omega^{1/2}),
\ |\tau^2 r_{ij}|\leq \frac{C_2}{\ep^2\delta}\omega^{1/2},
\end{equation*}
and $\supp r_j$ lying in $\omega^{1/2}\leq 3\ep\delta$, $|t-t_0|<3\delta$.
Thus,
\begin{equation*}
|r_0|\leq 3C_2(\delta+\frac{\delta}{\ep}),
\ |\tau r_i|\leq 3C_2(\delta+\frac{\delta}{\ep}),
\ |\tau^2 r_{ij}|\leq 3C_2\ep^{-1}.
\end{equation*}
Thus, the $R_0$ and $R_i$ terms can be treated exactly as in the
hyperbolic case, i.e.\ as in the proof of Proposition~\ref{prop:normal-prop}.
That is, as in the hyperbolic setting,
let $T\in\Psib^{-1}(X)$ be elliptic, $T^-\in\Psib^1(X)$ a parametrix, so
$T^-T=\Id+F$, $F\in\Psib^{-\infty}(X)$.
Then there exist $R'_0,R'_i\in\Psib^{-1}(X)$ such that for any $\gamma>0$,
\begin{equation*}\begin{split}
|\langle R_0 v,v\rangle|&\leq \|R_0 v\|\,\|v\|
\leq 2\sup |r_0|\,\|v\|^2+\|R_0' v\|\,\|v\|\\
&\leq 6C_2(\frac{\delta}{\ep}+\delta)\|v\|^2
+\gamma^{-1}\|R_0' v\|^2+\gamma \|v\|^2,
\end{split}\end{equation*}
\begin{equation*}\begin{split}
\|R_i w\|=\|R_i (T^- T -F)w\|&\leq\|(R_i T^-)(Tw)\|+\|R_iFw\|\\
&\leq 6C_2(\frac{\delta}{\ep}+\delta)\|Tw\|+
\|R_i' Tw\|+\|R_iFw\|
\end{split}\end{equation*}
for all $w$ with $Tw\in L^2(X)$, hence
\begin{equation*}\begin{split}
|\langle R_i D_{x_i} v,v\rangle|\leq
&6C_2(\frac{\delta}{\ep}+\delta)\|TD_{x_i} v\|\,\|v\|\\
&\qquad+2\gamma\|v\|^2+\gamma^{-1}\|R'_i TD_{x_i} v\|^2+\gamma^{-1}\|F_i D_{x_i}v\|^2,
\end{split}\end{equation*}
with

However, the $R_{ij}$ term needs to be treated separately, since we
need that microlocally
$\tau^{-1}D_{x_i}$ is small (bounded by a constant multiple
of $\delta$), and not merely bounded,
which is all we needed both in the proof of Proposition~\ref{prop:normal-prop}
and here for the $R_0$ and $R_i$ terms. This is accomplished
by the use of Lemma~\ref{lemma:Dt-Dx}. Namely, as in the hyperbolic
setting, there exist $R'_{ij}\in\Psib^{-1}(X)$ such that
\begin{equation*}
\|(T^-)^*R_{ij} w\|\leq 6C_2\ep^{-1}\|Tw\|+
\|R_{ij}' Tw\|+\|(T^-)^*R_{ij}Fw\|
\end{equation*}
for all $w$ with $Tw\in L^2(X)$.
Thus,
\begin{equation*}\begin{split}
|\langle R_{ij} D_{x_i} v,D_{x_j}v\rangle|&\leq
6C_2\ep^{-1}\|TD_{x_i}v\|\,\|TD_{x_j}v\|\\
&\qquad\qquad
+\gamma\|TD_{x_j}v\|^2+\gamma^{-1}\|R'_{ij}T D_{x_i}v\|^2
+\gamma^{-1}\|F_{ij} D_{x_i}v\|^2\\
&\qquad\qquad+\|R_{ij}D_{x_i}v\|\,\|FD_{x_j}v\|,
\end{split}\end{equation*}
with $F_{ij}\in\Psib^{-\infty}(X)$.
For $v=\tilde B_r u$, $\tilde B_r=\tilde B\Lambda_r$,
Lemma~\ref{lemma:Dt-Dx} thus gives
\begin{equation*}\begin{split}
|\langle R_{ij} D_{x_i} \tilde B_r u,D_{x_j}\tilde B_r u\rangle|&\leq
6C'_2\frac{\delta}{\ep}\|\tilde B_r u\|^2+\gamma\|\tilde B_r u\|^2\\
&\qquad\qquad+\gamma^{-1}\|R'_{ij}T D_{x_i}\tilde B_ru\|^2
+\gamma^{-1}\|F_{ij} D_{x_i}\tilde B_ru\|^2\\
&\qquad\qquad+\|R_{ij}D_{x_i}\tilde B_ru\|\,\|FD_{x_j}\tilde B_r u\|.
\end{split}\end{equation*}
For $\delta<\delta_0$, $\frac{\delta}{\ep}<C'_0$ sufficiently small,
we finish the
proof as in the hyperbolic setting, showing that
$\tilde B\Lambda_0 u\in L^2(X)$, and hence that $q_0\nin\WFbz^{1,s}(u)$.

Again, \eqref{eq:prop-22t} needs to be modified slightly to show
$q_0\nin\WFbz^{1,\infty}(u)$. Now we take, with $\nu\leq 1$,
\begin{equation*}
a=\chi_0(A_0^{-1}(1+\nu-\phi/\delta))\chi_1((t-t_0+\delta)/
\ep\delta+\nu)\chi_2(|\sigma|^2/\tau^2),
\end{equation*}
i.e.\ we replace $2$ by $1+\nu$ in in the argument of $\chi_0$, and
we replace $1$ by $\nu$ in the argument of $\chi_1$. In the iterative
step we decrease $\nu$ by an arbitrarily small amount, which suffices
to prove $q_0\nin\WFbz^{1,\infty}(u)$; see also the proof of
Proposition~\ref{prop:normal-prop} here,
and the proof of \cite[Proposition~24.5.1]{Hor}.
\end{proof}

The results of this section can be adapted to Neumann boundary conditions,
using the argument presented at the end of the previous section.

\section{Propagation of singularities}\label{sec:prop-sing}
An argument of Melrose and Sj\"ostrand \cite{Melrose-Sjostrand:I,
Melrose-Sjostrand:II}, see also \cite[Chapter~XXIV]{Hor} and
\cite{Lebeau:Propagation} allows us to conclude our main result
concerning the singularities of solutions of the wave equation.
The proof presented below essentially follows Lebeau's
paper \cite[Proposition~VII.1]{Lebeau:Propagation}. Correspondingly,
we only give the proof at $\cH$ in full detail; at $\cG$ the arguments
are sketched, but the details are {\em precisely} as in Lebeau's case.
We mostly discuss the Dirichlet boundary condition -- the results are
also valid for Neumann
boundary conditions, see Theorem~\ref{thm:prop-sing-N},
and the arguments presented need no modification at all in that case.
We thus have the following theorem.

\begin{thm}\label{thm:prop-sing}
Suppose that $u\in H^1_{0,\loc}(X)$. Then
$\WFbz^{1,\infty}(u)\setminus\WFbd^{-1,\infty}(Pu)\subset\dot\Sigma$, and it
is a union of maximally extended
generalized broken bicharacteristics of $P$ in $\dot\Sigma\setminus\WFbd^{-1,\infty}(Pu)$.

In fact, if $u\in H^{1,m}_{0,\loc}(X)$ for some $m\leq 0$,
then for all $s\in\Real\cup\{\infty\}$,
$\WFbz^{1,s}(u)\setminus\WFbd^{-1,s+1}(Pu)\subset\dot\Sigma$, and it
is a union of maximally extended
generalized broken bicharacteristics of $P$ in
$\dot\Sigma\setminus\WFbd^{-1,s+1}(Pu)$.
\end{thm}

\begin{rem}\label{rem:inhom-Dir}
Suppose that for each boundary hypersurface $H_j$, we
are given Dirichlet data $g_j\in\Cinf(H_j)$, which are compatible,
so at $H_i\cap H_j$, $g_i|_{H_i\cap H_j}
=g_j|_{H_i\cap H_j}$ for all $i,j$. Then there is $g\in\Cinf(X)$ with
$g|_{H_j}=g_j$. Now, if $u\in H^1_{\loc}(X)$ and $u|_{H_j}=g_j$, then
$v=u-g\in H^1_{0,\loc}(X)$. Thus, the theorem is applicable to $v$.
Since $Pv=Pu-Pg$ and $Pg\in\Cinf(X)$, $\WFbd^{-1,\infty}(Pu)=\WFbd^{-1,\infty}(Pv)$,
and similarly $\WFbz^{1,\infty}(u)=\WFbz^{1,\infty}(v)$, we
deduce that $\WFbz^{1,\infty}(u)\setminus\WFbd^{-1,\infty}(Pu)$
is a union of maximally extended
generalized broken bicharacteristics of $P$ in
$\dot\Sigma\setminus\WFbd^{-1,\infty}(Pu)$.
\end{rem}

\begin{rem}\label{rem:sing-soln}
As already expained in the introduction, we can relax the hypothesis
$u\in H^1_{0,\loc}(X)$ in the results of
Sections~\ref{sec:elliptic}-\ref{sec:glancing} to
$u\in H^{1,m}_{\bl,0,\loc}(X)$, $m\leq 0$ without changing the arguments,
except replacing the $H^1_{\loc}(X)$ norms by
the $H^{1,m}_{\bl,\loc}$ norms for the
`background terms', such as $\|u\|_{H^1_\loc(X)}$ in Lemma~\ref{lemma:Dt-dX},
and analogously for $\|Pu\|_{H^{-1}_\loc(X)}$. The microlocal norms,
in which we are gaining regularity, such as those of $Bu$ and $\tilde B Pu$
in Lemma~\ref{lemma:Dt-dX} are {\em unchanged!} Indeed, now
we merely need to apply Lemma~\ref{lemma:WFb-mic-q-p} in place of
Lemma~\ref{lemma:WFb-mic-q}.

The point of this generalization is to allow more singular (approximate)
solutions of the wave equation, such as its fundamental solution.
An alternative way to deal with these solutions is to regularize them
in time (which one can do without destroying, say, $Pu=0$),
and use the $H^1_{0,\loc}(X)$ results -- but stating (and proving) the
result for $u\in H^{1,m}_{\bl,0,\loc}(X)$ is the neater way to proceed.
\end{rem}

\begin{cor}\label{cor:prop-sing}
Suppose that $Pu=0$, $u\in H^1_{0,\loc}(X)$. Then $\WFb(u)\subset\dot\Sigma$,
and it
is a union of maximally extended
generalized broken bicharacteristics of $P$ in $\dot\Sigma$.
\end{cor}

The theorem for Neumann boundary conditions takes the following form.

\begin{thm}\label{thm:prop-sing-N}
Suppose that $u\in H^1_{\loc}(X)$ and $f\in \dot H^{-1}_{\loc}(X)$. Suppose
also that for all $v\in H^1_{c}(X)$,
\begin{equation}\label{eq:Neumann-def}
\langle D_t u,D_t v\rangle-\langle d_M u,d_M v\rangle=\langle f,v\rangle.
\end{equation}
Then
$\WFbz^{1,s}(u)\setminus\WFbd^{-1,s+1}(f)\subset\dot\Sigma$, and it
is a union of maximally extended
generalized broken bicharacteristics of $P$ in $\dot\Sigma\setminus\WFbd^{-1,s+1}(f)$.

In fact, if $u\in H^{1,m}_{\loc}(X)$ for some $m\leq 0$, and
\eqref{eq:Neumann-def} holds for all $v\in H^{1,-m}_{c}(X)$
then for all $s\in\Real\cup\{\infty\}$,
$\WFbz^{1,s}(u)\setminus\WFbd^{-1,s+1}(f)\subset\dot\Sigma$, and it
is a union of maximally extended
generalized broken bicharacteristics of $P$ in
$\dot\Sigma\setminus\WFbd^{-1,s+1}(f)$.
\end{thm}

\begin{proof}(Proof of Theorem~\ref{thm:prop-sing}.)
For notational simplicity, we state the proof for $\WFbz^{1,\infty}(u)$.
The case of general $s$ only requires notational changes. Note
that $\WFbz^{1,\infty}(u)\setminus\WFbd^{-1,\infty}(Pu)\subset\dot\Sigma$
by Proposition~\ref{prop:elliptic}, so we only need to prove that it is
a union of maximally extended
generalized broken bicharacteristics of $P$ in $\dot\Sigma\setminus\WFbd^{-1,\infty}(Pu)$.

We start by remarking that for every $V\subset\dot\Sigma$
and $q\in V$, the set $\calR$
of generalized broken bicharacteristics $\gamma$ defined on open intervals
including $0$, satisfying $\gamma(0)=q$, and with image in $V$, has
a natural partial order, namely if $\gamma:(\alpha,\beta)\to V$,
$\gamma':(\alpha',\beta')\to V$, then $\gamma\leq \gamma'$ if
the domains satisfy $(\alpha,\beta)\subset(\alpha',\beta')$ and
$\gamma=\gamma'|_{(\alpha,\beta)}$. Moreover, any non-empty
totally ordered subset has an upper bound: one can take the generalized broken
bicharacteristic with domain given by the union of the domains of those
in the totally ordered subset, and which extends these, as an upper bound.
Hence, by Zorn's lemma, if $\calR$ is not empty, it has a maximal element.
Note that we can also work with intervals of the form
$(\alpha,0]$, $\alpha<0$, instead of
open intervals.

We only need to prove that for every
$q_0\in\WFb^{1,\infty}(u)\setminus\WFbd^{-1,\infty}(Pu)$
there exists a generalized broken bicharacteristic
$\gamma:[-\ep_0,\ep_0]\to\dot\Sigma$, $\ep_0>0$,
with $\gamma(0)=q_0$ and such that
$\gamma(t)\in\WFb^{1,\infty}(u)\setminus\WFbd^{-1,\infty}(Pu)$ for $t\in[-\ep_0,\ep_0]$.
In fact, once this statement is shown, taking
$V=\WFb^{1,\infty}(u)\setminus\WFbd^{-1,\infty}(Pu)$, $q=q_0$, in the argument of the
previous paragraph, we see that $\calR$ is non-empty, hence has a maximal
element. We need to show that such an element, $\gamma:(\alpha,\beta)\to
\dot\Sigma$, is maximal in
$\dot\Sigma\setminus\WFbd^{-1,\infty}(Pu)$ as well, i.e.\ if we
take $V=\dot\Sigma\setminus\WFbd^{-1,\infty}(Pu)$, $q=q_0$
in the first paragraph. But if
$\gamma':(\alpha',\beta')\to
\dot\Sigma$ is any proper extension of $\gamma$, with say $\alpha'<\alpha$,
with image in $\dot\Sigma\setminus\WFbd^{-1,\infty}(Pu)$, then
$\gamma'(\alpha)\in\WFb^{1,\infty}(u)$ since $\WFb^{1,\infty}(u)$ is closed, and $\gamma$
maps into it, hence by our assumption there is a generalized broken
bicharacteristic $\gammat:(\alpha-\ep',\alpha+\ep')\to\WFb^{1,\infty}(u)
\setminus\WFbd^{-1,\infty}(Pu)$, $\ep'>0$, $\gammat(\alpha)=\gamma'(\alpha)$;
piecing together $\gammat|_{(\alpha-\ep',\alpha]}$ and $\gamma$,
directly from Definition~\ref{Def:gen-br-bichar}, gives a generalized
broken bicharacteristic which is a proper extension of $\gamma$,
with image in $\WFb^{1,\infty}(u)\setminus\WFbd^{-1,\infty}(Pu)$, contradicting
the maximality of $\gamma$.

Indeed, it suffices to show that
for any $i$, if
\begin{equation}\label{eq:prop-103}
q_0\in\WFb^{1,\infty}(u)\setminus\WFbd^{-1,\infty}(Pu)\Mand q_0\in T^*\cF_{i,\reg}
\end{equation}
then
\begin{equation}\label{eq:prop-104}\begin{split}
&\text{there exists a generalized broken bicharacteristic}
\ \gamma:[-\ep_0,0]\to\dot\Sigma,\ \ep_0>0,\\
&\qquad\qquad \gamma(0)=q_0,
\ \gamma(t)\in\WFb^{1,\infty}(u)\setminus\WFbd^{-1,\infty}(Pu),\ t\in[-\ep_0,0],
\end{split}\end{equation}
for the existence of a generalized broken bicharacteristic
on $[0,\ep_0]$ can be demonstrated similarly
by replacing the forward propagation
estimates by backward ones, and, directly from
Definition~\ref{Def:gen-br-bichar}, piecing together the two
generalized broken bicharacteristics gives one defined on $[-\ep_0,\ep_0]$.

We proceed to prove that \eqref{eq:prop-103} implies
\eqref{eq:prop-104} by induction on $i$. For $i=0$, this is certainly true
by H\"ormander's theorem on propagation of singularities, and if
$\codim \cF_i=1$, it follows from the Melrose-Sj\"ostrand theorem.

So suppose that \eqref{eq:prop-103}$\Rightarrow$\eqref{eq:prop-104}
has been proved for all $j$ with
$\cF_i\subsetneq \cF_j$ and that $q_0\in\cH\cap T^*\cF_{i,\reg}$
satisfies \eqref{eq:prop-103}.
We use the notation of the proof of
Proposition~\ref{prop:normal-prop} below.
Let $U\subset\cup_{\cF_i\subset \cF_j}T^*\cF_{j,\reg}$
be a neighborhood of $q_0=(0,y_0,t_0,\zeta_0,\tau_0)$ in
$\dot\Sigma$ which is given
by equations of the form $|x|<\delta'$, $|y-y_0|<\delta'$, $|t-t_0|<\delta'$,
$|\tau-\tau_0|<\delta'$, $|\zeta-\zeta_0|<\delta'$, $\delta'>0$,
such that $H_p\eta>0$
on $\pih^{-1}(U)$ and $U\cap\WFbd^{-1,\infty}(Pu)=\emptyset$.
Such a neighborhood exists since $q_0\nin\WFbd^{-1,\infty}(Pu)$
and $H_p\eta(\qt_0)
=\tau_0^2-|\zeta|^2>0$ for every $\qt_0\in\pih^{-1}(q_0)$.
Also let $U'$ be a subset of $U$ defined by replacing $\delta'$ by a
smaller $\delta''>0$, and let $\ep_0>0$ be such that for any
generalized broken bicharacteristic $\gamma$ with $\gamma(0)\in U'$,
$\gamma|_{[-\ep_0,\ep_0]}\in U$.
By Proposition~\ref{prop:normal-prop}, there is a sequence of points
$q_n\in\dot\Sigma$ such that $q_n\in\WFb^{1,\infty}(u)$,
$q_n\to q_0$ as $n\to\infty$, and
$\eta(q_n)<0$ for all $n$, so we may assume that
$q_n\in U'$ for all $n$.
By the inductive hypothesis, for each $n$, there exists a
generalized broken bicharcteristic
\begin{equation}
\gammat_n:(-\ep'_n,0]\to(\WFb^{1,\infty}(u)\setminus\WFbd^{-1,\infty}(Pu))
\cap\bigcup_{\cF_i\subsetneq \cF_j}T^*\cF_{j,\reg}
\end{equation}
with $\gammat_n(0)=q_n$. We now use
the argument of the first paragraph of the
proof (after the introductory remark about $s$)
with $V=(\WFb^{1,\infty}(u)\setminus\WFbd^{-1,\infty}(Pu))
\cap\bigcup_{\cF_i\subsetneq \cF_j}T^*\cF_{j,\reg}$, and $q=q_n$.
Thus, $\gammat_n\in\calR$, which is hence non-empty,
hence has a maximal element. We let
\begin{equation}
\gamma_n:(-\ep_n,0]\to(\WFb^{1,\infty}(u)\setminus\WFbd^{-1,\infty}(Pu))
\cap\bigcup_{\cF_i\subsetneq \cF_j}T^*\cF_{j,\reg}
\end{equation}
be a maximal element of $\calR$; it may happen that $-\ep_n=-\infty$.

We claim that $\ep_n\geq\ep_0$. For suppose that $\ep_n<\ep_0$.
By Corollary~\ref{cor:Lebeau-bichar-ext}, $\gamma_n$ extends
to a generalized broken bicharacteristic on $[-\ep_n,0]$, we continue
to denote this by $\gamma_n$. Since $\ep_n<\ep_0$, 
$\gamma_n$ is a generalized broken bicharacteristic
with image in $U$; indeed the closure of the image is still in $U$.
Taking into account that
$\eta$ is increasing on
generalized broken bicharacteristics in $U$
since $H_p\eta>0$ there, we conclude that
\begin{equation*}
-|\tau(\gamma_n(t))|^{-1}(x(\gamma_n(t))\cdot\xi(\gamma_n(t)))
=\eta(\gamma_n(t))\leq\eta(\gamma_n(0))<0
\end{equation*}
for $t\in[-\ep_n,0]$, hence $x(\gamma_n(t))\neq 0$. Thus, $\gamma_n(-\ep_n)
\in\cup_{\cF_i\subsetneq \cF_j}T^*\cF_{j,\reg}$. Moreover,
$\gamma_n(-\ep_n)\in\WFb^{1,\infty}(u)$ since $\WFb^{1,\infty}(u)$ is closed, and
$\gamma_n|_{(-\ep_n,0]}$ maps into it. Thus, by the inductive hypothesis,
there is a generalized broken bicharacteristic,
\begin{equation}
\gammat_n:(\alpha,-\ep_n]
\to(\WFb^{1,\infty}(u)\setminus\WFbd^{-1,\infty}(Pu))
\cap\bigcup_{\cF_i\subsetneq \cF_j}T^*\cF_{j,\reg}
\end{equation}
with $\alpha<-\ep_n$,
$\gammat_n(-\ep_n)=\gamma_n(-\ep_n)$. Hence, piecing together
$\gammat_n$ and $\gamma_n$ gives a generalized broken bicharacteristic
mapping into $(\WFb^{1,\infty}(u)\setminus\WFbd^{-1,\infty}(Pu))
\cap\bigcup_{\cF_i\subsetneq \cF_j}T^*\cF_{j,\reg}$ and extending $\gamma_n$,
which contradicts the maximal property of $\gamma_n$.
Thus, $\ep_n\geq \ep_0$ as claimed.

By Proposition~\ref{prop:Lebeau-compactness},
applied with $K=\WFb^{1,\infty}(u)$, there
is a subsequence of $\gamma_n|_{[-\ep_0,0]}$
converging uniformly to a generalized broken
bicharacteristic
\begin{equation*}
\gamma:[-\ep_0,0]\to\WFb^{1,\infty}(u).
\end{equation*}
In particular,
$\gamma(0)=q_0$ and $\gamma(t)\in\WFb^{1,\infty}(u)$ for all $t\in[-\ep_0,0]$,
providing the inductive step.

We now turn to $q_0\in\cG\cap T^*\cF_{i,\reg}$.
We repeat the argument of Melrose-Sj\"ostrand, as presented in
Lebeau's paper \cite[Proposition~VII.1]{Lebeau:Propagation}. We very
briefly outline the proof below; the detailed version follows Lebeau's
closely, with some changes in the notation.
Let $U\subset\cup_{\cF_i\subset \cF_j}T^*\cF_{j,\reg}\setminus\WFbd^{-1,\infty}(Pu)$
be a neighborhood of $q_0$, $U_0$ a smaller neighborhood,
as above. We take $\ep_0>0$
small. Suppose that
$0<\ep<\ep_0$, $q\in U_0$.
Let
\begin{equation}\begin{split}
\calR^1_{q,\ep}&=\{\text{generalized broken bicharacteristics}
\ \gamma:[-\ep,0]\to\WFb^{1,\infty}(u),\\
&\qquad\qquad\ \gamma(0)=q,\ \gamma(t)\nin\cG
\cap T^*\cF_{i,\reg}\Mfor t\in(-\ep,0]\},\\
\calR^2_{q,\ep}&=\{\text{generalized broken bicharacteristics}
\ \gamma:[-\ep',0]\to\WFb^{1,\infty}(u),\ \ep'\in(0,\ep),\\
&\qquad\qquad\ \gamma(0)=q,\ \gamma(t)\nin\cG
\cap T^*\cF_{i,\reg}\Mfor t\in(-\ep',0],\\
&\qquad\qquad\ \gamma(-\ep')
\in \cG
\cap T^*\cF_{i,\reg}\}.
\end{split}\end{equation}
Moreover, reflecting the inequalities in \eqref{eq:tgt-prop-est},
let
\begin{equation}
B(q,\ep)=\{q'\in\dot\Sigma:\ \max\{|\pi^e_i(q')-q|,
|x(q')|\}\leq\ep\}.
\end{equation}
Let $C_0>0$ be as in Proposition~\ref{prop:tgt-prop}.
For $q\in\cG\cap T^*\cF_{i,\reg}$, let
\begin{equation}
D(q,\ep)=B(\exp(-\ep H_p)(\pih^{-1}(q)),C_0\ep^2)\cap\WFb^{1,\infty}(u),
\end{equation}
and for $q\nin\cG\cap T^*\cF_{i,\reg}$, let
\begin{equation}\begin{split}
D(q,\ep)&=\{\gamma(-\ep):\ \gamma\in\calR^1_{q,\ep}\}\\
&\quad
\cup \{B(\exp(-(\ep-\ep') H_p)(\pih^{-1}(\gamma(\ep')),C_0(\ep-\ep')^2)
\cap\WFb^{1,\infty}(u):\ \gamma
\in\calR^2_{q,\ep}\}.
\end{split}\end{equation}
The reason for introducing $D(q,\ep)$ is that it is a good candidate
for the beginning point of a generalized broken bicharacteristic
segment in $\WFb^{1,\infty}(u)$, defined over an interval of length $\ep$, and
ending in $q$.

Indeed,
for $q\in\cG\cap T^*\cF_{i,\reg}\cap\WFb^{1,\infty}(u)$, we deduce from
Proposition~\ref{prop:tgt-prop} that $D(q,\ep)\neq\emptyset$.
For $q\in\WFb^{1,\infty}(u)\setminus(\cG\cap T^*\cF_{i,\reg})$,
by the inductive hypothesis, the previous part of the proof
concerning $\cH\cap T^*\cF_{i,\reg}$, and the
first two paragraphs (after the introductory remark about $s$) with
$V=\WFb^{1,\infty}(u)\setminus((\cG\cap T^*\cF_{i,\reg})\cup
\WFbd^{-1,\infty}(Pu)$, $q=q_0$, there is
a maximally extended generalized broken bicharacteristic $\gamma$
with image in $V$. By the argument of the second paragraph,
this is either defined on all of $[-\ep,0]$, or only
on $(-\ep',0]$ with $0<\ep'<\ep$, in which case $\gamma(-\ep')\in
\cG\cap T^*\cF_{i,\reg}$, hence again by
Proposition~\ref{prop:tgt-prop} we conclude that $D(q,\ep)\neq\emptyset$.
Thus, for all $q\in U\cap\WFb^{1,\infty}(u)$ we have deduced $D(q,\ep)\neq\emptyset$.

For each integer $N\geq 1$ now we define a sequence of $2^N+1$ points
$q_{j,N}$, $j\in\Nat$,
$0\leq j\leq 2^N$,
which will be used to construct points $\gamma(-j 2^{-N}\ep_0)$
on the desired generalized broken
bicharacteristic $\gamma:[-\ep_0,0]\to\WFb^{1,\infty}(u)$ through $q_0$.
Namely, let $\ep=2^{-N}\ep_0$, $q_{0,N}=q_0$, and choose $q_{j+1,N}
\in D(q_{j,N},\ep)$. Let $\calJ_N=\{-j 2^{-N}\ep_0:\ 0\leq j\leq 2^N\}
\subset[-\ep_0,0]$, $\calJ=\cup_{N=1}^\infty \calJ_N$. We write
$\gamma_N(t)=q_{j,N}$ for $t=-j 2^{-N}\ep_0$. For each $t\in \calJ$,
the sequence $\gamma_N(t)$ (defined for large $N$) stays in a compact set.
Hence there exists a subsequence $\gamma_{N_k}$ such that for all $t\in\calJ$,
$\gamma_{N_k}(t)$ converges to some $\gamma(t)$.

This defines $\gamma:[-\ep_0,0]\to\WFb^{1,\infty}(u)$ at elements of $\calJ$. One
can check exactly as in Lebeau's proof (which we have been following very
closely) that $\gamma$ extends to a continuous map defined on $[-\ep_0,0]$,
and that it is a generalized broken bicharacteristic. This completes
the inductive step for tangential points
$q_0\in\cG\cap T^*\cF_{i,\reg}$, hence the proof of the theorem.
\end{proof}

\bibliographystyle{plain}
\bibliography{sm}

\end{document}